
\catcode`\@=11


\message{Loading jyTeX fonts...}



\font\vptrm=cmr5 \font\vptmit=cmmi5 \font\vptsy=cmsy5 \font\vptbf=cmbx5

\skewchar\vptmit='177 \skewchar\vptsy='60 \fontdimen16
\vptsy=\the\fontdimen17 \vptsy

\def\vpt{\ifmmode\err@badsizechange\else
     \@mathfontinit
     \textfont0=\vptrm  \scriptfont0=\vptrm  \scriptscriptfont0=\vptrm
     \textfont1=\vptmit \scriptfont1=\vptmit \scriptscriptfont1=\vptmit
     \textfont2=\vptsy  \scriptfont2=\vptsy  \scriptscriptfont2=\vptsy
     \textfont3=\xptex  \scriptfont3=\xptex  \scriptscriptfont3=\xptex
     \textfont\bffam=\vptbf
     \scriptfont\bffam=\vptbf
     \scriptscriptfont\bffam=\vptbf
     \@fontstyleinit
     \def\rm{\vptrm\fam=\z@}%
     \def\bf{\vptbf\fam=\bffam}%
     \def\oldstyle{\vptmit\fam=\@ne}%
     \rm\fi}


\font\viptrm=cmr6 \font\viptmit=cmmi6 \font\viptsy=cmsy6
\font\viptbf=cmbx6

\skewchar\viptmit='177 \skewchar\viptsy='60 \fontdimen16
\viptsy=\the\fontdimen17 \viptsy

\def\vipt{\ifmmode\err@badsizechange\else
     \@mathfontinit
     \textfont0=\viptrm  \scriptfont0=\vptrm  \scriptscriptfont0=\vptrm
     \textfont1=\viptmit \scriptfont1=\vptmit \scriptscriptfont1=\vptmit
     \textfont2=\viptsy  \scriptfont2=\vptsy  \scriptscriptfont2=\vptsy
     \textfont3=\xptex   \scriptfont3=\xptex  \scriptscriptfont3=\xptex
     \textfont\bffam=\viptbf
     \scriptfont\bffam=\vptbf
     \scriptscriptfont\bffam=\vptbf
     \@fontstyleinit
     \def\rm{\viptrm\fam=\z@}%
     \def\bf{\viptbf\fam=\bffam}%
     \def\oldstyle{\viptmit\fam=\@ne}%
     \rm\fi}

\font\viiptrm=cmr7 \font\viiptmit=cmmi7 \font\viiptsy=cmsy7
\font\viiptit=cmti7 \font\viiptbf=cmbx7

\skewchar\viiptmit='177 \skewchar\viiptsy='60 \fontdimen16
\viiptsy=\the\fontdimen17 \viiptsy

\def\viipt{\ifmmode\err@badsizechange\else
     \@mathfontinit
     \textfont0=\viiptrm  \scriptfont0=\vptrm  \scriptscriptfont0=\vptrm
     \textfont1=\viiptmit \scriptfont1=\vptmit \scriptscriptfont1=\vptmit
     \textfont2=\viiptsy  \scriptfont2=\vptsy  \scriptscriptfont2=\vptsy
     \textfont3=\xptex    \scriptfont3=\xptex  \scriptscriptfont3=\xptex
     \textfont\itfam=\viiptit
     \scriptfont\itfam=\viiptit
     \scriptscriptfont\itfam=\viiptit
     \textfont\bffam=\viiptbf
     \scriptfont\bffam=\vptbf
     \scriptscriptfont\bffam=\vptbf
     \@fontstyleinit
     \def\rm{\viiptrm\fam=\z@}%
     \def\it{\viiptit\fam=\itfam}%
     \def\bf{\viiptbf\fam=\bffam}%
     \def\oldstyle{\viiptmit\fam=\@ne}%
     \rm\fi}


\font\viiiptrm=cmr8 \font\viiiptmit=cmmi8 \font\viiiptsy=cmsy8
\font\viiiptit=cmti8
\font\viiiptbf=cmbx8

\skewchar\viiiptmit='177 \skewchar\viiiptsy='60 \fontdimen16
\viiiptsy=\the\fontdimen17 \viiiptsy

\def\viiipt{\ifmmode\err@badsizechange\else
     \@mathfontinit
     \textfont0=\viiiptrm  \scriptfont0=\viptrm  \scriptscriptfont0=\vptrm
     \textfont1=\viiiptmit \scriptfont1=\viptmit \scriptscriptfont1=\vptmit
     \textfont2=\viiiptsy  \scriptfont2=\viptsy  \scriptscriptfont2=\vptsy
     \textfont3=\xptex     \scriptfont3=\xptex   \scriptscriptfont3=\xptex
     \textfont\itfam=\viiiptit
     \scriptfont\itfam=\viiptit
     \scriptscriptfont\itfam=\viiptit
     \textfont\bffam=\viiiptbf
     \scriptfont\bffam=\viptbf
     \scriptscriptfont\bffam=\vptbf
     \@fontstyleinit
     \def\rm{\viiiptrm\fam=\z@}%
     \def\it{\viiiptit\fam=\itfam}%
     \def\bf{\viiiptbf\fam=\bffam}%
     \def\oldstyle{\viiiptmit\fam=\@ne}%
     \rm\fi}


\def\getixpt{%
     \font\ixptrm=cmr9
     \font\ixptmit=cmmi9
     \font\ixptsy=cmsy9
     \font\ixptit=cmti9
     \font\ixptbf=cmbx9
     \skewchar\ixptmit='177 \skewchar\ixptsy='60
     \fontdimen16 \ixptsy=\the\fontdimen17 \ixptsy}

\def\ixpt{\ifmmode\err@badsizechange\else
     \@mathfontinit
     \textfont0=\ixptrm  \scriptfont0=\viiptrm  \scriptscriptfont0=\vptrm
     \textfont1=\ixptmit \scriptfont1=\viiptmit \scriptscriptfont1=\vptmit
     \textfont2=\ixptsy  \scriptfont2=\viiptsy  \scriptscriptfont2=\vptsy
     \textfont3=\xptex   \scriptfont3=\xptex    \scriptscriptfont3=\xptex
     \textfont\itfam=\ixptit
     \scriptfont\itfam=\viiptit
     \scriptscriptfont\itfam=\viiptit
     \textfont\bffam=\ixptbf
     \scriptfont\bffam=\viiptbf
     \scriptscriptfont\bffam=\vptbf
     \@fontstyleinit
     \def\rm{\ixptrm\fam=\z@}%
     \def\it{\ixptit\fam=\itfam}%
     \def\bf{\ixptbf\fam=\bffam}%
     \def\oldstyle{\ixptmit\fam=\@ne}%
     \rm\fi}


\font\xptrm=cmr10 \font\xptmit=cmmi10 \font\xptsy=cmsy10
\font\xptex=cmex10 \font\xptit=cmti10 \font\xptsl=cmsl10
\font\xptbf=cmbx10 \font\xpttt=cmtt10 \font\xptss=cmss10
\font\xptsc=cmcsc10 \font\xptbfs=cmb10 \font\xptbmit=cmmib10

\skewchar\xptmit='177 \skewchar\xptbmit='177 \skewchar\xptsy='60
\fontdimen16 \xptsy=\the\fontdimen17 \xptsy

\def\xpt{\ifmmode\err@badsizechange\else
     \@mathfontinit
     \textfont0=\xptrm  \scriptfont0=\viiptrm  \scriptscriptfont0=\vptrm
     \textfont1=\xptmit \scriptfont1=\viiptmit \scriptscriptfont1=\vptmit
     \textfont2=\xptsy  \scriptfont2=\viiptsy  \scriptscriptfont2=\vptsy
     \textfont3=\xptex  \scriptfont3=\xptex    \scriptscriptfont3=\xptex
     \textfont\itfam=\xptit
     \scriptfont\itfam=\viiptit
     \scriptscriptfont\itfam=\viiptit
     \textfont\bffam=\xptbf
     \scriptfont\bffam=\viiptbf
     \scriptscriptfont\bffam=\vptbf
     \textfont\bfsfam=\xptbfs
     \scriptfont\bfsfam=\viiptbf
     \scriptscriptfont\bfsfam=\vptbf
     \textfont\bmitfam=\xptbmit
     \scriptfont\bmitfam=\viiptmit
     \scriptscriptfont\bmitfam=\vptmit
     \@fontstyleinit
     \def\rm{\xptrm\fam=\z@}%
     \def\it{\xptit\fam=\itfam}%
     \def\sl{\xptsl}%
     \def\bf{\xptbf\fam=\bffam}%
     \def\tt{\xpttt}%
     \def\ss{\xptss}%
     \def\sc{\xptsc}%
     \def\bfs{\xptbfs\fam=\bfsfam}%
     \def\bmit{\fam=\bmitfam}%
     \def\oldstyle{\xptmit\fam=\@ne}%
     \rm\fi}


\def\getxipt{%
     \font\xiptrm=cmr10  scaled\magstephalf
     \font\xiptmit=cmmi10 scaled\magstephalf
     \font\xiptsy=cmsy10 scaled\magstephalf
     \font\xiptex=cmex10 scaled\magstephalf
     \font\xiptit=cmti10 scaled\magstephalf
     \font\xiptsl=cmsl10 scaled\magstephalf
     \font\xiptbf=cmbx10 scaled\magstephalf
     \font\xipttt=cmtt10 scaled\magstephalf
     \font\xiptss=cmss10 scaled\magstephalf
     \skewchar\xiptmit='177 \skewchar\xiptsy='60
     \fontdimen16 \xiptsy=\the\fontdimen17 \xiptsy}

\def\xipt{\ifmmode\err@badsizechange\else
     \@mathfontinit
     \textfont0=\xiptrm  \scriptfont0=\viiiptrm  \scriptscriptfont0=\viptrm
     \textfont1=\xiptmit \scriptfont1=\viiiptmit \scriptscriptfont1=\viptmit
     \textfont2=\xiptsy  \scriptfont2=\viiiptsy  \scriptscriptfont2=\viptsy
     \textfont3=\xiptex  \scriptfont3=\xptex     \scriptscriptfont3=\xptex
     \textfont\itfam=\xiptit
     \scriptfont\itfam=\viiiptit
     \scriptscriptfont\itfam=\viiptit
     \textfont\bffam=\xiptbf
     \scriptfont\bffam=\viiiptbf
     \scriptscriptfont\bffam=\viptbf
     \@fontstyleinit
     \def\rm{\xiptrm\fam=\z@}%
     \def\it{\xiptit\fam=\itfam}%
     \def\sl{\xiptsl}%
     \def\bf{\xiptbf\fam=\bffam}%
     \def\tt{\xipttt}%
     \def\ss{\xiptss}%
     \def\oldstyle{\xiptmit\fam=\@ne}%
     \rm\fi}


\font\xiiptrm=cmr12 \font\xiiptmit=cmmi12 \font\xiiptsy=cmsy10
scaled\magstep1 \font\xiiptex=cmex10  scaled\magstep1 \font\xiiptit=cmti12
\font\xiiptsl=cmsl12 \font\xiiptbf=cmbx12
\font\xiiptss=cmss12 \font\xiiptsc=cmcsc10 scaled\magstep1
\font\xiiptbfs=cmb10  scaled\magstep1 \font\xiiptbmit=cmmib10
scaled\magstep1

\skewchar\xiiptmit='177 \skewchar\xiiptbmit='177 \skewchar\xiiptsy='60
\fontdimen16 \xiiptsy=\the\fontdimen17 \xiiptsy

\def\xiipt{\ifmmode\err@badsizechange\else
     \@mathfontinit
     \textfont0=\xiiptrm  \scriptfont0=\viiiptrm  \scriptscriptfont0=\viptrm
     \textfont1=\xiiptmit \scriptfont1=\viiiptmit \scriptscriptfont1=\viptmit
     \textfont2=\xiiptsy  \scriptfont2=\viiiptsy  \scriptscriptfont2=\viptsy
     \textfont3=\xiiptex  \scriptfont3=\xptex     \scriptscriptfont3=\xptex
     \textfont\itfam=\xiiptit
     \scriptfont\itfam=\viiiptit
     \scriptscriptfont\itfam=\viiptit
     \textfont\bffam=\xiiptbf
     \scriptfont\bffam=\viiiptbf
     \scriptscriptfont\bffam=\viptbf
     \textfont\bfsfam=\xiiptbfs
     \scriptfont\bfsfam=\viiiptbf
     \scriptscriptfont\bfsfam=\viptbf
     \textfont\bmitfam=\xiiptbmit
     \scriptfont\bmitfam=\viiiptmit
     \scriptscriptfont\bmitfam=\viptmit
     \@fontstyleinit
     \def\rm{\xiiptrm\fam=\z@}%
     \def\it{\xiiptit\fam=\itfam}%
     \def\sl{\xiiptsl}%
     \def\bf{\xiiptbf\fam=\bffam}%
     \def\tt{\xiipttt}%
     \def\ss{\xiiptss}%
     \def\sc{\xiiptsc}%
     \def\bfs{\xiiptbfs\fam=\bfsfam}%
     \def\bmit{\fam=\bmitfam}%
     \def\oldstyle{\xiiptmit\fam=\@ne}%
     \rm\fi}


\def\getxiiipt{%
     \font\xiiiptrm=cmr12  scaled\magstephalf
     \font\xiiiptmit=cmmi12 scaled\magstephalf
     \font\xiiiptsy=cmsy9  scaled\magstep2
     \font\xiiiptit=cmti12 scaled\magstephalf
     \font\xiiiptsl=cmsl12 scaled\magstephalf
     \font\xiiiptbf=cmbx12 scaled\magstephalf
     \font\xiiipttt=cmtt12 scaled\magstephalf
     \font\xiiiptss=cmss12 scaled\magstephalf
     \skewchar\xiiiptmit='177 \skewchar\xiiiptsy='60
     \fontdimen16 \xiiiptsy=\the\fontdimen17 \xiiiptsy}

\def\xiiipt{\ifmmode\err@badsizechange\else
     \@mathfontinit
     \textfont0=\xiiiptrm  \scriptfont0=\xptrm  \scriptscriptfont0=\viiptrm
     \textfont1=\xiiiptmit \scriptfont1=\xptmit \scriptscriptfont1=\viiptmit
     \textfont2=\xiiiptsy  \scriptfont2=\xptsy  \scriptscriptfont2=\viiptsy
     \textfont3=\xivptex   \scriptfont3=\xptex  \scriptscriptfont3=\xptex
     \textfont\itfam=\xiiiptit
     \scriptfont\itfam=\xptit
     \scriptscriptfont\itfam=\viiptit
     \textfont\bffam=\xiiiptbf
     \scriptfont\bffam=\xptbf
     \scriptscriptfont\bffam=\viiptbf
     \@fontstyleinit
     \def\rm{\xiiiptrm\fam=\z@}%
     \def\it{\xiiiptit\fam=\itfam}%
     \def\sl{\xiiiptsl}%
     \def\bf{\xiiiptbf\fam=\bffam}%
     \def\tt{\xiiipttt}%
     \def\ss{\xiiiptss}%
     \def\oldstyle{\xiiiptmit\fam=\@ne}%
     \rm\fi}


\font\xivptrm=cmr12   scaled\magstep1 \font\xivptmit=cmmi12
scaled\magstep1 \font\xivptsy=cmsy10  scaled\magstep2
\font\xivptex=cmex10  scaled\magstep2 \font\xivptit=cmti12 scaled\magstep1
\font\xivptsl=cmsl12  scaled\magstep1 \font\xivptbf=cmbx12
scaled\magstep1
\font\xivptss=cmss12  scaled\magstep1 \font\xivptsc=cmcsc10
scaled\magstep2 \font\xivptbfs=cmb10  scaled\magstep2
\font\xivptbmit=cmmib10 scaled\magstep2

\skewchar\xivptmit='177 \skewchar\xivptbmit='177 \skewchar\xivptsy='60
\fontdimen16 \xivptsy=\the\fontdimen17 \xivptsy

\def\xivpt{\ifmmode\err@badsizechange\else
     \@mathfontinit
     \textfont0=\xivptrm  \scriptfont0=\xptrm  \scriptscriptfont0=\viiptrm
     \textfont1=\xivptmit \scriptfont1=\xptmit \scriptscriptfont1=\viiptmit
     \textfont2=\xivptsy  \scriptfont2=\xptsy  \scriptscriptfont2=\viiptsy
     \textfont3=\xivptex  \scriptfont3=\xptex  \scriptscriptfont3=\xptex
     \textfont\itfam=\xivptit
     \scriptfont\itfam=\xptit
     \scriptscriptfont\itfam=\viiptit
     \textfont\bffam=\xivptbf
     \scriptfont\bffam=\xptbf
     \scriptscriptfont\bffam=\viiptbf
     \textfont\bfsfam=\xivptbfs
     \scriptfont\bfsfam=\xptbfs
     \scriptscriptfont\bfsfam=\viiptbf
     \textfont\bmitfam=\xivptbmit
     \scriptfont\bmitfam=\xptbmit
     \scriptscriptfont\bmitfam=\viiptmit
     \@fontstyleinit
     \def\rm{\xivptrm\fam=\z@}%
     \def\it{\xivptit\fam=\itfam}%
     \def\sl{\xivptsl}%
     \def\bf{\xivptbf\fam=\bffam}%
     \def\tt{\xivpttt}%
     \def\ss{\xivptss}%
     \def\sc{\xivptsc}%
     \def\bfs{\xivptbfs\fam=\bfsfam}%
     \def\bmit{\fam=\bmitfam}%
     \def\oldstyle{\xivptmit\fam=\@ne}%
     \rm\fi}


\font\xviiptrm=cmr17 \font\xviiptmit=cmmi12 scaled\magstep2
\font\xviiptsy=cmsy10 scaled\magstep3 \font\xviiptex=cmex10
scaled\magstep3 \font\xviiptit=cmti12 scaled\magstep2
\font\xviiptbf=cmbx12 scaled\magstep2 \font\xviiptbfs=cmb10
scaled\magstep3

\skewchar\xviiptmit='177 \skewchar\xviiptsy='60 \fontdimen16
\xviiptsy=\the\fontdimen17 \xviiptsy

\def\xviipt{\ifmmode\err@badsizechange\else
     \@mathfontinit
     \textfont0=\xviiptrm  \scriptfont0=\xiiptrm  \scriptscriptfont0=\viiiptrm
     \textfont1=\xviiptmit \scriptfont1=\xiiptmit \scriptscriptfont1=\viiiptmit
     \textfont2=\xviiptsy  \scriptfont2=\xiiptsy  \scriptscriptfont2=\viiiptsy
     \textfont3=\xviiptex  \scriptfont3=\xiiptex  \scriptscriptfont3=\xptex
     \textfont\itfam=\xviiptit
     \scriptfont\itfam=\xiiptit
     \scriptscriptfont\itfam=\viiiptit
     \textfont\bffam=\xviiptbf
     \scriptfont\bffam=\xiiptbf
     \scriptscriptfont\bffam=\viiiptbf
     \textfont\bfsfam=\xviiptbfs
     \scriptfont\bfsfam=\xiiptbfs
     \scriptscriptfont\bfsfam=\viiiptbf
     \@fontstyleinit
     \def\rm{\xviiptrm\fam=\z@}%
     \def\it{\xviiptit\fam=\itfam}%
     \def\bf{\xviiptbf\fam=\bffam}%
     \def\bfs{\xviiptbfs\fam=\bfsfam}%
     \def\oldstyle{\xviiptmit\fam=\@ne}%
     \rm\fi}


\font\xxiptrm=cmr17  scaled\magstep1


\def\xxipt{\ifmmode\err@badsizechange\else
     \@mathfontinit
     \@fontstyleinit
     \def\rm{\xxiptrm\fam=\z@}%
     \rm\fi}


\font\xxvptrm=cmr17  scaled\magstep2


\def\xxvpt{\ifmmode\err@badsizechange\else
     \@mathfontinit
     \@fontstyleinit
     \def\rm{\xxvptrm\fam=\z@}%
     \rm\fi}




\message{Loading jyTeX macros...}

\message{modifications to plain.tex,}


\def\newcount{\alloc@0\count\countdef\insc@unt}
\def\newdimen{\alloc@1\dimen\dimendef\insc@unt}
\def\newskip{\alloc@2\skip\skipdef\insc@unt}
\def\newmuskip{\alloc@3\muskip\muskipdef\@cclvi}
\def\newbox{\alloc@4\box\chardef\insc@unt}
\def\newtoks{\alloc@5\toks\toksdef\@cclvi}
\def\newhelp#1#2{\newtoks#1\global#1\expandafter{\csname#2\endcsname}}
\def\newread{\alloc@6\read\chardef\sixt@@n}
\def\newwrite{\alloc@7\write\chardef\sixt@@n}
\def\newfam{\alloc@8\fam\chardef\sixt@@n}
\def\newinsert#1{\global\advance\insc@unt by\m@ne
     \ch@ck0\insc@unt\count
     \ch@ck1\insc@unt\dimen
     \ch@ck2\insc@unt\skip
     \ch@ck4\insc@unt\box
     \allocationnumber=\insc@unt
     \global\chardef#1=\allocationnumber
     \wlog{\string#1=\string\insert\the\allocationnumber}}
\def\newif#1{\count@\escapechar \escapechar\m@ne
     \expandafter\expandafter\expandafter
          \xdef\@if#1{true}{\let\noexpand#1=\noexpand\iftrue}%
     \expandafter\expandafter\expandafter
          \xdef\@if#1{false}{\let\noexpand#1=\noexpand\iffalse}%
     \global\@if#1{false}\escapechar=\count@}


\newlinechar=`\^^J
\overfullrule=0pt




\let\itfam=\undefined

\let\bffam=\undefined

\count18=3


\chardef\sharps="19


\mathchardef\alpha="710B \mathchardef\beta="710C
\mathchardef\gamma="710D \mathchardef\delta="710E
\mathchardef\epsilon="710F \mathchardef\zeta="7110
\mathchardef\eta="7111 \mathchardef\theta="7112 \mathchardef\iota="7113
\mathchardef\kappa="7114 \mathchardef\lambda="7115
\mathchardef\mu="7116 \mathchardef\nu="7117 \mathchardef\xi="7118
\mathchardef\pi="7119 \mathchardef\rho="711A \mathchardef\sigma="711B
\mathchardef\tau="711C \mathchardef\upsilon="711D
\mathchardef\phi="711E \mathchardef\chi="711F \mathchardef\psi="7120
\mathchardef\omega="7121 \mathchardef\varepsilon="7122
\mathchardef\vartheta="7123 \mathchardef\varpi="7124
\mathchardef\varrho="7125 \mathchardef\varsigma="7126
\mathchardef\varphi="7127 \mathchardef\imath="717B
\mathchardef\jmath="717C \mathchardef\ell="7160 \mathchardef\wp="717D
\mathchardef\partial="7140 \mathchardef\flat="715B
\mathchardef\natural="715C \mathchardef\sharp="715D



\def\angle{{\vbox{\ialign{$\m@th\scriptstyle##$\crcr
     \not\mathrel{\mkern14mu}\crcr
     \noalign{\nointerlineskip}
     \mkern2.5mu\leaders\hrule height.34\rp@\hfill\mkern2.5mu\crcr}}}}
\def\vdots{\vbox{\baselineskip4\rp@ \lineskiplimit\z@
     \kern6\rp@\hbox{.}\hbox{.}\hbox{.}}}
\def\ddots{\mathinner{\mkern1mu\raise7\rp@\vbox{\kern7\rp@\hbox{.}}\mkern2mu
     \raise4\rp@\hbox{.}\mkern2mu\raise\rp@\hbox{.}\mkern1mu}}
\def\overrightarrow#1{\vbox{\ialign{##\crcr
     \rightarrowfill\crcr
     \noalign{\kern-\rp@\nointerlineskip}
     $\hfil\displaystyle{#1}\hfil$\crcr}}}
\def\overleftarrow#1{\vbox{\ialign{##\crcr
     \leftarrowfill\crcr
     \noalign{\kern-\rp@\nointerlineskip}
     $\hfil\displaystyle{#1}\hfil$\crcr}}}
\def\overbrace#1{\mathop{\vbox{\ialign{##\crcr
     \noalign{\kern3\rp@}
     \downbracefill\crcr
     \noalign{\kern3\rp@\nointerlineskip}
     $\hfil\displaystyle{#1}\hfil$\crcr}}}\limits}
\def\underbrace#1{\mathop{\vtop{\ialign{##\crcr
     $\hfil\displaystyle{#1}\hfil$\crcr
     \noalign{\kern3\rp@\nointerlineskip}
     \upbracefill\crcr
     \noalign{\kern3\rp@}}}}\limits}
\def\big#1{{\hbox{$\left#1\vbox to8.5\rp@ {}\right.\n@space$}}}
\def\Big#1{{\hbox{$\left#1\vbox to11.5\rp@ {}\right.\n@space$}}}
\def\bigg#1{{\hbox{$\left#1\vbox to14.5\rp@ {}\right.\n@space$}}}
\def\Bigg#1{{\hbox{$\left#1\vbox to17.5\rp@ {}\right.\n@space$}}}
\def\@vereq#1#2{\lower.5\rp@\vbox{\baselineskip\z@skip\lineskip-.5\rp@
     \ialign{$\m@th#1\hfil##\hfil$\crcr#2\crcr=\crcr}}}
\def\rlh@#1{\vcenter{\hbox{\ooalign{\raise2\rp@
     \hbox{$#1\rightharpoonup$}\crcr
     $#1\leftharpoondown$}}}}
\def\bordermatrix#1{\begingroup\m@th
     \setbox\z@\vbox{%
          \def\cr{\crcr\noalign{\kern2\rp@\global\let\cr\endline}}%
          \ialign{$##$\hfil\kern2\rp@\kern\p@renwd
               &\thinspace\hfil$##$\hfil&&\quad\hfil$##$\hfil\crcr
               \omit\strut\hfil\crcr
               \noalign{\kern-\baselineskip}%
               #1\crcr\omit\strut\cr}}%
     \setbox\tw@\vbox{\unvcopy\z@\global\setbox\@ne\lastbox}%
     \setbox\tw@\hbox{\unhbox\@ne\unskip\global\setbox\@ne\lastbox}%
     \setbox\tw@\hbox{$\kern\wd\@ne\kern-\p@renwd\left(\kern-\wd\@ne
          \global\setbox\@ne\vbox{\box\@ne\kern2\rp@}%
          \vcenter{\kern-\ht\@ne\unvbox\z@\kern-\baselineskip}%
          \,\right)$}%
     \null\;\vbox{\kern\ht\@ne\box\tw@}\endgroup}
\def\endinsert{\egroup
     \if@mid\dimen@\ht\z@
          \advance\dimen@\dp\z@
          \advance\dimen@12\rp@
          \advance\dimen@\pagetotal
          \ifdim\dimen@>\pagegoal\@midfalse\p@gefalse\fi
     \fi
     \if@mid\bigskip\box\z@
          \bigbreak
     \else\insert\topins{\penalty100 \splittopskip\z@skip
               \splitmaxdepth\maxdimen\floatingpenalty\z@
               \ifp@ge\dimen@\dp\z@
                    \vbox to\vsize{\unvbox\z@\kern-\dimen@}%
               \else\box\z@\nobreak\bigskip
               \fi}%
     \fi
     \endgroup}


\def\cases#1{\left\{\,\vcenter{\m@th
     \ialign{$##\hfil$&\quad##\hfil\crcr#1\crcr}}\right.}
\def\matrix#1{\null\,\vcenter{\m@th
     \ialign{\hfil$##$\hfil&&\quad\hfil$##$\hfil\crcr
          \mathstrut\crcr
          \noalign{\kern-\baselineskip}
          #1\crcr
          \mathstrut\crcr
          \noalign{\kern-\baselineskip}}}\,}


\newif\ifraggedbottom

\def\raggedbottom{\ifraggedbottom\else
     \advance\topskip by\z@ plus60pt \raggedbottomtrue\fi}%
\def\normalbottom{\ifraggedbottom
     \advance\topskip by\z@ plus-60pt \raggedbottomfalse\fi}

\message{hacks,}


\toksdef\toks@i=1 \toksdef\toks@ii=2


\def\TeX{T\kern-.1667em \lower.5ex \hbox{E}\kern-.125em X\null}
\def\jyTeX{{\leavevmode
     \raise.587ex \hbox{\it\j}\kern-.1em \lower.048ex \hbox{\it y}\kern-.12em
     \TeX}}

\let\then=\iftrue
\def\ifnoarg#1\then{\def\hack@{#1}\ifx\hack@\empty}
\def\ifundefined#1\then{%
     \expandafter\ifx\csname\expandafter\blank\string#1\endcsname\relax}
\def\useif#1\then{\csname#1\endcsname}
\def\usename#1{\csname#1\endcsname}
\def\useafter#1#2{\expandafter#1\csname#2\endcsname}

\long\def\loop#1\repeat{\def\@iterate{#1\expandafter\@iterate\fi}\@iterate
     \let\@iterate=\relax}

\let\TeXend=\end
\def\begin#1{\begingroup\def\@@blockname{#1}\usename{begin#1}}
\def\end#1{\usename{end#1}\def\hack@{#1}%
     \ifx\@@blockname\hack@
          \endgroup
     \else\err@badgroup\hack@\@@blockname
     \fi}
\def\@@blockname{}

\def\defaultoption[#1]#2{%
     \def\hack@{\ifx\hack@ii[\toks@={#2}\else\toks@={#2[#1]}\fi\the\toks@}%
     \futurelet\hack@ii\hack@}

\def\markup#1{\let\@@marksf=\empty
     \ifhmode\edef\@@marksf{\spacefactor=\the\spacefactor\relax}\/\fi
     ${}^{\hbox{\subscriptfonts#1}}$\@@marksf}


\newtoks\shortyear
\newtoks\militaryhour
\newtoks\standardhour
\newtoks\minute
\newtoks\amorpm

\def\settime{\count@=\time\divide\count@ by60
     \militaryhour=\expandafter{\number\count@}%
     {\multiply\count@ by-60 \advance\count@ by\time
          \xdef\hack@{\ifnum\count@<10 0\fi\number\count@}}%
     \minute=\expandafter{\hack@}%
     \ifnum\count@<12
          \amorpm={am}
     \else\amorpm={pm}
          \ifnum\count@>12 \advance\count@ by-12 \fi
     \fi
     \standardhour=\expandafter{\number\count@}%
     \def\hack@19##1##2{\shortyear={##1##2}}%
          \expandafter\hack@\the\year}

\def\monthword#1{%
     \ifcase#1
          $\bullet$\err@badcountervalue{monthword}%
          \or January\or February\or March\or April\or May\or June%
          \or July\or August\or September\or October\or November\or December%
     \else$\bullet$\err@badcountervalue{monthword}%
     \fi}

\def\monthabbr#1{%
     \ifcase#1
          $\bullet$\err@badcountervalue{monthabbr}%
          \or Jan\or Feb\or Mar\or Apr\or May\or Jun%
          \or Jul\or Aug\or Sep\or Oct\or Nov\or Dec%
     \else$\bullet$\err@badcountervalue{monthabbr}%
     \fi}

\def\militarytime{\the\militaryhour:\the\minute}
\def\standardtime{\the\standardhour:\the\minute}


\def\@setnumstyle#1#2{\expandafter\global\expandafter\expandafter
     \expandafter\let\expandafter\expandafter
     \csname @\expandafter\blank\string#1style\endcsname
     \csname#2\endcsname}
\def\numstyle#1{\usename{@\expandafter\blank\string#1style}#1}
\def\ifblank#1\then{\useafter\ifx{@\expandafter\blank\string#1}\blank}

\def\blank#1{}

\def\Roman#1{\expandafter\uppercase\expandafter{\romannumeral#1}}
\def\alphabetic#1{%
     \ifcase#1
          $\bullet$\err@badcountervalue{alphabetic}%
          \or a\or b\or c\or d\or e\or f\or g\or h\or i\or j\or k\or l\or m%
          \or n\or o\or p\or q\or r\or s\or t\or u\or v\or w\or x\or y\or z%
     \else$\bullet$\err@badcountervalue{alphabetic}%
     \fi}
\def\Alphabetic#1{\expandafter\uppercase\expandafter{\alphabetic{#1}}}
\def\symbols#1{%
     \ifcase#1
          $\bullet$\err@badcountervalue{symbols}%
          \or*\or\dag\or\ddag\or\S\or$\|$%
          \or**\or\dag\dag\or\ddag\ddag\or\S\S\or$\|\|$%
     \else$\bullet$\err@badcountervalue{symbols}%
     \fi}


\catcode`\^^?=13 \def^^?{\relax}

\def\trimleading#1\to#2{\edef#2{#1}%
     \expandafter\@trimleading\expandafter#2#2^^?^^?}
\def\@trimleading#1#2#3^^?{\ifx#2^^?\def#1{}\else\def#1{#2#3}\fi}

\def\trimtrailing#1\to#2{\edef#2{#1}%
     \expandafter\@trimtrailing\expandafter#2#2^^? ^^?\relax}
\def\@trimtrailing#1#2 ^^?#3{\ifx#3\relax\toks@={}%
     \else\def#1{#2}\toks@={\trimtrailing#1\to#1}\fi
     \the\toks@}

\def\trim#1\to#2{\trimleading#1\to#2\trimtrailing#2\to#2}

\catcode`\^^?=15


\long\def\additemL#1\to#2{\toks@={\^^\{#1}}\toks@ii=\expandafter{#2}%
     \xdef#2{\the\toks@\the\toks@ii}}

\long\def\additemR#1\to#2{\toks@={\^^\{#1}}\toks@ii=\expandafter{#2}%
     \xdef#2{\the\toks@ii\the\toks@}}

\def\getitemL#1\to#2{\expandafter\@getitemL#1\hack@#1#2}
\def\@getitemL\^^\#1#2\hack@#3#4{\def#4{#1}\def#3{#2}}

\message{font macros,}


\newdimen\rp@
\newcount\@@sizeindex \@@sizeindex=0
\newcount\@@factori
\newcount\@@factorii
\newcount\@@factoriii
\newcount\@@factoriv

\countdef\maxfam=18
\newfam\itfam
\newfam\bffam
\newfam\bfsfam
\newfam\bmitfam

\def\@mathfontinit{\count@=4
     \loop\textfont\count@=\nullfont
          \scriptfont\count@=\nullfont
          \scriptscriptfont\count@=\nullfont
          \ifnum\count@<\maxfam\advance\count@ by\@ne
     \repeat}

\def\@fontstyleinit{%
     \def\it{\err@fontnotavailable\it}%
     \def\bf{\err@fontnotavailable\bf}%
     \def\bfs{\err@bfstobf}%
     \def\bmit{\err@fontnotavailable\bmit}%
     \def\sc{\err@fontnotavailable\sc}%
     \def\sl{\err@sltoit}%
     \def\ss{\err@fontnotavailable\ss}%
     \def\tt{\err@fontnotavailable\tt}}

\def\@parameterinit#1{\rm\rp@=.1em \@getscaling{#1}%
     \let\^^\=\@doscaling\scalingskipslist
     \setbox\strutbox=\hbox{\vrule
          height.708\baselineskip depth.292\baselineskip width\z@}}

\def\@getfactor#1#2#3#4{\@@factori=#1 \@@factorii=#2
     \@@factoriii=#3 \@@factoriv=#4}

\def\@getscaling#1{\count@=#1 \advance\count@ by-\@@sizeindex\@@sizeindex=#1
     \ifnum\count@<0
          \let\@mulordiv=\divide
          \let\@divormul=\multiply
          \multiply\count@ by\m@ne
     \else\let\@mulordiv=\multiply
          \let\@divormul=\divide
     \fi
     \edef\@@scratcha{\ifcase\count@                {1}{1}{1}{1}\or
          {1}{7}{23}{3}\or     {2}{5}{3}{1}\or      {9}{89}{13}{1}\or
          {6}{25}{6}{1}\or     {8}{71}{14}{1}\or    {6}{25}{36}{5}\or
          {1}{7}{53}{4}\or     {12}{125}{108}{5}\or {3}{14}{53}{5}\or
          {6}{41}{17}{1}\or    {13}{31}{13}{2}\or   {9}{107}{71}{2}\or
          {11}{139}{124}{3}\or {1}{6}{43}{2}\or     {10}{107}{42}{1}\or
          {1}{5}{43}{2}\or     {5}{69}{65}{1}\or    {11}{97}{91}{2}\fi}%
     \expandafter\@getfactor\@@scratcha}

\def\@doscaling#1{\@mulordiv#1by\@@factori\@divormul#1by\@@factorii
     \@mulordiv#1by\@@factoriii\@divormul#1by\@@factoriv}


\newskip\headskip
\newskip\footskip

\def\typesize=#1pt{\count@=#1 \advance\count@ by-10
     \ifcase\count@
          \@setsizex\or\err@badtypesize\or
          \@setsizexii\or\err@badtypesize\or
          \@setsizexiv
     \else\err@badtypesize
     \fi}

\def\@setsizex{\getixpt
     \def\subsubscriptfonts{\vpt}%
          \def\subsubscriptsize{\vpt\@parameterinit{-8}}%
     \def\subscriptfonts{\viipt}\def\subscriptsize{\viipt\@parameterinit{-4}}%
     \def\footnotefonts{\viiipt}\def\footnotesize{\viiipt\@parameterinit{-2}}%
     \def\smallfonts{\ixpt}\def\smallsize{\ixpt\@parameterinit{-1}}%
     \def\normalfonts{\xpt}\def\normalsize{\xpt\@parameterinit{0}}%
     \def\bigfonts{\xiipt}\def\bigsize{\xiipt\@parameterinit{2}}%
     \def\Bigfonts{\xivpt}\def\Bigsize{\xivpt\@parameterinit{4}}%
     \def\biggfonts{\xviipt}\def\biggsize{\xviipt\@parameterinit{6}}%
     \def\Biggfonts{\xxipt}\def\Biggsize{\xxipt\@parameterinit{8}}%
     \def\tinyfonts{\vpt}\def\tinysize{\vpt\@parameterinit{-8}}%
     \def\HUGEFONTS{\xxvpt}\def\HUGESIZE{\xxvpt\@parameterinit{10}}%
     \normalsize\fixedskipslist}

\def\@setsizexii{\getxipt
     \def\subsubscriptfonts{\vipt}%
          \def\subsubscriptsize{\vipt\@parameterinit{-6}}%
     \def\subscriptfonts{\viiipt}%
          \def\subscriptsize{\viiipt\@parameterinit{-2}}%
     \def\footnotefonts{\xpt}\def\footnotesize{\xpt\@parameterinit{0}}%
     \def\smallfonts{\xipt}\def\smallsize{\xipt\@parameterinit{1}}%
     \def\normalfonts{\xiipt}\def\normalsize{\xiipt\@parameterinit{2}}%
     \def\bigfonts{\xivpt}\def\bigsize{\xivpt\@parameterinit{4}}%
     \def\Bigfonts{\xviipt}\def\Bigsize{\xviipt\@parameterinit{6}}%
     \def\biggfonts{\xxipt}\def\biggsize{\xxipt\@parameterinit{8}}%
     \def\Biggfonts{\xxvpt}\def\Biggsize{\xxvpt\@parameterinit{10}}%
     \def\tinyfonts{\vpt}\def\tinysize{\vpt\@parameterinit{-8}}%
     \def\HUGEFONTS{\xxvpt}\def\HUGESIZE{\xxvpt\@parameterinit{10}}%
     \normalsize\fixedskipslist}

\def\@setsizexiv{\getxiiipt
     \def\subsubscriptfonts{\viipt}%
          \def\subsubscriptsize{\viipt\@parameterinit{-4}}%
     \def\subscriptfonts{\xpt}\def\subscriptsize{\xpt\@parameterinit{0}}%
     \def\footnotefonts{\xiipt}\def\footnotesize{\xiipt\@parameterinit{2}}%
     \def\smallfonts{\xiiipt}\def\smallsize{\xiiipt\@parameterinit{3}}%
     \def\normalfonts{\xivpt}\def\normalsize{\xivpt\@parameterinit{4}}%
     \def\bigfonts{\xviipt}\def\bigsize{\xviipt\@parameterinit{6}}%
     \def\Bigfonts{\xxipt}\def\Bigsize{\xxipt\@parameterinit{8}}%
     \def\biggfonts{\xxvpt}\def\biggsize{\xxvpt\@parameterinit{10}}%
     \def\Biggfonts{\err@sizetoolarge\Biggfonts\HUGEFONTS}%
          \def\Biggsize{\err@sizetoolarge\Biggsize\HUGESIZE}%
     \def\tinyfonts{\vpt}\def\tinysize{\vpt\@parameterinit{-8}}%
     \def\HUGEFONTS{\xxvpt}\def\HUGESIZE{\xxvpt\@parameterinit{10}}%
     \normalsize\fixedskipslist}

\def\subsubscriptfonts{\vpt} \def\subsubscriptsize{\vpt\@parameterinit{-8}}
\def\subscriptfonts{\viipt}  \def\subscriptsize{\viipt\@parameterinit{-4}}
\def\footnotefonts{\viiipt}  \def\footnotesize{\viiipt\@parameterinit{-2}}
\def\smallfonts{\err@sizenotavailable\smallfonts}
                             \def\smallsize{\ixpt\@parameterinit{-1}}
\def\normalfonts{\xpt}       \def\normalsize{\xpt\@parameterinit{0}}
\def\bigfonts{\xiipt}        \def\bigsize{\xiipt\@parameterinit{2}}
\def\Bigfonts{\xivpt}        \def\Bigsize{\xivpt\@parameterinit{4}}
\def\biggfonts{\xviipt}      \def\biggsize{\xviipt\@parameterinit{6}}
\def\Biggfonts{\xxipt}       \def\Biggsize{\xxipt\@parameterinit{8}}
\def\tinyfonts{\vpt}         \def\tinysize{\vpt\@parameterinit{-8}}
\def\HUGEFONTS{\xxvpt}       \def\HUGESIZE{\xxvpt\@parameterinit{10}}

\message{document layout,}


\newtoks\everyoutput \everyoutput={}
\newdimen\depthofpage
\newcount\pagenum \pagenum=0

\newdimen\oddtopmargin  \newdimen\eventopmargin
\newdimen\oddleftmargin \newdimen\evenleftmargin
\newtoks\oddhead        \newtoks\evenhead
\newtoks\oddfoot        \newtoks\evenfoot

\def\topmargin{\afterassignment\@seteventop\oddtopmargin}
\def\leftmargin{\afterassignment\@setevenleft\oddleftmargin}
\def\head{\afterassignment\@setevenhead\oddhead}
\def\foot{\afterassignment\@setevenfoot\oddfoot}

\def\@seteventop{\eventopmargin=\oddtopmargin}
\def\@setevenleft{\evenleftmargin=\oddleftmargin}
\def\@setevenhead{\evenhead=\oddhead}
\def\@setevenfoot{\evenfoot=\oddfoot}

\def\pagenumstyle#1{\@setnumstyle\pagenum{#1}}

\newif\ifdraft
\def\draft{\drafttrue\leftmargin=.5in \overfullrule=5pt }

\def\outputstyle#1{\global\expandafter\let\expandafter
          \@outputstyle\csname#1output\endcsname
     \usename{#1setup}}

\output={\@outputstyle}

\def\normaloutput{\the\everyoutput
     \global\advance\pagenum by\@ne
     \ifodd\pagenum
          \voffset=\oddtopmargin \hoffset=\oddleftmargin
     \else\voffset=\eventopmargin \hoffset=\evenleftmargin
     \fi
     \advance\voffset by-1in  \advance\hoffset by-1in
     \count0=\pagenum
     \expandafter\shipout\pagebox
     \ifnum\outputpenalty>-\@MM\else\dosupereject\fi}

\newdimen\fullhsize
\newbox\leftpage
\newcount\leftpagenum
\newcount\outputpagenum \outputpagenum=0
\let\leftorright=L

\def\twoupoutput{\the\everyoutput
     \global\advance\pagenum by\@ne
     \if L\leftorright
          \global\setbox\leftpage=\leftline{\pagebox}%
          \global\leftpagenum=\pagenum
          \global\let\leftorright=R%
     \else\global\advance\outputpagenum by\@ne
          \ifodd\outputpagenum
               \voffset=\oddtopmargin \hoffset=\oddleftmargin
          \else\voffset=\eventopmargin \hoffset=\evenleftmargin
          \fi
          \advance\voffset by-1in  \advance\hoffset by-1in
          \count0=\leftpagenum \count1=\pagenum
          \shipout\vbox{\hbox to\fullhsize
               {\box\leftpage\hfil\leftline{\pagebox}}}%
          \global\let\leftorright=L%
     \fi
     \ifnum\outputpenalty>-\@MM
     \else\dosupereject
          \if R\leftorright
               \globaldefs=\@ne\head={\hfil}\foot={\hfil}\globaldefs=\z@
               \null\newpage
          \fi
     \fi}

\def\pagebox{\vbox{\makeheadline\pagebody\makefootline}}

\def\makeheadline{%
     \vbox to\z@{\baselinestretch=\@m
          \vskip\topskip\vskip-.708\baselineskip\vskip-\headskip
          \line{\vbox to\ht\strutbox{}%
               \ifodd\pagenum\the\oddhead\else\the\evenhead\fi}%
          \vss}%
     \nointerlineskip}

\def\pagebody{\vbox to\vsize{%
     \boxmaxdepth\maxdepth
     \ifvoid\topins\else\unvbox\topins\fi
     \depthofpage=\dp255
     \unvbox255
     \ifraggedbottom\kern-\depthofpage\vfil\fi
     \ifvoid\footins
     \else\vskip\skip\footins
          \footnoterule
          \unvbox\footins
          \vskip-\footnoteskip
     \fi}}

\def\makefootline{\baselineskip=\footskip
     \line{\ifodd\pagenum\the\oddfoot\else\the\evenfoot\fi}}


\newskip\abovechapterskip
\newskip\belowchapterskip
\newskip\abovesectionskip
\newskip\belowsectionskip
\newskip\abovesubsectionskip
\newskip\belowsubsectionskip

\def\chapterstyle#1{\global\expandafter\let\expandafter\@chapterstyle
     \csname#1text\endcsname}
\def\sectionstyle#1{\global\expandafter\let\expandafter\@sectionstyle
     \csname#1text\endcsname}
\def\subsectionstyle#1{\global\expandafter\let\expandafter\@subsectionstyle
     \csname#1text\endcsname}

\def\chapter#1{%
     \ifdim\lastskip=17sp \else\chapterbreak\vskip\abovechapterskip\fi
     \@chapterstyle{\ifblank\chapternumstyle\then
          \else\newchapternum=\next\chapternumformat\ \fi#1}%
     \nobreak\vskip\belowchapterskip\vskip17sp }

\def\section#1{%
     \ifdim\lastskip=17sp \else\sectionbreak\vskip\abovesectionskip\fi
     \@sectionstyle{\ifblank\sectionnumstyle\then
          \else\newsectionnum=\next\sectionnumformat\ \fi#1}%
     \nobreak\vskip\belowsectionskip\vskip17sp }

\def\subsection#1{%
     \ifdim\lastskip=17sp \else\subsectionbreak\vskip\abovesubsectionskip\fi
     \@subsectionstyle{\ifblank\subsectionnumstyle\then
          \else\newsubsectionnum=\next\subsectionnumformat\ \fi#1}%
     \nobreak\vskip\belowsubsectionskip\vskip17sp }


\let\TeXunderline=\underline
\let\TeXoverline=\overline
\def\underline#1{\relax\ifmmode\TeXunderline{#1}\else
     $\TeXunderline{\hbox{#1}}$\fi}
\def\overline#1{\relax\ifmmode\TeXoverline{#1}\else
     $\TeXoverline{\hbox{#1}}$\fi}

\def\baselinestretch{\afterassignment\@baselinestretch\count@}
\def\@baselinestretch{\baselineskip=\normalbaselineskip
     \divide\baselineskip by\@m\baselineskip=\count@\baselineskip
     \setbox\strutbox=\hbox{\vrule
          height.708\baselineskip depth.292\baselineskip width\z@}%
     \bigskipamount=\the\baselineskip
          plus.25\baselineskip minus.25\baselineskip
     \medskipamount=.5\baselineskip
          plus.125\baselineskip minus.125\baselineskip
     \smallskipamount=.25\baselineskip
          plus.0625\baselineskip minus.0625\baselineskip}

\def\\{\ifhmode\ifnum\lastpenalty=-\@M\else\hfil\penalty-\@M\fi\fi
     \ignorespaces}
\def\newpage{\vfil\break}

\def\lefttext#1{\par{\@text\leftskip=\z@\rightskip=\centering
     \noindent#1\par}}
\def\righttext#1{\par{\@text\leftskip=\centering\rightskip=\z@
     \noindent#1\par}}
\def\centertext#1{\par{\@text\leftskip=\centering\rightskip=\centering
     \noindent#1\par}}
\def\@text{\parindent=\z@ \parfillskip=\z@ \everypar={}%
     \spaceskip=.3333em \xspaceskip=.5em
     \def\\{\ifhmode\ifnum\lastpenalty=-\@M\else\penalty-\@M\fi\fi
          \ignorespaces}}

\def\beginleft{\par\@text\leftskip=\z@ \rightskip=\centering}
     
\def\beginright{\par\@text\leftskip=\centering\rightskip=\z@ }
     
\def\begincenter{\par\@text\leftskip=\centering\rightskip=\centering}

\def\beginnarrow{\defaultoption[\parindent]\@beginnarrow}
\def\@beginnarrow[#1]{\par\advance\leftskip by#1\advance\rightskip by#1}

\begingroup
\catcode`\[=1 \catcode`\{=11 \gdef\beginignore[\endgroup\bgroup
     \catcode`\e=0 \catcode`\\=12 \catcode`\{=11 \catcode`\f=12 \let\or=\relax
     \let\nd{ignor=\fi \let\}=\egroup
     \iffalse}
\endgroup

\long\def\marginnote#1{\leavevmode
     \edef\@marginsf{\spacefactor=\the\spacefactor\relax}%
     \ifdraft\strut\vadjust{%
          \hbox to\z@{\hskip\hsize\hskip.1in
               \vbox to\z@{\vskip-\dp\strutbox
                    \marginnoteformat
                    \vskip-\ht\strutbox
                    \noindent\strut#1\par
                    \vss}%
               \hss}}%
     \fi
     \@marginsf}


\newtoks\everybye \everybye={\par\vfil}
\outer\def\bye{\the\everybye
     \footnotecheck
     \prelabelcheck
     \streamcheck
     \supereject
     \TeXend}

\message{footnotes,}

\newcount\footnotenum \footnotenum=0
\newskip\footnoteskip
\let\@footnotelist=\empty

\def\footnotenumstyle#1{\@setnumstyle\footnotenum{#1}%
     \useafter\ifx{@footnotenumstyle}\symbols
          \global\let\@footup=\empty
     \else\global\let\@footup=\markup
     \fi}

\def\footnote{\footnotecheck\defaultoption[]\@footnote}
\def\@footnote[#1]{\@footnotemark[#1]\@footnotetext}

\def\footnotemark{\defaultoption[]\@footnotemark}
\def\@footnotemark[#1]{\let\@footsf=\empty
     \ifhmode\edef\@footsf{\spacefactor=\the\spacefactor\relax}\/\fi
     \ifnoarg#1\then
          \global\advance\footnotenum by\@ne
          \@footup{\footnotenumformat}%
          \edef\@@foota{\footnotenum=\the\footnotenum\relax}%
          \expandafter\additemR\expandafter\@footup\expandafter
               {\@@foota\footnotenumformat}\to\@footnotelist
          \global\let\@footnotelist=\@footnotelist
     \else\markup{#1}%
          \additemR\markup{#1}\to\@footnotelist
          \global\let\@footnotelist=\@footnotelist
     \fi
     \@footsf}

\def\footnotetext{%
     \ifx\@footnotelist\empty\err@extrafootnotetext\else\@footnotetext\fi}
\def\@footnotetext{%
     \getitemL\@footnotelist\to\@@foota
     \global\let\@footnotelist=\@footnotelist
     \insert\footins\bgroup
     \footnoteformat
     \splittopskip=\ht\strutbox\splitmaxdepth=\dp\strutbox
     \interlinepenalty=\interfootnotelinepenalty\floatingpenalty=\@MM
     \noindent\llap{\@@foota}\strut
     \bgroup\aftergroup\@footnoteend
     \let\@@scratcha=}
\def\@footnoteend{\strut\par\vskip\footnoteskip\egroup}

\def\footnoterule{\normalfonts
     \kern-.3em \hrule width2in height.04em \kern .26em }

\def\footnotecheck{%
     \ifx\@footnotelist\empty
     \else\err@extrafootnotemark
          \global\let\@footnotelist=\empty
     \fi}

\message{labels,}

\let\@@labeldef=\xdef
\newif\if@labelfile
\newwrite\@labelfile
\let\@prelabellist=\empty

\def\label#1#2{\trim#1\to\@@labarg\edef\@@labtext{#2}%
     \edef\@@labname{lab@\@@labarg}%
     \useafter\ifundefined\@@labname\then\else\@yeslab\fi
     \useafter\@@labeldef\@@labname{#2}%
     \ifstreaming
          \expandafter\toks@\expandafter\expandafter\expandafter
               {\csname\@@labname\endcsname}%
          \immediate\write\streamout{\noexpand\label{\@@labarg}{\the\toks@}}%
     \fi}
\def\@yeslab{%
     \useafter\ifundefined{if\@@labname}\then
          \err@labelredef\@@labarg
     \else\useif{if\@@labname}\then
               \err@labelredef\@@labarg
          \else\global\usename{\@@labname true}%
               \useafter\ifundefined{pre\@@labname}\then
               \else\useafter\ifx{pre\@@labname}\@@labtext
                    \else\err@badlabelmatch\@@labarg
                    \fi
               \fi
               \if@labelfile
               \else\global\@labelfiletrue
                    \immediate\write\sixt@@n{--> Creating file \jobname.lab}%
                    \immediate\openout\@labelfile=\jobname.lab
               \fi
               \immediate\write\@labelfile
                    {\noexpand\prelabel{\@@labarg}{\@@labtext}}%
          \fi
     \fi}

\def\putlab#1{\trim#1\to\@@labarg\edef\@@labname{lab@\@@labarg}%
     \useafter\ifundefined\@@labname\then\@nolab\else\usename\@@labname\fi}
\def\@nolab{%
     \useafter\ifundefined{pre\@@labname}\then
          \undefinedlabelformat
          \err@needlabel\@@labarg
          \useafter\xdef\@@labname{\undefinedlabelformat}%
     \else\usename{pre\@@labname}%
          \useafter\xdef\@@labname{\usename{pre\@@labname}}%
     \fi
     \useafter\newif{if\@@labname}%
     \expandafter\additemR\@@labarg\to\@prelabellist}

\def\prelabel#1{\useafter\gdef{prelab@#1}}

\def\ifundefinedlabel#1\then{%
     \expandafter\ifx\csname lab@#1\endcsname\relax}
\def\useiflab#1\then{\csname iflab@#1\endcsname}

\def\prelabelcheck{{%
     \def\^^\##1{\useiflab{##1}\then\else\err@undefinedlabel{##1}\fi}%
     \@prelabellist}}

\message{equation numbering,}

\newcount\chapternum
\newcount\sectionnum
\newcount\subsectionnum
\newcount\equationnum
\newcount\subequationnum
\newcount\figurenum
\newcount\subfigurenum
\newcount\tablenum
\newcount\subtablenum

\newif\if@subeqncount
\newif\if@subfigcount
\newif\if@subtblcount

\def\newchapternum{\newsectionnum=\z@\@resetnum\chapternum}
\def\newsectionnum{\newsubsectionnum=\z@\@resetnum\sectionnum}
\def\newsubsectionnum{\newequationnum=\z@\newfigurenum=\z@\newtablenum=\z@
     \@resetnum\subsectionnum}
\def\newequationnum{\newsubequationnum=\z@\@resetnum\equationnum}
\def\newsubequationnum{\@resetnum\subequationnum}
\def\newfigurenum{\newsubfigurenum=\z@\@resetnum\figurenum}
\def\newsubfigurenum{\@resetnum\subfigurenum}
\def\newtablenum{\newsubtablenum=\z@\@resetnum\tablenum}
\def\newsubtablenum{\@resetnum\subtablenum}

\def\@resetnum#1{\global\advance#1by1 \edef\next{\the#1\relax}\global#1}

\newchapternum=0

\def\chapternumstyle#1{\@setnumstyle\chapternum{#1}}
\def\sectionnumstyle#1{\@setnumstyle\sectionnum{#1}}
\def\subsectionnumstyle#1{\@setnumstyle\subsectionnum{#1}}
\def\equationnumstyle#1{\@setnumstyle\equationnum{#1}}
\def\subequationnumstyle#1{\@setnumstyle\subequationnum{#1}%
     \ifblank\subequationnumstyle\then\global\@subeqncountfalse\fi
     \ignorespaces}
\def\figurenumstyle#1{\@setnumstyle\figurenum{#1}}
\def\subfigurenumstyle#1{\@setnumstyle\subfigurenum{#1}%
     \ifblank\subfigurenumstyle\then\global\@subfigcountfalse\fi
     \ignorespaces}
\def\tablenumstyle#1{\@setnumstyle\tablenum{#1}}
\def\subtablenumstyle#1{\@setnumstyle\subtablenum{#1}%
     \ifblank\subtablenumstyle\then\global\@subtblcountfalse\fi
     \ignorespaces}

\def\eqnlabel#1{%
     \if@subeqncount
          \newsubequationnum=\next
     \else\newequationnum=\next
          \ifblank\subequationnumstyle\then
          \else\global\@subeqncounttrue
               \newsubequationnum=\@ne
          \fi
     \fi
     \label{#1}{\puteqnformat}(\puteqn{#1})%
     \ifdraft\rlap{\hskip.1in{\tt#1}}\fi}

\let\puteqn=\putlab

\def\equation#1#2{\useafter\gdef{eqn@#1}{#2\eqno\eqnlabel{#1}}}
\def\Equation#1{\useafter\gdef{eqn@#1}}

\def\putequation#1{\useafter\ifundefined{eqn@#1}\then
     \err@undefinedeqn{#1}\else\usename{eqn@#1}\fi}

\def\eqnseriesstyle#1{\gdef\@eqnseriesstyle{#1}}
\def\begineqnseries{\subequationnumstyle{\@eqnseriesstyle}%
     \defaultoption[]\@begineqnseries}
\def\@begineqnseries[#1]{\edef\@@eqnname{#1}}
\def\endeqnseries{\subequationnumstyle{blank}%
     \expandafter\ifnoarg\@@eqnname\then
     \else\label\@@eqnname{\puteqnformat}%
     \fi
     \aftergroup\ignorespaces}

\def\figlabel#1{%
     \if@subfigcount
          \newsubfigurenum=\next
     \else\newfigurenum=\next
          \ifblank\subfigurenumstyle\then
          \else\global\@subfigcounttrue
               \newsubfigurenum=\@ne
          \fi
     \fi
     \label{#1}{\putfigformat}\putfig{#1}%
     {\def\marginnoteformat{\tt}\marginnote{#1}}}

\let\putfig=\putlab

\def\figseriesstyle#1{\gdef\@figseriesstyle{#1}}
\def\beginfigseries{\subfigurenumstyle{\@figseriesstyle}%
     \defaultoption[]\@beginfigseries}
\def\@beginfigseries[#1]{\edef\@@figname{#1}}
\def\endfigseries{\subfigurenumstyle{blank}%
     \expandafter\ifnoarg\@@figname\then
     \else\label\@@figname{\putfigformat}%
     \fi
     \aftergroup\ignorespaces}

\def\tbllabel#1{%
     \if@subtblcount
          \newsubtablenum=\next
     \else\newtablenum=\next
          \ifblank\subtablenumstyle\then
          \else\global\@subtblcounttrue
               \newsubtablenum=\@ne
          \fi
     \fi
     \label{#1}{\puttblformat}\puttbl{#1}%
     {\def\marginnoteformat{\tt}\marginnote{#1}}}

\let\puttbl=\putlab

\def\tblseriesstyle#1{\gdef\@tblseriesstyle{#1}}
\def\begintblseries{\subtablenumstyle{\@tblseriesstyle}%
     \defaultoption[]\@begintblseries}
\def\@begintblseries[#1]{\edef\@@tblname{#1}}
\def\endtblseries{\subtablenumstyle{blank}%
     \expandafter\ifnoarg\@@tblname\then
     \else\label\@@tblname{\puttblformat}%
     \fi
     \aftergroup\ignorespaces}

\message{reference numbering,}

\newcount\referencenum \referencenum=0
\newcount\@@prerefcount \@@prerefcount=0
\newcount\@@thisref
\newcount\@@lastref
\newcount\@@loopref
\newcount\@@refseq
\newdimen\refnumindent
\let\@undefreflist=\empty

\def\referencenumstyle#1{\@setnumstyle\referencenum{#1}}

\def\referencestyle#1{\usename{@ref#1}}

\def\@refsequential{%
     \gdef\@refpredef##1{\global\advance\referencenum by\@ne
          \let\^^\=0\label{##1}{\^^\{\the\referencenum}}%
          \useafter\gdef{ref@\the\referencenum}{{##1}{\undefinedlabelformat}}}%
     \gdef\@reference##1##2{%
          \ifundefinedlabel##1\then
          \else\def\^^\####1{\global\@@thisref=####1\relax}\putlab{##1}%
               \useafter\gdef{ref@\the\@@thisref}{{##1}{##2}}%
          \fi}%
     \gdef\endputreferences{%
          \loop\ifnum\@@loopref<\referencenum
                    \advance\@@loopref by\@ne
                    \expandafter\expandafter\expandafter\@printreference
                         \csname ref@\the\@@loopref\endcsname
          \repeat
          \par}}

\def\@refpreordered{%
     \gdef\@refpredef##1{\global\advance\referencenum by\@ne
          \additemR##1\to\@undefreflist}%
     \gdef\@reference##1##2{%
          \ifundefinedlabel##1\then
          \else\global\advance\@@loopref by\@ne
               {\let\^^\=0\label{##1}{\^^\{\the\@@loopref}}}%
               \@printreference{##1}{##2}%
          \fi}
     \gdef\endputreferences{%
          \def\^^\####1{\useiflab{####1}\then
               \else\reference{####1}{\undefinedlabelformat}\fi}%
          \@undefreflist
          \par}}

\def\beginprereferences{\par
     \def\reference##1##2{\global\advance\referencenum by1\@ne
          \let\^^\=0\label{##1}{\^^\{\the\referencenum}}%
          \useafter\gdef{ref@\the\referencenum}{{##1}{##2}}}}
\def\endprereferences{\global\@@prerefcount=\the\referencenum\par}

\def\beginputreferences{\par
     \refnumindent=\z@\@@loopref=\z@
     \loop\ifnum\@@loopref<\referencenum
               \advance\@@loopref by\@ne
               \setbox\z@=\hbox{\referencenum=\@@loopref
                    \referencenumformat\enskip}%
               \ifdim\wd\z@>\refnumindent\refnumindent=\wd\z@\fi
     \repeat
     \putreferenceformat
     \@@loopref=\z@
     \loop\ifnum\@@loopref<\@@prerefcount
               \advance\@@loopref by\@ne
               \expandafter\expandafter\expandafter\@printreference
                    \csname ref@\the\@@loopref\endcsname
     \repeat
     \let\reference=\@reference}

\def\@printreference#1#2{\ifx#2\undefinedlabelformat\err@undefinedref{#1}\fi
     \noindent\ifdraft\rlap{\hskip\hsize\hskip.1in \tt#1}\fi
     \llap{\referencenum=\@@loopref\referencenumformat\enskip}#2\par}

\def\reference#1#2{{\par\refnumindent=\z@\putreferenceformat\noindent#2\par}}

\def\putref#1{\trim#1\to\@@refarg
     \expandafter\ifnoarg\@@refarg\then
          \toks@={\relax}%
     \else\@@lastref=-\@m\def\@@refsep{}\def\@more{\@nextref}%
          \toks@={\@nextref#1,,}%
     \fi\the\toks@}
\def\@nextref#1,{\trim#1\to\@@refarg
     \expandafter\ifnoarg\@@refarg\then
          \let\@more=\relax
     \else\ifundefinedlabel\@@refarg\then
               \expandafter\@refpredef\expandafter{\@@refarg}%
          \fi
          \def\^^\##1{\global\@@thisref=##1\relax}%
          \global\@@thisref=\m@ne
          \setbox\z@=\hbox{\putlab\@@refarg}%
     \fi
     \advance\@@lastref by\@ne
     \ifnum\@@lastref=\@@thisref\advance\@@refseq by\@ne\else\@@refseq=\@ne\fi
     \ifnum\@@lastref<\z@
     \else\ifnum\@@refseq<\thr@@
               \@@refsep\def\@@refsep{,}%
               \ifnum\@@lastref>\z@
                    \advance\@@lastref by\m@ne
                    {\referencenum=\@@lastref\putrefformat}%
               \else\undefinedlabelformat
               \fi
          \else\def\@@refsep{--}%
          \fi
     \fi
     \@@lastref=\@@thisref
     \@more}

\message{streaming,}

\newif\ifstreaming

\def\streamto{\defaultoption[\jobname]\@streamto}
\def\@streamto[#1]{\global\streamingtrue
     \immediate\write\sixt@@n{--> Streaming to #1.str}%
     \newwrite\streamout\immediate\openout\streamout=#1.str }

\def\streamfrom{\defaultoption[\jobname]\@streamfrom}
\def\@streamfrom[#1]{\newread\streamin\openin\streamin=#1.str
     \ifeof\streamin
          \expandafter\err@nostream\expandafter{#1.str}%
     \else\immediate\write\sixt@@n{--> Streaming from #1.str}%
          \let\@@labeldef=\gdef
          \ifstreaming
               \edef\@elc{\endlinechar=\the\endlinechar}%
               \endlinechar=\m@ne
               \loop\read\streamin to\@@scratcha
                    \ifeof\streamin
                         \streamingfalse
                    \else\toks@=\expandafter{\@@scratcha}%
                         \immediate\write\streamout{\the\toks@}%
                    \fi
                    \ifstreaming
               \repeat
               \@elc
               \input #1.str
               \streamingtrue
          \else\input #1.str
          \fi
          \let\@@labeldef=\xdef
     \fi}

\def\streamcheck{\ifstreaming
     \immediate\write\streamout{\pagenum=\the\pagenum}%
     \immediate\write\streamout{\footnotenum=\the\footnotenum}%
     \immediate\write\streamout{\referencenum=\the\referencenum}%
     \immediate\write\streamout{\chapternum=\the\chapternum}%
     \immediate\write\streamout{\sectionnum=\the\sectionnum}%
     \immediate\write\streamout{\subsectionnum=\the\subsectionnum}%
     \immediate\write\streamout{\equationnum=\the\equationnum}%
     \immediate\write\streamout{\subequationnum=\the\subequationnum}%
     \immediate\write\streamout{\figurenum=\the\figurenum}%
     \immediate\write\streamout{\subfigurenum=\the\subfigurenum}%
     \immediate\write\streamout{\tablenum=\the\tablenum}%
     \immediate\write\streamout{\subtablenum=\the\subtablenum}%
     \immediate\closeout\streamout
     \fi}


\def\err@badtypesize{%
     \errhelp={The limited availability of certain fonts requires^^J%
          that the base type size be 10pt, 12pt, or 14pt.^^J}%
     \errmessage{--> Illegal base type size}}

\def\err@badsizechange{\immediate\write\sixt@@n
     {--> Size change not allowed in math mode, ignored}}

\def\err@sizetoolarge#1{\immediate\write\sixt@@n
     {--> \noexpand#1 too big, substituting HUGE}}

\def\err@sizenotavailable#1{\immediate\write\sixt@@n
     {--> Size not available, \noexpand#1 ignored}}

\def\err@fontnotavailable#1{\immediate\write\sixt@@n
     {--> Font not available, \noexpand#1 ignored}}

\def\err@sltoit{\immediate\write\sixt@@n
     {--> Style \noexpand\sl not available, substituting \noexpand\it}%
     \it}

\def\err@bfstobf{\immediate\write\sixt@@n
     {--> Style \noexpand\bfs not available, substituting \noexpand\bf}%
     \bf}

\def\err@badgroup#1#2{%
     \errhelp={The block you have just tried to close was not the one^^J%
          most recently opened.^^J}%
     \errmessage{--> \noexpand\end{#1} doesn't match \noexpand\begin{#2}}}

\def\err@badcountervalue#1{\immediate\write\sixt@@n
     {--> Counter (#1) out of bounds}}

\def\err@extrafootnotemark{\immediate\write\sixt@@n
     {--> \noexpand\footnotemark command
          has no corresponding \noexpand\footnotetext}}

\def\err@extrafootnotetext{%
     \errhelp{You have given a \noexpand\footnotetext command without first
          specifying^^Ja \noexpand\footnotemark.^^J}%
     \errmessage{--> \noexpand\footnotetext command has no corresponding
          \noexpand\footnotemark}}

\def\err@labelredef#1{\immediate\write\sixt@@n
     {--> Label "#1" redefined}}

\def\err@badlabelmatch#1{\immediate\write\sixt@@n
     {--> Definition of label "#1" doesn't match value in \jobname.lab}}

\def\err@needlabel#1{\immediate\write\sixt@@n
     {--> Label "#1" cited before its definition}}

\def\err@undefinedlabel#1{\immediate\write\sixt@@n
     {--> Label "#1" cited but never defined}}

\def\err@undefinedeqn#1{\immediate\write\sixt@@n
     {--> Equation "#1" not defined}}

\def\err@undefinedref#1{\immediate\write\sixt@@n
     {--> Reference "#1" not defined}}

\def\err@nostream#1{%
     \errhelp={You have tried to input a stream file that doesn't exist.^^J}%
     \errmessage{--> Stream file #1 not found}}

\message{jyTeX initialization}

\everyjob{\immediate\write16{--> jyTeX version \fmtversion}%
     \edef\@@jobname{\jobname}%
     \edef\jobname{\@@jobname}%
     \settime
     \openin0=\jobname.lab
     \ifeof0
     \else\closein0
          \immediate\write16{--> Getting labels from file \jobname.lab}%
          \input\jobname.lab
     \fi}


\def\fixedskipslist{%
     \^^\{\topskip}%
     \^^\{\splittopskip}%
     \^^\{\maxdepth}%
     \^^\{\skip\topins}%
     \^^\{\skip\footins}%
     \^^\{\headskip}%
     \^^\{\footskip}}

\def\scalingskipslist{%
     \^^\{\p@renwd}%
     \^^\{\delimitershortfall}%
     \^^\{\nulldelimiterspace}%
     \^^\{\scriptspace}%
     \^^\{\jot}%
     \^^\{\normalbaselineskip}%
     \^^\{\normallineskip}%
     \^^\{\normallineskiplimit}%
     \^^\{\baselineskip}%
     \^^\{\lineskip}%
     \^^\{\lineskiplimit}%
     \^^\{\bigskipamount}%
     \^^\{\medskipamount}%
     \^^\{\smallskipamount}%
     \^^\{\parskip}%
     \^^\{\parindent}%
     \^^\{\abovedisplayskip}%
     \^^\{\belowdisplayskip}%
     \^^\{\abovedisplayshortskip}%
     \^^\{\belowdisplayshortskip}%
     \^^\{\abovechapterskip}%
     \^^\{\belowchapterskip}%
     \^^\{\abovesectionskip}%
     \^^\{\belowsectionskip}%
     \^^\{\abovesubsectionskip}%
     \^^\{\belowsubsectionskip}}


\def\twoupsetup{
     \topmargin=.75in
     \leftmargin=.5in
     \vsize=6.9in
     \hsize=4.75in
     \fullhsize=10in
     \let\draft=\relax}

\outputstyle{normal}                             

\def\marginnoteformat{\subscriptsize             
     \hsize=1in \baselinestretch=1000 \everypar={}%
     \tolerance=5000 \hbadness=5000 \parskip=0pt \parindent=0pt
     \leftskip=0pt \rightskip=0pt \raggedright}

\head={\ifdraft\normalfonts\it\hfil DRAFT\hfil   
     \llap{\number\day\ \monthword\month\ \militarytime}\else\hfil\fi}
\foot={\hfil\normalfonts\numstyle\pagenum\hfil}  

\normalbaselineskip=12pt                         
\normallineskip=0pt                              
\normallineskiplimit=0pt                         
\normalbaselines                                 

\topskip=.85\baselineskip \splittopskip=\topskip \headskip=2\baselineskip
\footskip=\headskip

\pagenumstyle{arabic}                            

\parskip=0pt                                     
\parindent=20pt                                  

\baselinestretch=1000                            


\chapterstyle{left}                              
\chapternumstyle{blank}                          
\def\chapterbreak{\newpage}                      
\abovechapterskip=0pt                            
\belowchapterskip=1.5\baselineskip               
     plus.38\baselineskip minus.38\baselineskip
\def\chapternumformat{\numstyle\chapternum.}     

\sectionstyle{left}                              
\sectionnumstyle{blank}                          
\def\sectionbreak{\vskip0pt plus4\baselineskip\penalty-100
     \vskip0pt plus-4\baselineskip}              
\abovesectionskip=1.5\baselineskip               
     plus.38\baselineskip minus.38\baselineskip
\belowsectionskip=\the\baselineskip              
     plus.25\baselineskip minus.25\baselineskip
\def\sectionnumformat{
     \ifblank\chapternumstyle\then\else\numstyle\chapternum.\fi
     \numstyle\sectionnum.}

\subsectionstyle{left}                           
\subsectionnumstyle{blank}                       
\def\subsectionbreak{\vskip0pt plus4\baselineskip\penalty-100
     \vskip0pt plus-4\baselineskip}              
\abovesubsectionskip=\the\baselineskip           
     plus.25\baselineskip minus.25\baselineskip
\belowsubsectionskip=.75\baselineskip            
     plus.19\baselineskip minus.19\baselineskip
\def\subsectionnumformat{
     \ifblank\chapternumstyle\then\else\numstyle\chapternum.\fi
     \ifblank\sectionnumstyle\then\else\numstyle\sectionnum.\fi
     \numstyle\subsectionnum.}


\footnotenumstyle{symbols}                       
\footnoteskip=0pt                                
\def\footnotenumformat{\numstyle\footnotenum}    
\def\footnoteformat{\footnotesize                
     \everypar={}\parskip=0pt \parfillskip=0pt plus1fil
     \leftskip=1em \rightskip=0pt
     \spaceskip=0pt \xspaceskip=0pt
     \def\\{\ifhmode\ifnum\lastpenalty=-10000
          \else\hfil\penalty-10000 \fi\fi\ignorespaces}}


\def\undefinedlabelformat{$\bullet$}             


\equationnumstyle{arabic}                        
\subequationnumstyle{blank}                      
\figurenumstyle{arabic}                          
\subfigurenumstyle{blank}                        
\tablenumstyle{arabic}                           
\subtablenumstyle{blank}                         

\eqnseriesstyle{alphabetic}                      
\figseriesstyle{alphabetic}                      
\tblseriesstyle{alphabetic}                      

\def\puteqnformat{\hbox{
     \ifblank\chapternumstyle\then\else\numstyle\chapternum.\fi
     \ifblank\sectionnumstyle\then\else\numstyle\sectionnum.\fi
     \ifblank\subsectionnumstyle\then\else\numstyle\subsectionnum.\fi
     \numstyle\equationnum
     \numstyle\subequationnum}}
\def\putfigformat{\hbox{
     \ifblank\chapternumstyle\then\else\numstyle\chapternum.\fi
     \ifblank\sectionnumstyle\then\else\numstyle\sectionnum.\fi
     \ifblank\subsectionnumstyle\then\else\numstyle\subsectionnum.\fi
     \numstyle\figurenum
     \numstyle\subfigurenum}}
\def\puttblformat{\hbox{
     \ifblank\chapternumstyle\then\else\numstyle\chapternum.\fi
     \ifblank\sectionnumstyle\then\else\numstyle\sectionnum.\fi
     \ifblank\subsectionnumstyle\then\else\numstyle\subsectionnum.\fi
     \numstyle\tablenum
     \numstyle\subtablenum}}


\referencestyle{sequential}                      
\referencenumstyle{arabic}                       
\def\putrefformat{\numstyle\referencenum}        
\def\referencenumformat{\numstyle\referencenum.} 
\def\putreferenceformat{
     \everypar={\hangindent=1em \hangafter=1 }%
     \def\\{\hfil\break\null\hskip-1em \ignorespaces}%
     \leftskip=\refnumindent\parindent=0pt \interlinepenalty=1000 }


\normalsize


\def\fmtversion{2.6M (June 1992)}

\catcode`\@=12

\typesize=10pt \magnification=1200 \baselineskip17truept
\footnotenumstyle{arabic} \hsize=6truein\vsize=8.5truein
\sectionnumstyle{blank}
\chapternumstyle{blank}
\chapternum=1
\sectionnum=1
\pagenum=0

\def\begintitle{\pagenumstyle{blank}\parindent=0pt
\begin{narrow}[0.4in]}
\def\endtitle{\end{narrow}\newpage\pagenumstyle{arabic}}


\def\beginexercise{\vskip 20truept\parindent=0pt\begin{narrow}[10
truept]}
\def\endexercise{\vskip 10truept\end{narrow}}


\def\eql#1{\eqno\eqnlabel{#1}}
\def\ref{\reference}
\def\peq{\puteqn}
\def\pref{\putref}

\def\mgn{\marginnote}
\def\bex{\begin{exercise}}
\def\eex{\end{exercise}}


\font\open=msbm10 

 
\def\StretchRtArr#1{{\count255=0\loop\relbar\joinrel\advance\count255 by1
\ifnum\count255<#1\repeat\rightarrow}}
\def\StretchLtArr#1{\,{\leftarrow\!\!\count255=0\loop\relbar
\joinrel\advance\count255 by1\ifnum\count255<#1\repeat}}

\def\StretchLRtArr#1{\,{\leftarrow\!\!\count255=0\loop\relbar\joinrel\advance
\count255 by1\ifnum\count255<#1\repeat\rightarrow\,\,}}

\def\mbox#1{{\leavevmode\hbox{#1}}}

\def\hspace#1{{\phantom{\mbox#1}}}
\def\oZ{\mbox{\open\char90}}

\def\b0{{\bf0}}
\def\bxi{{\bmit\xi}}
\def\bze{{\bmit\zeta}}

\def\de{\delta}

\def\th{\theta}

\def\bze{{\bmit\zeta}}
\def\De{\Delta}

\def\caG{{\cal G}}

\def\sc{{\rm sc }}


\def\frac#1/#2{\leavevmode\kern.1em
\raise.5ex\hbox{\the\scriptfont0 #1}\kern-.1em/\kern-.15em
\lower.25ex\hbox{\the\scriptfont0 #2}}
\def\sfrac#1/#2{\leavevmode\kern.1em
\raise.5ex\hbox{\the\scriptscriptfont0 #1}\kern-.1em/\kern-.15em
\lower.25ex\hbox{\the\scriptscriptfont0 #2}}

\def\gtorder{\mathrel{\raise.3ex\hbox{$>$}\mkern-14mu
             \lower0.6ex\hbox{$\sim$}}}
\def\ltorder{\mathrel{\raise.3ex\hbox{$<$}\mkern-14mu
             \lower0.6ex\hbox{$\sim$}}}

\def\semidirprod{\rlap{\ss C}\raise1pt\hbox{$\mkern.75mu\times$}}
\def\for{\lower6pt\hbox{$\Big|$}}
\def\fish{\kern-.25em{\phantom{abcde}\over \phantom{abcde}}\kern-.25em}


\def\boxit#1{\vbox{\hrule\hbox{\vrule\kern3pt
        \vbox{\kern3pt#1\kern3pt}\kern3pt\vrule}\hrule}}
\def\dalemb#1#2{{\vbox{\hrule height .#2pt
        \hbox{\vrule width.#2pt height#1pt \kern#1pt \vrule
                width.#2pt} \hrule height.#2pt}}}

\def\frac#1#2{{{#1}\over{#2}}}

\def\noin{\noindent}

\def\comb#1#2{{\left(#1\atop#2\right)}}

\def\etc{{\it etc. }}

\def\eg{{\it e.g.}}
\def\ie{{\it i.e. }}
\def\cf{{\it cf }}



\def\cosech{{\rm cosech }}
\def\sech{{\rm sech }}

\def\3j#1#2#3#4#5#6{\left\lgroup\matrix{#1&#2&#3\cr#4&#5&#6\cr}
\right\rgroup}

\def\m?{\mgn{?}}

\def\beq{\begin{eqnarray}}
\def\eeq{\end{eqnarray}}


\def\aop#1#2#3{{\it Ann. Phys.} {\bf {#1}} ({#2}) #3}
\def\cjp#1#2#3{{\it Can. J. Phys.} {\bf {#1}} ({#2}) #3}
\def\cmp#1#2#3{{\it Comm. Math. Phys.} {\bf {#1}} ({#2}) #3}
\def\cqg#1#2#3{{\it Class. Quant. Grav.} {\bf {#1}} ({#2}) #3}

\def\ijmp#1#2#3{{\it Int. J. Mod. Phys.} {\bf {#1}} ({#2}) #3}

\def\jmp#1#2#3{{\it J. Math. Phys.} {\bf {#1}} ({#2}) #3}
\def\jpa#1#2#3{{\it J. Phys.} {\bf A{#1}} ({#2}) #3}
\def\jpamt#1#2#3{{\it J. Phys.A:Math.Theor.} {\bf{#1}} ({#2}) #3}
\def\jpc#1#2#3{{\it J. Phys.} {\bf C{#1}} ({#2}) #3}
\def\lnm#1#2#3{{\it Lect. Notes Math.} {\bf {#1}} ({#2}) #3}

\def\np#1#2#3{{\it Nucl. Phys.} {\bf B{#1}} ({#2}) #3}
\def\npa#1#2#3{{\it Nucl. Phys.} {\bf A{#1}} ({#2}) #3}
\def\pl#1#2#3{{\it Phys. Lett.} {\bf {#1}} ({#2}) #3}

\def\prp#1#2#3{{\it Phys. Rep.} {\bf {#1}} ({#2}) #3}
\def\pr#1#2#3{{\it Phys. Rev.} {\bf {#1}} ({#2}) #3}
\def\prA#1#2#3{{\it Phys. Rev.} {\bf A{#1}} ({#2}) #3}

\def\prD#1#2#3{{\it Phys. Rev.} {\bf D{#1}} ({#2}) #3}
\def\prE#1#2#3{{\it Phys. Rev.} {\bf E{#1}} ({#2}) #3}
\def\prl#1#2#3{{\it Phys. Rev. Lett.} {\bf #1} ({#2}) #3}

\def\rmp#1#2#3{{\it Rev. Mod. Phys.} {\bf {#1}} ({#2}) #3}

\def\zfp#1#2#3{{\it Z. f. Phys.} {\bf {#1}} ({#2}) #3}

\def\cras#1#2#3{{\it Comptes Rend. Acad. Sci. (Paris)} {\bf{#1}} (#2) #3}
\def\prs#1#2#3{{\it Proc. Roy. Soc.} {\bf A{#1}} ({#2}) #3}
\def\pcps#1#2#3{{\it Proc. Camb. Phil. Soc.} {\bf{#1}} ({#2}) #3}
\def\mpcps#1#2#3{{\it Math. Proc. Camb. Phil. Soc.} {\bf{#1}} ({#2}) #3}

\def\amsh#1#2#3{{\it Abh. Math. Sem. Ham.} {\bf {#1}} ({#2}) #3}
\def\am#1#2#3{{\it Acta Mathematica} {\bf {#1}} ({#2}) #3}
\def\aim#1#2#3{{\it Adv. in Math.} {\bf {#1}} ({#2}) #3}
\def\ajm#1#2#3{{\it Am. J. Math.} {\bf {#1}} ({#2}) #3}
\def\amm#1#2#3{{\it Am. Math. Mon.} {\bf {#1}} ({#2}) #3}

\def\aom#1#2#3{{\it Ann. of Math.} {\bf {#1}} ({#2}) #3}
\def\cjm#1#2#3{{\it Can. J. Math.} {\bf {#1}} ({#2}) #3}
\def\bams#1#2#3{{\it Bull.Am.Math.Soc.} {\bf {#1}} ({#2}) #3}

\def\cmh#1#2#3{{\it Comm. Math. Helv.} {\bf {#1}} ({#2}) #3}

\def\dmj#1#2#3{{\it Duke Math. J.} {\bf {#1}} ({#2}) #3}
\def\invm#1#2#3{{\it Invent. Math.} {\bf {#1}} ({#2}) #3}

\def\jdg#1#2#3{{\it J. Diff. Geom.} {\bf {#1}} ({#2}) #3}

\def\joa#1#2#3{{\it J. of Algebra} {\bf {#1}} ({#2}) #3}
\def\jram#1#2#3{{\it J. f. Reine u. Angew. Math.} {\bf {#1}} ({#2}) #3}
\def\jims#1#2#3{{\it J. Indian. Math. Soc.} {\bf {#1}} ({#2}) #3}
\def\jlms#1#2#3{{\it J. Lond. Math. Soc.} {\bf {#1}} ({#2}) #3}
\def\jmpa#1#2#3{{\it J. Math. Pures. Appl.} {\bf {#1}} ({#2}) #3}
\def\ma#1#2#3{{\it Math. Ann.} {\bf {#1}} ({#2}) #3}

\def\mz#1#2#3{{\it Math. Zeit.} {\bf {#1}} ({#2}) #3}
\def\ojm#1#2#3{{\it Osaka J.Math.} {\bf {#1}} ({#2}) #3}

\def\pems#1#2#3{{\it Proc. Edin. Math. Soc.} {\bf {#1}} ({#2}) #3}

\def\plb#1#2#3{{\it Phys. Letts.} {\bf {B#1}} ({#2}) #3}
\def\pla#1#2#3{{\it Phys. Letts.} {\bf {A#1}} ({#2}) #3}
\def\plms#1#2#3{{\it Proc. Lond. Math. Soc.} {\bf {#1}} ({#2}) #3}
\def\pgma#1#2#3{{\it Proc. Glasgow Math. Ass.} {\bf {#1}} ({#2}) #3}
\def\qjm#1#2#3{{\it Quart. J. Math.} {\bf {#1}} ({#2}) #3}
\def\qjpam#1#2#3{{\it Quart. J. Pure and Appl. Math.} {\bf {#1}} ({#2}) #3}

\def\rmjm#1#2#3{{\it Rocky Mountain J. Math.} {\bf {#1}} ({#2}) #3}

\def\tams#1#2#3{{\it Trans.Am.Math.Soc.} {\bf {#1}} ({#2}) #3}

\begin{title}
\vglue 0.5truein
\vskip15truept
\centertext {\Bigfonts \bf Central differences, Euler numbers} \vskip2truept
\vskip10truept\centertext{\Bigfonts \bf and symbolic methods}
 \vskip 20truept
\centertext{J.S.Dowker\footnote{dowker@man.ac.uk; dowkeruk@yahoo.co.uk}} \vskip
7truept \centertext{\it Theory Group,} \centertext{\it School of Physics and Astronomy,}
\centertext{\it The University of Manchester,} \centertext{\it Manchester, England} \vskip
7truept \centertext{}

\vskip 7truept

\vskip40truept
\begin{narrow}
I relate some coefficients encountered when computing the functional determinants on
spheres to the central differentials of nothing. In doing this I use some historic works, in
particular, transcribing the elegant symbolic formalism of Jeffery (1861) into central
difference form which has computational advantages for Euler numbers, as discovered by
Shovelton (1915). I derive sum rules for these, and for the central differentials, the proof of
which involves an interesting expression for powers of sech x as multiple derivatives. I
present a more general, symbolic treatment of central difference calculus which allows
known, and unknown, things to be obtained in an elegant and compact fashion gaining, at no
cost, the expansion of the powers of the inverse sinh, a basic central function. Systematic
use is made of the operator 2 asinh(D/2). Umbral calculus is employed to compress the
operator formalism. For example the orthogonality/completeness of the factorial numbers,
of the first and second kinds, translates, umbrally, to T(t(x))=x. The classic expansions of
multiple angle cosh and sinh in terms of powers of sinh are thence obtained with minimal
effort.
\end{narrow}
\vskip 5truept
\vskip 60truept
\vfil
\end{title}
\pagenum=0
\newpage

\section{\bf 1. Introduction}

In a work, [\pref{Dowren}], involving quantum field theory calculations of the effective
action on spheres (as the logarithm of the operator determinant of the field propagation
operator), in order to show the equivalence of two explicit expressions for this `logdet', I
was obliged to look at an expansion of the powers of the $\sech$ function. Some details are
in the following section.

As background, the power series expansion of $\sec x$ goes back to Euler, the coefficients
being named Euler numbers in the 1820's to 1840's, \eg\ by Scherk and Raabe. They were,
and still are, called secant numbers. Their evaluation is distinct from those of $\tan x$ which
are related to the Bernoulli numbers. The series expansion of {\it powers} of $\sech\, x$
yields {\it generalised} Euler numbers.

The study of Bernoulli numbers, their relatives and their generalisations remains strong
after nearly three hundred years. For this reason, it is difficult to know whether any
displayed result is really new and, not unexpectedly, the expansion in [\pref{Dowren}]
turns out to be known and I wish here to go into more detail about this interesting topic
which has involved me in an analysis of the existing, and sometimes ancient, literature. The
Dilcher online bibliography contains most (not all) of these references. Whether anyone
other than the historically minded has actually used, or even examined, this old material is
open to question.

The relevant calculus turns out to be that of central differences and in section 5 I transcribe
some of the symbolic formalism given by Jeffery, [\pref{Jeffery}], who uses forward
differences, into central calculus. This allows various sum rules and identities to be derived,
which might be new, as well as some (known) series expansions.

In sections 6 and 7 I further develope  the symbolic approach to central calculus in a more
general way presenting some standard and non--standard applications in a (possibly) novel
fashion. In section 8, I combine this approach with the umbral notation, generating
functions, and other things, emerging effortlessly.

\section{\bf2. The expansion and central factorials}
The integral needed in [\pref{Dowren}] is,
  $$
   \int_0^\infty dx\,
    {1\over (x^2+1)\cosh^{2r+1} \pi x/2}\,,
    \eql{jay}
  $$
and, in order to allow partial integration, one requires the expansion,
  $$
    \sech^{2r+1}x={1\over(2r)!}\sum_{\rho=0}^r(-1)^\rho\,\caG^r_\rho \,
    {d^{2\rho}\over dx^{2\rho}}\,\sech\, x\,,
    \eql{sechs}
  $$
for odd powers in terms of even derivatives.\footnote{ I have improved the definition of the
coefficients in order to make them integers.} The coefficients satisfy the recursion,
(difference equation),
   $$
    \caG^r_\rho=(2r-1)^2\,\caG^{r-1}_\rho+\caG^{r-1}_{\rho-1}\,,
    \eql{erecur}
  $$
together with $\caG^r_{-1}=0=\caG^r_{r+1}$ and the initial value $\caG^0_0=1$.
Special values are $\caG^r_r=1$. They are directly related to central difference factorial
numbers (\eg\ Steffensen, [\pref{Steffensen}]) as is easily seen by repeated action of the
basic, and old, relation,
  $$
  \sech^nx=-{1\over(n-1)(n-2)}\big(D^2-(n-2)^2\big)\,\sech^{n-2} x
  \,,\quad D={d\over dx}\,,
  $$
which yields,
  $$\eqalign{
  \sech^{2r+1}x&={(-1)^r\over(2r)!}\big(D^2-(2r-1)^2\big)\big(D^2-(2r-3)^2\big)
  \ldots\big(D^2-1\big)\,\sech\, x\cr
  &\equiv2^{2r}\bigg({D\over2}\bigg)^{[2r+1]-1}\,\sech\, x\,,
  }
  \eql{iter}
  $$
in Steffensen's notation. The coefficients of the expansion of the central factorial in ordinary
powers of $D$ are, by definition, the central factorial numbers of the first kind (or, up to a
factor, the `central differentials of nothing')\footnote{ One might call them `central Stirling
numbers of the first kind'.}. The complete iteration, (\peq{iter}), shows that these are
equal to the sum of the squares of the first $r$ odd integers taken a certain number of times
(\cf Ely [\pref{Ely}], Jeffery [\pref{Jeffery}]) but I will not continue with this combinatorial
view just now.

The precise connection between the coefficients in (\peq{sechs}) and the central differentials
of nothing is,\footnote{ A slightly more complicated treatment, but essentially the same, is
given by Zhang, [\pref{Zhang}].}
  $$\eqalign{
  \caG^r_\rho&=(-1)^{r-\rho}\,2^{2r}{(2r)!\over(2\rho+1)!}\,D^{2\rho+1}\,\b0^{[2r+1]}\cr
  &=(-1)^{r-\rho}\,2^{2r}(2r)!\,\,t(2r+1,2\rho+1)\,,
  }
  \eql{gees}
  $$
where I have introduced the notation of Riordan, [\pref{Riordan}], for the central factorial
numbers, $t$. The recursion, (\peq{erecur}), is then equivalent to those in Steffensen,
[\pref{Steffensen}], and in [\pref{Riordan}] and [\pref{BSSV}].

Writing (\peq{sechs}) this way round was the most convenient for the purposes in
[\pref{Dowren}]. One could, of course, invert it and, in this form, a related expansion is
given by Stern, [\pref{Stern}], which I wish to relate to (\peq{sechs}). He actually gives
two expansions, one for even derivatives,
  $$
  {d^{2\rho}\over dx^ {2\rho}}\,\sech\, x=(-1)^\rho \,\sech x\sum_{k=0}^\rho(-1)^k
  a^{(2\rho)}_k\tanh^{2k}x\,,
  \eql{sterne}
  $$
and one for odd,
  $$\eqalign{
  {d^{2\rho+1}\over dx^ {2\rho+1}}\,\sech\, x&=(-1)^{\rho+1} \cosech\, x
  \sum_{k=0}^\rho(-1)^k a^{(2\rho+1)}_k\tanh^{2k+2}x\cr
  &=(-1)^{\rho+1} \sech\, x\tanh x
  \sum_{k=0}^\rho(-1)^k a^{(2\rho+1)}_k\tanh^{2k}x\,.
  }
  \eql{sterno}
  $$
One sees the appearance of `derivative polynomials' which have been subject to more
recent investigations \eg\ [\pref{Cvijovic1}].
  
Stern gives the coupled recursions,
  $$\eqalign{
   a^{(2n+1)}_k&=(2k+1)a^{(2n)}_k+(2k+2)a^{(2n)}_{k+1}\cr
   a^{(2n+2)}_k&=2k\,a^{(2n+1)}_{k-1}+(2k+1)a^{(2n+1)}_k\,.\cr
   }
   \eql{strecur}
  $$
  
The coefficients $\caG^r_\rho$ can be connected to the $ a^{(2\rho)}_k$ by substituting
(\peq{sterne}) into (\peq{sechs}) giving,
  $$
   (1-\tanh^2x)^r=\sum_{\rho=0}^r\caG^r_\rho\,
    \sum_{k=0}^\rho(-1)^k a^{(2\rho)}_k\tanh^{2k}x\,,
  $$
whence,
  $$
   \sum_{\rho=s}^r\caG^r_\rho\,
   a^{(2\rho)}_s=(2r)!\,\comb rs\,.
   \eql{sp}
  $$

Similarly, substituting (\peq{sterno}) into the derivative of (\peq{sechs}) produces
  $$
  \sum_{\rho=s}^r\caG^r_\rho\,
   a^{(2\rho+1)}_s=(2r+1)!\comb rs\,.
   \eql{sp2}
  $$

The second recurrence in (\peq{strecur}) is equivalent, directly, to the binomial recurrence
  $$
  \comb r{s+1}={r-s\over s+1}\,\comb rs\,,\quad \comb r0=1\,.
  $$

\section{\bf3. Generalised Euler numbers}
Generalised Euler numbers can be defined by the power series (I use N\"orlund's notation,
[\pref{Norlund,Norlund1}]),
  $$
  \sech ^r x=\sum_\nu {x^{2\nu}\over(2\nu)!} E^{(r)}_{2\nu}\,.
  \eql{geneul}
  $$
The original Euler numbers correspond to $n=1$, $E_{2\nu}\equiv E^{(1)}_{2\nu}$.

Setting $x$ to zero in (\peq{sterne}) and using the MacLaurin expansion, yields\mgn{Check
sign}
  $$
  E_{2\nu}=a^{(2\nu)}_0=a^{(2\nu-1)}_0\,,
  $$
where the second equality follows on differentiating (\peq{sterno}),  and so Stern refers to
the other coefficients, $a^*_k$, as {\it Euler numbers of higher order} but this name is
more often used for the numbers defined by (\peq{geneul}). Stern's term reappears in
some recent works, \eg\ [\pref{Cvijovic1}].

By differentiating (\peq{sechs}) $2\nu$ times, one can easily obtain an expression for the
odd generalised numbers in terms of the ordinary ones,
  $$
    E^{(2r+1)}_{2\nu}={1\over(2r)!}\sum_{\rho=0}^r(-1)^\rho
    \,\caG^r_\rho\,E_{2\nu+2\rho}\,,\quad \forall\, \nu\,.
    \eql{reln}
  $$
This relation can be found in Ely, [\pref{Ely}], where the coefficients, $\caG$, are expressed
in the combinatorial fashion.

Schwatt, [\pref{Schwatt}], pp.68--71, has an extensive algebraic treatment of expansions
of powers of $\sech\,x$ and also of the derivatives of these in the Stern form,
(\peq{sterne}), (\peq{sterno}).

One could look upon (\peq{reln}) as a means of computing the (odd) generalised numbers
since the recursion (\peq{erecur}) for $\caG$ is easily solved and the first order Euler
numbers are fairly well studied. In particular an efficient recursion method is adopted by
Knuth and Buckholtz, [\pref{KnandB}], which it is worth recalling and then comparing with
earlier work. (See also Brent and Harvey, [\pref{BrandH}]. )
\section{\bf4. Computing Euler numbers}

Knuth and Buckholtz use trigonometric functions and, by differentiation, deduce the relation,
  $$
  {d^{n}\over dz^ {n}}\,\sec z=\sec z(E_{n0}+E_{n1}\tan z+E_{n2}\tan^2 z+\ldots)\,,
  \eql{KB}
  $$
and the recursion,
  $$
  E_{n+1,k}=kE_{n,k-1}+(k+1)E_{n,k+1}\,,\quad E_{0k}=\de_{0k}\,.
  \eql{KBrecurs}
  $$
Since  $E_n=E_{n0}$ this allows an efficient computation of the Euler numbers. \footnote{ A
simple algorithm implementing (\peq{KB}) is given by Brent and Harvey. Programming this
in DERIVE (not a fast CAS) I calculated the first 1000 Euler numbers in 24 seconds using an
Athlon 610e processor. The values could be used in a look--up table for fast access.}

The earlier Stern recursions are the same as (\peq{KB}) which is easily seen by splitting
(\peq{KB}) into two by the parity behaviour under $z\to-z$. Then one sees that,
  $$
   E_{2\rho,2k}=a^{2\rho}_k\,,\quad E_{2\rho,2k+1}=0\,,
  $$
and
  $$
  E_{2\rho+1,2k+1}=a^{(2\rho+1)}_k\,,\quad  E_{2\rho+1,2k}=0\,.
  $$
  
Stern uses the recursions to derive many arithmetical properties, as do Knuth and
Buckholtz, who were unaware of Stern's work.

Historical computations of the secant numbers were performed by Euler and by Scherk (in
1825). I do not give these continental references as they can be found in the 1893 text by
Saalsch\"utz, [\pref{Saalschutz}], as well as a list of the first 14 numbers. Glaisher
extended the list to 27 and Joffe to 50. The effort seems to be more of a computational
challenge than necessity.

Herschel, [\pref{Herschel,Herschel3}], (1816 and 1820) gave finite difference expressions
suitable for numerical computation (in terms of forward differences of nothing) and also
related Euler numbers to Bernoulli numbers via tangent numbers.

Herschel's theorem applied to $\sech\, x$ yields,
  $$
   E_{2n}=(-1)^n{2\over (1+\De)+(1+\De)^{-1}}\,{\b0}^{2n}\,,
   \eql{hersch}
  $$
which is manipulated in various ways. However, a more elegant, and practical, formulation
was found by Shovelton, [\pref{Shovelton}], in terms of the {\it central} difference
operator, $\de$, as, \footnote{ Shovelton developes his results via the Bernoulli
polynomials (the `Bernoullian function') and introduces central differences as a simplifying
tactic.}
  $$
  E_{2n}=(-1)^n{1\over 1+\de^2/2}\,{\b0}^{2n}\equiv (-1)^n
 \bigg(1-{1\over2}\,\de^2+{1\over2^2}\,\de^4+\ldots+
 {(-1)^n\over2^{n}}\,\de^{2n}\bigg){\b0}^{2n}\,,
  \eql{Shov}
  $$
which can be obtained directly from (\peq{hersch}) by noting that $\de=\De\,E^{-1/2}$
($E$ being the translation operator, $E=1+\De$) so that (\eg\ [\pref{Shovelton}],
[\pref{Steffensen}] pp.10,11, [\pref{Bickley}]),
  $$
  \de=E^{1/2}-E^{-1/2}\,,\quad i.e.\,\,\quad 2+\de^2=E+E^{-1}\,,
  $$
and also,
  $$
 ( E^{1/2}+E^{-1/2})^2=4+\de^2\,.
 \eql{mean}
  $$

After use of the equivalence $ f(E^n)\,\b0^s=n^s\,f(E)\,\b0^s\,,\forall n$  (Herschel,
[\pref{Herschel2}] VII \S(8).) equations (\peq{mean}) and (\peq{hersch}) allow the
alternative compact expression, also given by Shovelton, [\pref{Shovelton}],
  $$\eqalign{
  E_{2n}&={(-1)^n\over 2^{2n}}{1\over (1+\de^2/4)^{1/2}}\,{\b0}^{2n}\cr
  &\equiv {(-1)^n\over2^{2n}}
 \bigg[1+\comb{-{1\over2}}1{\de^2\over2^2}
 +\comb{-{1\over2}}2{\de^4\over2^4}+\ldots+
 \comb{-{1\over2}}n{\de^{2n}\over2^{2n}}\bigg]{\b0}^{2n}\,.
 }
  \eql{Shov2}
  $$
(See Butzer {\it et al}, [\pref{BSSV}].)

 Joffe, [\pref{Joffe}], undertook the numerical evaluation of $E_{2n}$ using
(\peq{Shov}) in terms of the central differences of nothing which were computed using the
recursion formula. An early table of values is provided by Shovelton, [\pref{Shovelton}].

The expansion of the powers of $\sec$ is formally no harder and by the same process one
has the coefficients, [\pref{Shovelton}],
  $$\eqalign{
    E^{(r)}_{2n}&=(-1)^n{1\over (1+\de^2/2)^r}\,{\b0}^{2n}\cr
    &\equiv (-1)^n
 \bigg[1+{1\over2}\comb {-r}1\,\de^2+{1\over2^2}\comb {-r}2\,\de^4+\ldots+
 {1\over2^{n}}\comb {-r}n \,\de^{2n}\bigg]{\b0}^{2n}\,,
 }
  \eql{Shov3}
  $$
which has the advantage over (\peq{reln}) of displaying the polynomial nature of
$E^{(r)}_{2n}$ in $r$.

The form corresponding to (\peq{Shov2}) is,
  $$\eqalign{
    E^{(r)}_{2n}&={(-1)^n\over2^{2n}}{1\over (1+\de^2/4)^{r/2}}\,{\b0}^{2n}\cr
    &\equiv {(-1)^n\over2^{2n}}
 \bigg[1+\comb {-{r\over2}}1\,{\de^2\over2^2}+
 \comb {-{r\over2}}2\,{\de^4\over2^4}+\ldots+
 \comb {-{r\over2}}n\,{\de^{2n}\over2^{2n}}\bigg]{\b0}^{2n}\,.
 }
  \eql{Shov4}
  $$

The central differences of nothing are, to a factor, the central factorial numbers of the {\it
second} kind, \eg\ Steffensen, [\pref{Steffensen}], and [\pref{Riordan}] and
[\pref{BSSV}].

Earlier expansions of $\sech ^n\,x$ exist but are generally expressed in terms of the
forward differences of nothing. For example Jeffery, [\pref{Jeffery2}], gives the neat
looking, but computationally somewhat awkward, form,
  $$
    \sec^nx=\bigg({2\over2+\De}\bigg)^n\,\cos(2.\b0+n)x\,.
  $$
  
An expression, essentially equivalent to (\peq{Shov3}), has been given more recently by
Liu and Zhang, [\pref{LiandZ}].
\section{\bf5. Jeffery's theorem, central differences and Euler relations}

The early works of Herschel, Jeffery, Boole and others are characterised by the extensive
use of symbolic methods. Shovelton does likewise, but in terms, partly, of central
differences. The 1909 textbook by Thiele, [\pref{Thiele}], is a rare example of their
application by a continental worker. Another is the classic text by Steffensen,
[\pref{Steffensen}], both geared to interpolation questions. There is also an 1865 work by
Hansen that I have not been able to access. Bickley, [\pref{Bickley}], gives a useful
summary, without references, of standard operator relations. Chapter 6 in Riordan,
[\pref{Riordan}], is a good, modern pedagogic treatment.

The elegant papers by Jeffery contain applications of four basic theorems concerning the
transformation of symbolic expressions. (Herschel's much earlier work also involves a lot of
these manipulations.) The formalism is expressed in terms of forward differences.
Shovelton's results show that considerable simplification ensues if central differences are
brought into play and I wish here to transcribe some of Jeffery's formulae into central form
and apply them to obtain Euler number relations as an illustration.

For this purpose, I return to the alternative formula, (\peq{Shov2}), for the Euler numbers,
and write it equivalently as,
  $$\eqalign{
  E_{2n}&={(-1)^n\over2^{2n}}\,{2\over E^{1/2}+E^{-1/2}}\,\b0^{2n}\cr
  &={(-1)^n\over2^{2n}}\, \,\sech {D\over2}\,\,\b0^{2n}\,.
  }
  $$

Rather than expanding the operator in powers of $\de^2$ directly (as in (\peq{Shov2})), I
substitute  for the $\b0^{2n}$ in central factorials using,
  $$
  x^{2n}=\sum_{\nu=1}^n x^{[2\nu]}\, {\de^{2\nu}\b0^{2n}\over(2\nu)!}\,,
  $$
to give the intermediate formula,
 $$\eqalign{
  E_{2n}&={(-1)^n\over2^{2n}} \,\sum_{\nu=1}^n \sech {D\over2}\,\,\b0^{[2\nu]}\cdot
  {\de^{2\nu}\b0^{2n}\over(2\nu)!}\,.
  }
  \eql{eul}
  $$

I am following Jeffery, [\pref{Jeffery2}] \S3, and, as there, a theorem is employed to
re--express the action of the pre--operator, here, $\sech\,D/2$. The theorem is the central
difference form of Theorem $D$ in Jeffery, [\pref{Jeffery}], and reads (I prove this in
section 6),
  $$
   f(D)\,\b0^{[n]}={d\,f(D)\over d\de}\,\b0^{[n-1]}=\ldots={d^nf(D)\over d\de^n}\,1\,,
   \eql{D}
  $$
which can be run left to right or right to left. Applied to powers of $\sech\, D/2$, taking it all
the way to the right trivially gives the expansion (\peq{Shov4}), and so nothing new.

More interesting is to work right to left. For example, set $n=2\nu+1$ and choose,
   $$\eqalign{
    f(D)&=\int d\de\,\sech\,{D\over2}=\int dD\,\cosh{D\over2}\,\sech\,{D\over2}\cr
    &=D\,,
    }
    \eql{effd}
   $$
(since $\de=2\sinh\,D/2$) so that (\peq{eul}) becomes,
  $$\eqalign{
  E_{2n}&={(-1)^n\over2^{2n}} \,\sum_{\nu=1}^n D\,\b0^{[2\nu+1]}\cdot
  {\de^{2\nu}\b0^{2n}\over(2\nu)!}\,,
  }
  \eql{eul2}
  $$
which is identical to (\peq{Shov2}) as the value of the particular central differential in
(\peq{eul2}) is standard, [\pref{Steffensen}] (see below). This result is given in
[\pref{BSSV}] but only by comparing the actual values.

Another conclusion can be drawn from the identity,
   $$
    \sech {D\over2}\,\b0^{[2\nu]}=D\,\b0^{[2\nu+1]}\,,
    \eql{D1}
   $$
by comparing with the expansion (\peq{geneul}) to give the sum rule for Euler numbers,
  $$\eqalign{
  \sum_{n=1}^s{E_{2n}\over(2n)!\,2^{2n}}\,D^{2n}\b0^{[2s]}&=D\,\b0^{[2s+1]}
  ={(2s)!\over 2^{2s}}\comb{-{1\over2}}{s}\cr&=
  (-1)^s\bigg({1.3.\ldots.(2s-1)\over2^s}\bigg)^2\,,
  }
  \eql{SR}
  $$
which might be novel. Given the central differentials, it could be used to compute the Euler
numbers, but there is no need.

Related manipulations can be performed  for the generalised Euler numbers,
(\peq{Shov4}), although the simple relation, (\peq{D1}), becomes more complicated. As
an example, in order to further illustrate the theorem, (\peq{D}), I look at the
$E^{(3)}_{2n}$ numbers. Instead of (\peq{effd}), $f(D)$ in (\peq{D}) is given by,
   $$\eqalign{
    f(D)&=\int dD\,\cosh{D\over2}\,\sech^3\,{D\over2}\cr
    &=2\tanh D/2\,,
    }
    \eql{effd2}
   $$
which results in,
  $$
    \sech^3 {D\over2}\,\b0^{[2\nu]}=2\tanh {D\over2}\,\b0^{[2\nu+1]}\,.
    \eql{D2}
   $$
The $\tanh$ could be expanded in order to evaluate the right--hand side but it is better to
apply (\peq{D}) for a second time by choosing,
  $$\eqalign{
    f(D)&=2\int dD\,\cosh{D\over2}\,\tanh{D\over2}\cr
    &=4\cosh D/2\,,
    }
    \eql{effd3}
   $$
yielding the successive equivalences,
  $$
    \sech^3 {D\over2}\,\b0^{[2\nu]}=
    2\tanh {D\over2}\,\b0^{[2\nu+1]}=4\cosh {D\over2}\,\b0^{[2\nu+2]}\,,
    \eql{D3}
   $$
the last of which can now easily be evaluated as a finite sum of central differentials,
  $$
  4\cosh {D\over2}\,\b0^{[2\nu+2]}
  =4\sum_{k=1}^{\nu+1}{D^{2k}\b0^{[2\nu+2]}\over2^{2k}(2k)!}\,,
  $$
generalising (\peq{D1}).

Spelling out the details, the resulting expression for the generalised Euler numbers follows
as (\peq{eul}),
  $$\eqalign{
  E^{(3)}_{2n}&={(-1)^n\over2^{2n}} \,\sum_{\nu=1}^n \sech^3 {D\over2}\,\b0^{[2\nu]}\cdot
  {\de^{2\nu}\b0^{2n}\over(2\nu)!}\cr
  &={(-1)^n\over2^{2n-2}} \,\sum_{\nu=1}^n
  \sum_{k=1}^{\nu+1}{D^{2k}\b0^{[2\nu+2]}\over2^{2k}(2k)!}\cdot
  {\de^{2\nu}\b0^{2n}\over(2\nu)!}\,,
  }
  \eql{eul3}
  $$
which is considerably more complicated than the direct expansion, (\peq{Shov4}), but
allows us, therefrom, to infer the (even) sum rules,
  $$
   \sum_{k=1}^{\nu+1}{D^{2k}\b0^{[2\nu+2]}\over2^{2k-2}(2k)!}
   ={(2\nu)!\over2^{2\nu}}\comb{-{3\over2}}\nu\,,
   \eql{SR1}
  $$
which can be confirmed numerically.

As another example, I give the (odd) sum rule that results from looking at $E^{(5)}_{2n}$,
  $$
   {4\over3}\sum_{k=0}^{\nu+1}{D^{2k+1}\b0^{[2\nu+3]}\over(2k+1)!}
   ={(2\nu)!\over2^{2\nu}}\comb{-{5\over2}}\nu\,,
   \eql{SR2}
  $$
the formula corresponding to (\peq{D3}) being, after three integrations,
  $$
    \sech^5 {D\over2}\,\b0^{[2\nu]}={4\over3}\sinh {D}\,\b0^{[2\nu+3]}\,.
    \eql{D4}
   $$

The general case depends on integrating the relation,\footnote{ Constants of integration
give zero.}
  $$
    \sech^{2r+1}x={(\pm1)^{r+1}\over r\,(2r-1)!!}\bigg({d\over d\sinh x}
    \bigg)^{r+1}e^{\pm rx}
   \,\,,\quad r=1,2,\ldots\,,
   \eql{ident3}
  $$
which I cannot find listed. I need not write out the corresponding sum rules as examples
(\peq{SR1}) and (\peq{SR2}) are sufficient.

It is clear that there are a large number of sum rules derived by using the different
transformations, such as (\peq{D3}) and (\peq{D4}), arising from (\peq{D}) and by
computing higher Euler numbers, $E^{(2r+1)}_{2n}$. They would seem to have only a
curiosity or checking value. I have not seen them elsewhere.

With (\peq{ident3}) I have essentially returned to the direct expansion, (\peq{Shov4}),
since the multiple--angle $\cosh$ and $\sinh$ can be expanded in powers of $\sinh x$ in a
familiar way, the coefficients usually being derived by recursion or by schoolbook
trigonometry (\eg\ Loney, [\pref{Loney}], Bromwich, [\pref{Bromwich}], Workman,
[\pref{Workman}], Schwatt, [\pref{Schwatt}] p.124).

The expansions can actually be combined neatly to the unfamiliar form,
  $$\eqalign{
   e^{rx}&=\sum_{s=0}^\infty 2^{s}{(r/2)^{[s]}\over s!}\,\sinh ^{s}x\,,
   \quad r\in\oZ\,,\cr
   }
   \eql{lonexp}
  $$
where I have written the coefficients in terms of central factorials, for future use in section
8. (See [\pref{BSSV}], \S3.1.)
  
Eq.(\peq{ident3}), after some brief algebra, coincides with the direct expansion, (\cf\
(\peq{Shov4})),\footnote{ This could be considered a proof of (\peq{ident3}), given the
expansions, (\peq{lonexp}), or these could be derived, by continuation.}
  $$
   \sech^{2r+1}x=\sum_{s=0}^\infty
   \comb{-r-{1\over2}}s\sinh^{2s}x\,,\quad \sinh x<1\,,
  $$
\cf also Schwatt, [\pref{Schwatt}] p.70.\footnote{ Schwatt also refers to Shovelton.}
Questions of convergence have been mostly ignored, and, in any case, for the finite series
involved in (\peq{Shov4}) are irrelevant.

\section{\bf6. Central difference theorems}

I develope some  basic concepts in central difference calculus using a mainly symbolic
approach. I do this by first proving the theorem (\peq{D}) and then passing on to other
notions and some expansion formulae. Of the few discussions of this topic, that by Riordan,
[\pref{Riordan}], is most valuable. The more recent review by Butzer {\it et al},
[\pref{BSSV}], is also useful.

The proof of (\peq{D}) is not obvious so I expose it here, following the discussion by Jeffery,
[\pref{Jeffery}], for forward difference calculus. For comparison I write out the theorems as
given by him,
  $$\eqalign{
  &A:\quad f(\De)\,\b0^n=n{f(\De)\over D}\,\b0^{n-1}=\ldots
  =n!{f(\De)\over D^n}1\cr
  &B:\quad f(\De)\,\b0^n={df(\De)\over dD}\,\b0^{n-1}=\ldots
  ={d^n f(\De)\over dD^n}1\cr
  &C:\quad f(D)\,\b0^{(n)}=n{f(D)\over \De}\,\b0^{(n-1)}=\ldots
  =n!{ f(D)\over \De^n}1\cr
  &D:\quad f(D)\,\b0^{(n)}={df(D)\over d\De}\b0^{(n-1)}=\ldots
  ={d^n f(D)\over d\De^n}1\,.\cr
  }
  \eql{JABCD}
  $$

I will not prove them but only say that $A$ and $C$ express the fact that the symbols $\De$
and $D$ equally commute when acting upon $x$ at $x=0$ as they do when $x$ is current.
\footnote{ As noted by Jeffery, $A$ and $B$ appear as Misc. Ex. (1) on p.242 of Boole,
[\pref{Boole}]. Since this was a Senate House Examination problem, it seems likely it was
due to Herschel.}

In proving $D$, Jeffery makes essential use of the operator $\bze\equiv\log(1+D)$  which is
a sort of dual, or functional inverse to the relation $D=e^\De-1$ in the sense that,
  $$
   e^{\log(1+D)}-1=D\,.
  $$

I start from the very basic property of differences as the action on factorials,
  $$
  F(\de)\,\b0^{[r]}=F(D)\,\b0^r\,,
  \eql{factact}
  $$
expressed in terms of central quantities for which a basic relation is,
   $$
     \de=2\sinh {D\over2}\,,\quad i.e.\quad D=2\sinh^{-1} {\de\over2}\,.
     \eql{cd}
   $$
I can therefore define a new function $f$ by the functional relation,
  $$
  F(\de)=F(2\sinh D/2)\equiv f(D)=f(2\sinh^{-1} \de/2)\,,
  $$
and just replace $F(\de)$ by $f(D)$ on the left--hand side of (\peq{factact}) while, on the
right--hand side, replace $\de$ by $D$ in $f$, {\it regarded as a function of} $\de$. So,
corresponding to $\bze$, if I define the operator $\bxi$, by
   $$
    \bxi=2\sinh^{-1}{D\over2}\,.
    \eql{xidef}
   $$
I obtain the handy reciprocal result,\footnote{Jeffery derives his formula in a different way,
but I find the present route easier.}
  $$
  f(D)\,\b0^{[r]}=f(\bxi)\,\b0^r\,.
  \eql{factact2}
  $$
I have not seen this approach in the literature.

The proof of the central version of Theorem $D$ now follows more or less directly. Applying
Theorem $B$, which holds generally, on appropriate interpretation of $f(\De)$, one finds,
  $$\eqalign{
  f(D)\,\b0^{[r]}&=f(\bxi)\,\b0^r={f'(\bxi)\over (1+D^2/4)^{1/2}}\,\b0^{r-1}\cr
  &=f'(\bxi)\bigg[1+\comb{-{1\over2}}1{D^2\over2^2}
  +\comb{-{1\over2}}2{D^4\over2^4}+\ldots\bigg]{\b0}^{r-1}\cr
 &=f'(\bxi)\bigg[\b0^{r-1}+\comb{-{1\over2}}1{\b0^{r-3}\over2^2}
  +\comb{-{1\over2}}2{\b0^{r-5}\over2^4}+\ldots\bigg]\cr
  &=f'(D)\bigg[\b0^{[r-1]}+\comb{-{1\over2}}1{\b0^{[r-3]}\over2^2}
  +\comb{-{1\over2}}2{\b0^{[r-5]}\over2^4}+\ldots\bigg]\cr
  &={f'(D)\over(1+\de^2/4)^{1/2}}\,\b0^{[r-1]}={df(D)\over d\de}\,\b0^{[r-1]}\,.
  }
  \eql{thmD}
  $$
This is equation (\peq{D}).

Continuing with some central formalism, I display the expansion of the central factorial,
  $$\eqalign{
  x^{[n]}&=\sum_{\nu=1}^n x^{\nu}\, {D^{\nu}\b0^{[n]}\over\nu!}=
  \sum_{\nu=1}^n x^{\nu}\, {\bxi^{\nu}\b0^{n}\over\nu!}\cr
  &=e^{x\,\bxi}\,\b0^n\,,
  }
  \eql{cdf}
  $$
using (\peq{factact2}) to rewrite the coefficients and also formally summing the series
($1.\b0^n=0$ for $n\ne0$). (Compare Jeffery, [\pref{Jeffery}], \S5.)

Furthermore, from the explicit definition of $\bxi$, (\peq{xidef}), the expansion of the
inverse $\sinh$ function is contained in,
  $$
   \big(2\sinh^{-1}{x\over2}\big)^n=\bxi^n\,e^{\b0\,x}\,,
   \eql{invsh}
  $$
showing rapidly that the coefficients are just the central factorial numbers of the first kind
(\cf\ [\pref{BSSV}] \S4.2 for a more expansive treatment.)

The same expansion is derived by Butzer {\it et al}, [\pref{BSSV}] \S4.1, but less
attractively in my view as the proof involves manipulating several series. See also
Charalambides, [\pref{Charal}] and the classic Schwatt, [\pref{Schwatt}], pp.6,7,125, who
credits Cauchy.

As a further illustration of the symbolic formalism, I derive an intrinsic relation between the
first and second kind central factorial numbers, \ie between the $\bxi^m\,\b0^n$ and
$\de^r\,\b0^s$.

I start from the definition of the second kind,
  $$\eqalign{
  x^{n}=\sum_{\nu=1}^n x^{[\nu]}\, {\de^{\nu}\b0^{n}\over\nu!}
  =\sum_{\nu=1}^n e^{x\bxi}\,\b0^\nu\, {\de^{\nu}\b0^{n}\over\nu!}\,,\cr
  }
  $$
and equate powers of $x$ to get the compact statement of `orthogonality',
  $$
   \bxi^m\,e^{\b0\,\de}\,\b0^n=n!\,\de_m^n\,.
   \eql{inv}
  $$
The $\de$ on the right--hand side is a Kronecker delta. Equation (\peq{inv}) is really a
statement about inverses, and tantamount to completeness (of polynomial bases), the
reciprocal relation being,
  $$
   \de^n\,e^{\b0\,\bxi}\,\b0^m=m!\,\de_n^m\,.
   \eql{inv2}
  $$
\cf\ Riordan, [\pref{Riordan}] p.213.
  
  A more formal expression of completeness is,
  $$\eqalign{
   e^{\b0.\de}e^{\b0.\bxi}&=1_\xi\cr
   e^{\b0.\bxi}e^{\b0.\de}&=1_\de\,,
   }
   \eql{comp}
  $$
with,
  $$
   f(\bxi)\,1_\xi=f(\bxi)\,,\quad  f(\de)\,1_\de=f(\de)\,.
  $$

Two further useful relations are
  $$
    f(D)\,\b0^r=f(\bxi)\,e^{\b0.\de}\,\b0^r
    \eql{cor}
  $$
and
  $$
    f(\bxi)=f(D)\,e^{\b0.\bxi}\,.
    \eql{delt}
  $$
The first is a corollary of (\peq{inv}) and the second a trivial consequence of the basic
derivative, $D^m\,\b0^n=m!\,\de_m^n$.

Corresponding to (\peq{inv}) and (\peq{inv2}) there are two important expansions
connecting multiple derivatives and multiple differences. I derive these for central
differences following the forward difference treatment of Boole, [\pref{Boole}] p.24, which
is an application of Maclaurin's theorem in its secondary form,
  $$
  \phi(t)=\phi(D)\,e^{\b0.t}\,,
  \eql{mt2}
  $$
where $D$ acts on $\b0$. From the identities, (\peq{cd}), $\de$ is a function of $D$ and so,
setting $t=D$ and $\phi=\de^n$,
  $$
  \de^n=\de^n\,e^{\b0.D}\,,
  \eql{mcd}
  $$
where the $\de$ on the right--hand side acts on $\b0$. This is the required expression in
symbolic form. Expansion of the exponential gives a sum of powers of $D$ and, acting on a
function of $x$, yields the concrete representation,\footnote{ I ignore questions of
convergence and assume that the function, $f$, possesses the qualities required for all the
equations to exist.}
  $$
  \de^n u(x)=\sum_{m=n}^\infty {\de^n\b0^m\over m!}\,D^m u(x)\,.
  \eql{mcdc}
  $$

Inversely, consider $\phi=D^m$ as a function of $\de$, and so set $t=\de$. According to
(\peq{mt2}), $\de$ in $D$ then has to be replaced by $D$. This gives the $\bxi$ operator,
(\peq{xidef}), and so one has the required symbolic expression,
  $$
   D^m=\bxi^m\,e^{\b0.\de}\,,
   \eql{mcx}
  $$
yielding the concrete, Newton series form,
  $$
  D^m u(x)=\sum_{n=m}^\infty {\bxi^m\b0^n\over n!}\,\,\de^m u(x)\,,
  $$
which could be derived from (\peq{mcdc}) using (\peq{inv}). These expansions are
standard, see, \eg, [\pref{BSSV}] \S6.1.

I note in passing that the explicit forms of the operators, (\peq{cd}), (\peq{xidef}), are not
required when deriving the symbolic equations. Writing the relation, (\peq{cd}), as
$\de=\Psi(D)$, $\bxi$ is defined by the functional inverse,
$\bxi=\Psi^{-1}(D)=\Psi^{-2}(\de)$. A choice of $\Psi(*)$ other than (\peq{xidef}), or
$\log(1+*)$, would entail a factorial different from the usual ones.
\section{\bf7. Expansions}
As Jeffery points out, knowledge of the factorial numbers is useful in many expansions. For
example,
  $$\eqalign{
  f\bigg(2\sinh^{-1}\big(\sinh^{-1}{x\over2}\big)\bigg)&=
  f\bigg(2\sinh^{-1}{\bxi\over2}\bigg)\,e^{\b0.x}\cr
  &=f\bigg(2\sinh^{-1}{D\over2}\bigg)\,e^{\b0.\bxi}\,e^{\b0.x}\cr
  &=f(\bxi)\,e^{\b0.\bxi}\,e^{\b0.x}\cr
  &=f(\bxi)\,\bigg(\b0^{[1]}{x\over1!}+\b0^{[2]}\,{x^2\over2!}\ldots\bigg)\,,
  }
  \eql{exp2}
  $$
where (\peq{delt}) and (\peq{cdf}) have been used to give the second and fourth lines. This
process can be continued.

As a calculated example, consider the simplest case of (\peq{exp2}) which is
   $$\eqalign{
  2\sinh^{-1}\big(\sinh^{-1}{x\over2}\big)
  =\bxi\,\b0^{[1]}\cdot{x\over1!}+\bxi\,\b0^{[2]}\cdot{x^2\over2!}\ldots\,.
  }
  \eql{exp3}
  $$
and require the coefficient of $x^r/r!$,
  $$
  \bxi\,\b0^{[r]}=\sum_{s=0}^r{\bxi\,\b0^s\over s!}\cdot\bxi^s\,\b0^r\,.
  $$
(This is really going back to the third line of  (\peq{exp2}). A simple iteration of
(\peq{invsh}) also yields the same result.) The summation index $s$ must be odd, and so
therefore must $r$ (as one knows). Hence the coefficient of $x^{2k+1}/(2k+1)!$ is
  $$\eqalign{
  \bxi\,\b0^{[2k+1]}&=\sum_{j=0}^k{\bxi\,\b0^{2j+1}\over (2j+1)!}\cdot
  \bxi^{2j+1}\,\b0^{2k+1}\cr
  &=\sum_{j=0}^k\comb{-{1\over2}}{j}{1\over  2^{2j}(2j+1)}\,,
  \bxi^{2j+1}\,\b0^{2k+1}\cr
  }
  $$
which can be evaluated straightforwardly.

\section{\bf8. Using Representative Notation}

A simpler case to (\peq{exp2}) follows either as an extension of (\peq{invsh}) or directly,
  $$
  f\big(2\sinh^{-1}{x\over2}\big)=f(\bxi)\,e^{\b0.x}=\sum_{r=0}^\infty\,
  {f(\bxi)\b0^r\over r!}\,x^r\,,
  \eql{fxi}
  $$
and allows me to bring in Blissard's seductive `Representative Notation'  (or classical umbral
calculus, [\pref{RandT2}]) by writing the primary form of Maclaurin's theorem as,
  $$
  f(\bxi)=\sum_{s=0}^\infty{f^{(s)}(0)\over s!}\,\bxi^s=e^{f\,\bxi}\,,
  $$
where the power $f^n$ is to be replaced by the derivative, $f^{(n)}(0)$, after the
expansion of the exponential. Then the elegant result,
  $$\eqalign{
  f\big(2\sinh^{-1}{x\over2}\big)&=\sum_{r=0}^\infty
  {e^{f\,\bxi}\b0^r\over r!}\,x^r\cr
  &=\sum_{r=0}^\infty{f^{[r]}\over r!}\,x^r\,,\cr
  }
  $$
follows on using (\peq{cdf}). One could write this result more enigmatically as,
  $$
 f(\bxi)\,e^{\b0.x} =e^{fx}\,,
 \eql{snap}
  $$
where now the power, $f^r$, first represents the central factorial, $f^{[r]}$, in which the
powers of $f$ are still to be replaced by  the derivatives, $f^{(n)}(0)$, as above.

In the form,
  $$\eqalign{
  f(y)=\sum_{r=0}^\infty{f^{[r]}\over r!}\,\big(2\sinh{y\over2}\big)^r\,,\cr
  }
  $$
this expansion is derived, in a complicated way, by Teixeira, [\pref{Teixeira}], using {\it
existing} trigonometric expansions, \cf\ (\peq{lonexp}). He sets $y=i\th$ and splits the
formula into real and imaginary parts. The derivation here, in contrast, is rapid and from
finite difference fundamentals, which is how the trigonometric expansions can be, and are,
more logically obtained, [\pref{JandS}]. A more elaborate, conventional discussion can be
found in [\pref{BSSV}] \S5.

Riordan employs the umbral `method' and I use his notation, (see Eq.(\peq{gees})), for the
factorial numbers, which has operational advantages). My object is to combine umbral
techniques with the Jeffery operator formalism (in central form).

I have already introduced the representative, $t$, by,
   $$
    t^n(x)=t_n(x)\equiv\sum_{\nu=0}^n t(n,\nu)\,x^\nu=x^{[n]}=e^{x\bxi}\,\b0^n\,,
    \eql{tee}
   $$
so that, by definition,
  $$
    x^n=e^{t\de}\,\b0^n\,,\quad t=t(x)\,.
    \eql{def2}
  $$

I now define its reciprocal, $T$, by, \cf\ [\pref{Riordan}],
  $$
    T^n(x)=T_n(x)\equiv\sum_{\nu=0}^n T(n,\nu)\,x^\nu=e^{x\de}\,\b0^n\,.
    \eql{Tee}
   $$
On the right I have given the Jeffery central equivalents which enable me to obtain, by
summation, the umbral expressions,
  $$
e^{yt}=e^{x\bxi}\,e^{\b0 y}
 \eql{cons}
 $$
 $$
  e^{yT}=e^{x\de}\,e^{\b0 y}\,,
  \eql{cons2}
  $$
where $t=t(x)$ and $T=T(x)$.

From the definitions of the operators $\bxi$ and $\de$, the generating functions,
(\peq{cons}) and (\peq{cons2}), take the explicit forms, \cf\ (\peq{fxi}),
  $$
    e^{yt(x)}= e^{2x\sinh^{-1}y/2}\,,
     \eql{Bl4}
  $$
and
  $$
     e^{yT(x)}=e^{2x\sinh y/2}\,,
     \eql{Bl3}
  $$
both of which Riordan (see also Carlitz and Riordan, [\pref{CaandR}]) derives by series
combinatorial manipulation but which here follow almost painlessly.

The two expressions are related by reciprocity as can be shown in the following way.

I first derive the umbral statement of reciprocity. Umbrally, from  Eq.(\peq{cons2}) by
setting $x\to t(x)$ one arrives at the formally neat expression of
orthogonality/completeness,
   $$
    T\big(t(*)\big)=*\,.
    \eql{recip}
   $$
(This also easily results from (\peq{Tee}) followed by (\peq{tee}), which gives
$T^n\big(t(x)\big)=x^n$.)

The reciprocal relation follows on iteration of (\peq{recip}),
  $$
    t\big(T\big(t(*)\big)\big)=t(*)\,,
    \eql{recip2}
   $$
or, invoking completeness,
  $$
    t\big(T(*)\big)=*\,.
    \eql{recip3}
   $$
One might therefore formally set $T=t^{-1}$, \etc Eqs. (\peq{recip}) and (\peq{recip3})
are the umbral equivalents of the explicit (\peq{inv}) and (\peq{inv2}).
 
Now, making the replacement $x\to T(x)$ and employing (\peq{recip3}) turns (\peq{Bl4})
into (\peq{Bl3}) as promised.

I emphasise the rather interesting fact that Eq.(\peq{Bl4}) is identical with the classical
trigonometric expansion,  (\peq{lonexp}), which was expressed in a relevant way for this
reason. Adjusting notation, (\peq{Bl4}), under $y\to2\sinh x$ and $x\to r/2$, becomes,
neatly,
  $$
     e^{rx}=e^{2t\sinh x}\,,\quad t=t(r/2)\,,
     \eql{Bl1}
  $$
which is the umbral representation of (\peq{lonexp}). A trigonometric expansion has
therefore been obtained with very little effort. \footnote{ Eq.(\peq{lonexp}) is to be found
as Eq.(28) on p.215 of  [\pref{Riordan}] derived in a combinatorial way. (Riordan refers to
P\'olya and Szeg\"o.) The approach could also be considered as an alternative derivation of
the expansions (\peq{lonexp}). Riordan also exhibits (\peq{Bl1}). (To make the connection,
set $x\to z/2$, $r\to 2x$ and $t(r/2)\to t(x)$.)}

Finally I make a formal functional generalisation in accordance with the remarks at the end
of section 6. If $\bxi=\Psi^{-1}(D)$, the generating functions are
   $$
     e^{yt(x)}=e^{x\Psi^{-1}(y)}\,,
     \eql{Bl3}
  $$
and
  $$
     e^{yT(x)}=e^{x\Psi(y)}\,,
     \eql{Bl3}
  $$
with $t$ and $T$ now standing for mutually inverse generic `factorial' numbers.
 
\section{\bf9. Conclusion}
The compactness and rapidity of the central difference symbolic calculus has, I think, been
demonstrated. I hope also to have shown that there is a considerable body of ignored
historical work that could usefully be exploited. Whether the particular results are
interesting is, of course, a matter of personal taste. Certainly, the symbolic approach has an
intrinsic elegance and power. I draw attention to the convenience of the operator, $\bxi$. As
shown in the last section it allows the double generating functions of Carlitz and Riordan to
be obtained with ease and thence some classic trigonometric expansions, likewise.
\newpage
\vglue 15truept

\noin{\bf References.} \vskip5truept

\begin{putreferences}
   \ref{BaandB}{Balian,R. and Bloch,C. \aop{60}{1970}{401}.}
   \ref{Zhang}{Zhang,W.P. {\it The Fibonacci Quarterly} {\bf 36} (1998) 154.}
   \ref{LiandZ}{Liu,G.D. and Zhang,W.P. {\it Acta Math. Sinica} {\bf24} (2008) 343.}
   \ref{JandS}{Joliffe, A.E. and Shovelton, S.T. \plms{13}{1913}{29}.}
   \ref{CaandR}{Carlitz,L. and Riordan,J. \cjm {15}{1963}{94}.}
   \ref{RandT2}{Rota, G-.C. and Taylor, B.D. {\it SIAM J.Math.Anal.} {\bf 25} (1994) 694.}
   \ref{Bickley}{Bickley, W.G. {\it J. Math. and Phys.}(MIT) {\bf 27} (1948) 183.}
   \ref{Schwatt}{Schwatt, I.J. {\it An introduction to the operations with series}
    (Univ.Pennsylvania Press, Philadelphia, 1924).}
    \ref{Loney}{Loney, S.L. {\it Plane Trigonometry} (CUP, Cambridge, 1893).}
   \ref{Teixeira}{Teixeira, F.G. \jram{116}{1896}{14}.}
   \ref{Riordan}{Riordan,J. {\it Combinatorial Identities} (Wiley, New York, 1968).}
   \ref{Charal}{Charalambides, Ch.A. {\it The Fibonacci Quarterly} {\bf 19.5} (1981) 451.}
    \ref{Stern}{Stern,W. \jram {79}{1875}{67}.}
    \ref{Thiele}{Thiele,T.N. {\it Interpolationsrechnung} (Teubner, Leipzig, 1909).}
    \ref{Dowren}{Dowker,J.S. \jpamt {46}{2013}{2254}.}
    \ref{KnandB}{Knuth,D.E. and Buckholtz,T.J. {\it Math.Comp.} {\bf 21} (1967) 663.}
    \ref{BSSV}{Butzer,P.L., Schmidt,M., Stark,E.L. and Vogt,I. {\it Numer.Funct.Anal.Optim.}
    {\bf 10} (1989) 419.}
    \ref{Ely}{Ely,G.S. \ajm{5}{1882}{337}.}
    \ref{BrandH}{Brent,R.P. and Harvey,D., {\it Fast Computation of Bernoulli, Tangent and
    Secant numbers}, ArXiv:1108.0286.}
   \ref{BaandBe}{Bayad and Beck  ArXiv.}
   \ref{Blissard}{Blissard, \qjm{}{}{}.}
   \ref{Gould2}{Gould,H.W. \amm{79}{1972}{} .}
   \ref{Worpitzky}{Worpitsky,J. \jram{}{}{}.}
   \ref{Saalschutz}{Saalsch\"utz,L. {\it Vorlesungen \"uber die Bernoullischen Zahlen.}
   (Springer--Verlag, Berlin, 1893).}
   \ref{Workman}{Workman, W.P. {\it Memoranda Mathematica}, (OUP, Oxford, 1912).}
   \ref{Nielsen}{Nielsen,N. {\it Handbuch der Theorie der Gammafunktion} (Teubner,
   Leipzig, 1906).}
   \ref{Steffensen}{Steffensen,J.F. {\it Interpolation}, (Williams and Wilkins,
    Baltimore, 1927).}
    \ref{Joffe}{Joffe,S.A. \qjm{47}{1916}{103}.}
    \ref{Shovelton}{Shovelton,S.T. \qjm{46}{1915}{220}.}
   \ref{Horner}{Horner,J. \qjm{4}{1861}{111}.}
   \ref{Jeffery}{Jeffery, H.M. \qjm{4}{1861}{364}.}
   \ref{Jeffery2}{Jeffery, H.M. \qjm{5}{1862}{91}.}
   \ref{Jeffery3}{Jeffery, H.M. \qjm{6}{1864}{179}.}
   \ref{Grunert}{Grunert, \jram{25}{1843}{240}.}
   \ref{AgandD}{Agoh,T. and Dilcher,K. {\it J. Number Theory} {\bf 124} (2007) 105.}
   \ref{Franssens}{Franssens,G.R. {\it Analysis Mathematica} {\bf 33} (2007) 17.}
   \ref{Andrews}{Andrews,G.E. {\it Ramaujan J.} {\bf 7} (2003) 385.}
    \ref{Glaisher2}{Glaisher,J.W.L. \qjm {40}{1909}{275}.}
    \ref{RandF}{Rubinstein, B.Y. and Fel,L.G., {\it Ramanujan J.} {\bf11}(2006)331.}
    \ref{Rubinstein}{Rubinstein, B.Y. {\it Ramanujan J.} {\bf15}(2008)177.}
   \ref{Dickson3}{Dickson,L.E. {\it History of the Theory of Numbers} vol.2.
   (Carnegie Institute, Washington, 1920).}
   \ref{Brioschi}{Brioschi,F. {\it Annali di sc. mat. e fis.} {\bf8} (1857) 5.}
   \ref{Fedosov}{Fedosov,B.V. {\it Sov.Math.Dokl.} {\bf 5} (1963) 1992}
   \ref{Chapman}{Chapman}
    \ref{Jordan}{Jordan,C. {\it Calculus of Finite Differences}, (Budapest, 1939).}
   \ref{PandB}{van der Pol, B. and Bremmer,H. {\it Operational Calculus} (CUP,
   Cambridge, 1964).}
   \ref{TandD}{Tauber,S. and Dean,D. {\it J.SIAM} {\bf 8} (1960) 174.}
   \ref{Traub}{Traub,J.F. {\it Math. of Comp.} {\bf 19} (1965) 177.}
   \ref{Chapman}{Chapman,S. \plms{9}{1911}{369}. }
   \ref{Moore}{Moore,D.H. \amm{69}{1962}{132}.}
   \ref{HaandR}{Hardy,G.H. and Riesz,M. {\it General theory of Dirichlet's series}
   (CUP, Cambridge, 1915).}
   \ref{Gupta}{Gupta,H. {\it Tables of partitions} (Royal.Soc.Math.Tables 4)
   (Cambridge, 1958).}
   \ref{Bromwich}{Bromwich, T.J.I'A. {\it Infinite Series},
  (Macmillan, London, 1947).}
  \ref{Hobson2}{Hobson,E.W. {\it Theory of functions of a real variable}, vol.2.
  (C.U.P., Cambridge,1907).}
   \ref{Bromwich2}{Bromwich, T.J.I'A. {\it Infinite Series}, 1st Edn.
  (Macmillan, London,1907).}
   \ref{Knopp3}{Knopp,K. {\it Theory of Infinite Series} (Blackie, London, 1928).}
   \ref{Vandiver}{Vandiver,H.S. \tams{51}{1942}{502}.}
   \ref{Adamchik}{Adamchik,V.S.{\it Appl.Math. and Comp.} {\bf 187} (2007) 3.}
   \ref{Cvijovic}{Cvijovi\'{c},D. {\it Appl.Math. and Comp.} {\bf 215} (2009) 3002.}
   \ref{Cvijovic1}{Cvijovi\'{c},D. {\it Comp.and Math. with Appl.} {\bf 62} (2011) 1879.}
   \ref{Carlitz2}{Carlitz,L. {\it Math.Magazine} {\bf 32} (1959) 247.}
   \ref{Israilov}{Israilov.M.I. {\it Sibirsk Mat. Zh.} {\bf 22} (1981) 121.}
   \ref{SandZ}{Sills,A.V. and Zeilberger,D. {\it Formulae for the number of partitions of
    $n$ into at most $m$ parts (using the quasi-polynomial ansatz)} ArXiv:1108.4391.}
   \ref{BandV}{Brion,M. and Vergne,M. {\it J.Am.Math.Soc.} {\bf 10} (1997) 371.}
   \ref{Knopf}{Knopf,P.M. {\it Math. Magazine} {\bf 76} (2003) 364.}
   \ref{Gould}{Gould,H.W. {\it Am. Math. Monthly} {\bf 86} (1978) 450.}
   \ref{Hoffman}{Hoffman,M.E. \amm{102}{1995}{23}.}
   \ref{Hoffman2}{Hoffman,M.E. {\it Electronic J.Comb.} {\bf 6} (1999) \#R21.}
   \ref{Boyadzhiev}{Boyadzhiev, K.N. {\it Fibonacci Quarterly} {\bf45} (2008) 291.}
   \ref{Sills}{Sills}
   \ref{Munagi}{Munagi,A.O. {\it Electronic J. Comb. Number Th.} {\bf 7} (2007) \#A25.}
   \ref{Dowsyl}{Dowker,J.S. {\it Relations between Ehrhart polynomials, the heat kernel
   and \break Sylvester waves} ArXiv:1108.1760}
   \ref{Alfonsin}{Alfonsin,J.L.R. {\it The Diophantine Frobenius Problem} (O.U.P.,
   Oxford, 2005).}
   \ref{GeandS}{Gel'fand,I.M. and Shilov,G.E. {\it Generalized Functions} Vol.1
   (Academic Press, N.Y., 1964).}
    \ref{Boole}{Boole, G. {\it Calculus of Finite Differences}, (MacMillan, Cambridge,
1860).}
   \ref{Boole2}{Boole, G. {\it Calculus of Finite Differences}, 2nd Edn. (MacMillan, Cambridge,
1872).}
     \ref{Cayley5}{Cayley, A.  {\it Phil. Trans. Roy. Soc. Lond.} {\bf 148} (1858) 47.}
     \ref{Milne-Thomson}{Milne-Thomson, L.M. {\it The Calculus of Finite Differences},
     (MacMillan, London, 1933).}
    \ref{Herschel}{Herschel, J.F.W.  {\it Phil. Trans. Roy. Soc. Lond.} {\bf 106} (1816) 25.}
    \ref{Herschel2}{Herschel, J.F.W.  {\it Phil. Trans. Roy. Soc. Lond.} {\bf 108} (1818) 144.}
     \ref{Herschel4}{Herschel, J.F.W.  {\it Phil. Trans. Roy. Soc. Lond.} {\bf 140} (1850) 399.}
      \ref{Herschel5}{Herschel, J.F.W.  {\it Phil. Trans. Roy. Soc. Lond.} {\bf 150} (1860) 319.}
    \ref{Herschel3}{Herschel,J.F.W. {\it Examples of the Applications of the Calculus of Finite
    Differences} (Deighton, Cambridge, 1820).}
    \ref{Brinkley}{Brinkley {\it Phil. Trans. Roy. Soc. Lond.} {\bf } (1807) 114 .}
    \ref{Littlewood2}{Littlewood,D.E. {\it The Theory of Group Characters}
    (Clarendon Press, Oxford, 1950).}
    \ref{Wright}{Wright,E.M. \amm{68}{1961}{144}.}
\ref{Carlitz}{Carlitz,L. \dmj{27}{1960}{401}.}
     \ref{Netto}{Netto,E. {\it Lehrbuch der Combinatorik} 2nd Edn. (Teubner, Leipzig, 1927).}
    \ref{FdeB}{Fa\`{a} de Bruno, F. {\it Th\'eorie des Formes Binaires} (Brero, Turin,1876).}
    \ref{Ehrhart}{Ehrhart,E. \jram{227}{25}{1967}.}
    \ref{Bell}{Bell,E.T.\ajm{65}{1943}{382}.}
    \ref{BandR}{Beck,M. and Robins,S. {\it Computing the Continuous Discretely,}
    (Springer, New York, 2007).}
    \ref{BandR2}{Beck, M. and Robins,S. {\it Discrete and Comp. Geom.} {\bf 27}(2002) 443.}
    \ref{Harmer}{Harmer,M. {\it J.Australian Math.Soc.} {\bf 84}(2008)217.}
    \ref{RandF}{Rubinstein, B.Y. and Fel,L.G., {\it Ramanujan J.} {\bf11}(2006)331.}
    \ref{BGK}{Beck,M., Gessel, I.M. and Komatsu,T. {\it Electronic Journal of Combinatorics}
    {\bf8}(2001) 1.}
    \ref{Sylvester}{Sylvester,J.J. \qjpam{1}{1858}{81}.}
    \ref{Sylvester2}{Sylvester,J.J. \qjpam{1}{1858}{142}.}
    \ref{Sylvester3}{Sylvester,J.J. \ajm{5}{1882}{79}.}
    \ref{Sylvester4}{Sylvester,J.J. \plms{28}{1896}{33}.}
    \ref{Dowgta}{J.S.Dowker, {\it Group theory aspects of spectral problems on spherical
    factors}, ArXiv.Math.DG: 0907.1309.}
    \ref{BDR}{Beck,M., Diaz and Robins,S. {\it J.Numb.Theory} {\bf 96} (2002) 1.}
    \ref{PandS}{P\'{o}lya, G. and Szeg\H{o},G. {\it Aufgaben und Lehrs\"atze aus der Analysis}
    (Springer--Verlag, Berlin, 1925).}
    \ref{EOS}{Elizalde,E., Odintsov, S.D. and Saharian, A.A. \prD{79}{2009}{065023}.}
    \ref{Cavalcanti}{Cavalcanti,R.M. \prD{69}{2004}{065015}.}
    \ref{MWK}{Milton, K.A., Wagner,J. and Kirsten,K. \prD{80}{2009}{125028}.}
    \ref{EBM2}{Ellingsen,S.A., Brevik,I. and Milton,K.A. \prE{81}{2010}{065031}.}
    \ref{EBM}{Ellingsen,S.A., Brevik,I. and Milton,K.A. \prE{80}{2009}{021125}.}
    \ref{BEM}{Brevik,I., Ellingsen,S.A. and Milton,K.A. \prE{79}{2009}{041120}.}
    \ref{FKW}{Fulling,S.A, Kaplan L. and Wilson,J.H. \prA{76}{2007}{012118}.}
    \ref{Lukosz}{Lukosz,W, {\it Physica} {\bf 56} (1971) 109; \zfp{258}{1973}{99}
    ;\zfp{262}{1973}{327}.}
    \ref{Gromes}{Gromes, D. \mz{94}{1966}{110}.}
    \ref{FandK1}{Kirsten,K. and Fulling,S.A. \prD{79}{2009}{065019} .}
    \ref{FandK2}{Fucci,G. and Kirsten,K, JHEP (2011), 1103:016.}
    \ref{dowgjms}{Dowker,J.S. {\it Determinants and conformal anomalies
    of GJMS operators on spheres}, ArXiv: 1007.3865.}
    \ref{Dowcascone}{dowker,J.S. \prD{36}{1987}{3095}.}
    \ref{Dowcos}{dowker,J.S. \prD{36}{1987}{3742}.}
    \ref{Dowtherm}{Dowker,J.S. \prD{18}{1978}{1856}.}
    \ref{Dowgeo}{Dowker,J.S. \cqg{11}{1994}{L55}.}
    \ref{ApandD2}{Dowker,J.S. and Apps,J.S. \cqg{12}{1995}{1363}.}
   \ref{HandW}{Hertzberg,M.P. and Wilczek,F. {\it Some calculable contributions to
   Entanglement Entropy}, ArXiv:1007.0993.}
   \ref{KandB}{Kamela,M. and Burgess,C.P. \cjp{77}{1999}{85}.}
   \ref{Dowhyp}{Dowker,J.S. \jpa{43}{2010}{445402}; ArXiv:1007.3865.}
   \ref{LNST}{Lohmayer,R., Neuberger,H, Schwimmer,A. and Theisen,S.
   \plb{685}{2010}{222}.}
   \ref{Allen2}{Allen,B. PhD Thesis, University of Cambridge, 1984.}
   \ref{MyandS}{Myers,R.C. and Sinha,A. {\it Seeing a c-theorem with
   holography}, ArXiv:1006.1263}
   \ref{MyandS2}{Myers,R.C. and Sinha,A. {\it Holographic c-theorems in
   arbitrary dimensions},\break ArXiv: 1011.5819.}
   \ref{RyandT}{Ryu,S. and Takayanagi,T. JHEP {\bf 0608}(2006)045.}
   \ref{CaandH}{Casini,H. and Huerta,M. {\it Entanglement entropy
   for the n--sphere},\break arXiv:1007.1813.}
   \ref{CaandH3}{Casini,H. and Huerta,M. \jpa {42}{2009}{504007}.}
   \ref{Solodukhin}{Solodukhin,S.N. \plb{665}{2008}{305}.}
   \ref{Solodukhin2}{Solodukhin,S.N. \plb{693}{2010}{605}.}
   \ref{CaandW}{Callan,C.G. and Wilczek,F. \plb{333}{1994}{55}.}
   \ref{FandS1}{Fursaev,D.V. and Solodukhin,S.N. \plb{365}{1996}{51}.}
   \ref{FandS2}{Fursaev,D.V. and Solodukhin,S.N. \prD{52}{1995}{2133}.}
   \ref{Fursaev}{Fursaev,D.V. \plb{334}{1994}{53}.}
   \ref{Donnelly2}{Donnelly,H. \ma{224}{1976}{161}.}
   \ref{ApandD}{Apps,J.S. and Dowker,J.S. \cqg{15}{1998}{1121}.}
   \ref{FandM}{Fursaev,D.V. and Miele,G. \prD{49}{1994}{987}.}
   \ref{Dowker2}{Dowker,J.S.\cqg{11}{1994}{L137}.}
   \ref{Dowker1}{Dowker,J.S.\prD{50}{1994}{6369}.}
   \ref{FNT}{Fujita,M.,Nishioka,T. and Takayanagi,T. JHEP {\bf 0809}
   (2008) 016.}
   \ref{Hund}{Hund,F. \zfp{51}{1928}{1}.}
   \ref{Elert}{Elert,W. \zfp {51}{1928}{8}.}
   \ref{Poole2}{Poole,E.G.C. \qjm{3}{1932}{183}.}
   \ref{Bellon}{Bellon,M.P. {\it On the icosahedron: from two to three
   dimensions}, arXiv:0705.3241.}
   \ref{Bellon2}{Bellon,M.P. \cqg{23}{2006}{7029}.}
   \ref{McLellan}{McLellan,A,G. \jpc{7}{1974}{3326}.}
   \ref{Boiteaux}{Boiteaux, M. \jmp{23}{1982}{1311}.}
   \ref{HHandK}{Hage Hassan,M. and Kibler,M. {\it On Hurwitz
   transformations} in {Le probl\`eme de factorisation de Hurwitz}, Eds.,
   A.Ronveaux and D.Lambert (Fac.Univ.N.D. de la Paix, Namur, 1991),
   pp.1-29.}
   \ref{Weeks2}{Weeks,Jeffrey \cqg{23}{2006}{6971}.}
   \ref{LandW}{Lachi\`eze-Rey,M. and Weeks,Jeffrey, {\it Orbifold construction of
   the modes on the Poincar\'e dodecahedral space}, arXiv:0801.4232.}
   \ref{Cayley4}{Cayley,A. \qjpam{58}{1879}{280}.}
   \ref{JMS}{Jari\'c,M.V., Michel,L. and Sharp,R.T. {\it J.Physique}
   {\bf 45} (1984) 1. }
   \ref{AandB}{Altmann,S.L. and Bradley,C.J.  {\it Phil. Trans. Roy. Soc. Lond.}
   {\bf 255} (1963) 199.}
   \ref{CandP}{Cummins,C.J. and Patera,J. \jmp{29}{1988}{1736}.}
   \ref{Sloane}{Sloane,N.J.A. \amm{84}{1977}{82}.}
   \ref{Gordan2}{Gordan,P. \ma{12}{1877}{147}.}
   \ref{DandSh}{Desmier,P.E. and Sharp,R.T. \jmp{20}{1979}{74}.}
   \ref{Kramer}{Kramer,P., \jpa{38}{2005}{3517}.}
   \ref{Klein2}{Klein, F.\ma{9}{1875}{183}.}
   \ref{Hodgkinson}{Hodgkinson,J. \jlms{10}{1935}{221}.}
   \ref{ZandD}{Zheng,Y. and Doerschuk, P.C. {\it Acta Cryst.} {\bf A52}
   (1996) 221.}
   \ref{EPM}{Elcoro,L., Perez--Mato,J.M. and Madariaga,G.
   {\it Acta Cryst.} {\bf A50} (1994) 182.}
    \ref{PSW2}{Prandl,W., Schiebel,P. and Wulf,K.
   {\it Acta Cryst.} {\bf A52} (1999) 171.}
    \ref{FCD}{Fan,P--D., Chen,J--Q. and Draayer,J.P.
   {\it Acta Cryst.} {\bf A55} (1999) 871.}
   \ref{FCD2}{Fan,P--D., Chen,J--Q. and Draayer,J.P.
   {\it Acta Cryst.} {\bf A55} (1999) 1049.}
   \ref{Honl}{H\"onl,H. \zfp{89}{1934}{244}.}
   \ref{PSW}{Patera,J., Sharp,R.T. and Winternitz,P. \jmp{19}{1978}{2362}.}
   \ref{LandH}{Lohe,M.A. and Hurst,C.A. \jmp{12}{1971}{1882}.}
   \ref{RandSA}{Ronveaux,A. and Saint-Aubin,Y. \jmp{24}{1983}{1037}.}
   \ref{JandDeV}{Jonker,J.E. and De Vries,E. \npa{105}{1967}{621}.}
   \ref{Rowe}{Rowe, E.G.Peter. \jmp{19}{1978}{1962}.}
   \ref{KNR}{Kibler,M., N\'egadi,T. and Ronveaux,A. {\it The Kustaanheimo-Stiefel
   transformation and certain special functions} \lnm{1171}{1985}{497}.}
   \ref{GLP}{Gilkey,P.B., Leahy,J.V. and Park,J-H, \jpa{29}{1996}{5645}.}
   \ref{Kohler}{K\"ohler,K.: Equivariant Reidemeister torsion on
   symmetric spaces. Math.Ann. {\bf 307}, 57-69 (1997)}
   \ref{Kohler2}{K\"ohler,K.: Equivariant analytic torsion on ${\bf P^nC}$.
   Math.Ann.{\bf 297}, 553-565 (1993) }
   \ref{Kohler3}{K\"ohler,K.: Holomorphic analytic torsion on Hermitian
   symmetric spaces. J.Reine Angew.Math. {\bf 460}, 93-116 (1995)}
   \ref{Zagierzf}{Zagier,D. {\it Zetafunktionen und Quadratische
   K\"orper}, (Springer--Verlag, Berlin, 1981).}
   \ref{Stek}{Stekholschkik,R. {\it Notes on Coxeter transformations and the McKay
   correspondence.} (Springer, Berlin, 2008).}
   \ref{Pesce}{Pesce,H. \cmh {71}{1996}{243}.}
   \ref{Pesce2}{Pesce,H. {\it Contemp. Math} {\bf 173} (1994) 231.}
   \ref{Sutton}{Sutton,C.J. {\it Equivariant isospectrality
   and isospectral deformations on spherical orbifolds}, ArXiv:math/0608567.}
   \ref{Sunada}{Sunada,T. \aom{121}{1985}{169}.}
   \ref{GoandM}{Gornet,R, and McGowan,J. {\it J.Comp. and Math.}
   {\bf 9} (2006) 270.}
   \ref{Suter}{Suter,R. {\it Manusc.Math.} {\bf 122} (2007) 1-21.}
   \ref{Lomont}{Lomont,J.S. {\it Applications of finite groups} (Academic
   Press, New York, 1959).}
   \ref{DandCh2}{Dowker,J.S. and Chang,Peter {\it Analytic torsion on
   spherical factors and tessellations}, arXiv:math.DG/0904.0744 .}
   \ref{Mackey}{Mackey,G. {\it Induced representations}
   (Benjamin, New York, 1968).}
   \ref{Koca}{Koca, {\it Turkish J.Physics}.}
   \ref{Brylinski}{Brylinski, J-L., {\it A correspondence dual to McKay's}
    ArXiv alg-geom/9612003.}
   \ref{Rossmann}{Rossman,W. {\it McKay's correspondence
   and characters of finite subgroups of\break SU(2)} {\it Progress in Math.}
      Birkhauser  (to appear) .}
   \ref{JandL}{James, G. and Liebeck, M. {\it Representations and
   characters of groups} (CUP, Cambridge, 2001).}
   \ref{IandR}{Ito,Y. and Reid,M. {\it The Mckay correspondence for finite
   subgroups of SL(3,C)} Higher dimensional varieties, (Trento 1994),
   221-240, (Berlin, de Gruyter 1996).}
   \ref{BandF}{Bauer,W. and Furutani, K. {\it J.Geom. and Phys.} {\bf
   58} (2008) 64.}
   \ref{Luck}{L\"uck,W. \jdg{37}{1993}{263}.}
   \ref{LandR}{Lott,J. and Rothenberg,M. \jdg{34}{1991}{431}.}
   \ref{DoandKi} {Dowker.J.S. and Kirsten, K. {\it Analysis and Appl.}
   {\bf 3} (2005) 45.}
   \ref{dowtess1}{Dowker,J.S. \cqg{23}{2006}{1}.}
   \ref{dowtess2}{Dowker,J.S. {\it J.Geom. and Phys.} {\bf 57} (2007) 1505.}
   \ref{MHS}{De Melo,T., Hartmann,L. and Spreafico,M. {\it Reidemeister
   Torsion and analytic torsion of discs}, ArXiv:0811.3196.}
   \ref{Vertman}{Vertman, B. {\it Analytic Torsion of a  bounded
   generalized cone}, ArXiv:0808.0449.}
   \ref{WandY} {Weng,L. and You,Y., {\it Int.J. of Math.}{\bf 7} (1996)
   109.}
   \ref{ScandT}{Schwartz, A.S. and Tyupkin,Yu.S. \np{242}{1984}{436}.}
   \ref{AAR}{Andrews, G.E., Askey,R. and Roy,R. {\it Special functions}
  (CUP, Cambridge, 1999).}
   \ref{Tsuchiya}{Tsuchiya, N.: R-torsion and analytic torsion for spherical
   Clifford-Klein manifolds.: J. Fac.Sci., Tokyo Univ. Sect.1 A, Math.
   {\bf 23}, 289-295 (1976).}
   \ref{Tsuchiya2}{Tsuchiya, N. J. Fac.Sci., Tokyo Univ. Sect.1 A, Math.
   {\bf 23}, 289-295 (1976).}
  \ref{Lerch}{Lerch,M. \am{11}{1887}{19}.}
  \ref{Lerch2}{Lerch,M. \am{29}{1905}{333}.}
  \ref{TandS}{Threlfall, W. and Seifert, H. \ma{104}{1930}{1}.}
  \ref{RandS}{Ray, D.B., and Singer, I. \aim{7}{1971}{145}.}
  \ref{RandS2}{Ray, D.B., and Singer, I. {\it Proc.Symp.Pure Math.}
  {\bf 23} (1973) 167.}
  \ref{Jensen}{Jensen,J.L.W.V. \aom{17}{1915-1916}{124}.}
  \ref{Rosenberg}{Rosenberg, S. {\it The Laplacian on a Riemannian Manifold}
  (CUP, Cambridge, 1997).}
  \ref{Nando2}{Nash, C. and O'Connor, D-J. {\it Int.J.Mod.Phys.}
  {\bf A10} (1995) 1779.}
  \ref{Fock}{Fock,V. \zfp{98}{1935}{145}.}
  \ref{Levy}{Levy,M. \prs {204}{1950}{145}.}
  \ref{Schwinger2}{Schwinger,J. \jmp{5}{1964}{1606}.}
  \ref{Muller}{M\"uller, \lnm{}{}{}.}
  \ref{VMK}{Varshalovich.}
  \ref{DandWo}{Dowker,J.S. and Wolski, A. \prA{46}{1992}{6417}.}
  \ref{Zeitlin1}{Zeitlin,V. {\it Physica D} {\bf 49} (1991).  }
  \ref{Zeitlin0}{Zeitlin,V. {\it Nonlinear World} Ed by
   V.Baryakhtar {\it et al},  Vol.I p.717,  (World Scientific, Singapore, 1989).}
  \ref{Zeitlin2}{Zeitlin,V. \prl{93}{2004}{264501}. }
  \ref{Zeitlin3}{Zeitlin,V. \pla{339}{2005}{316}. }
  \ref{Groenewold}{Groenewold, H.J. {\it Physica} {\bf 12} (1946) 405.}
  \ref{Cohen}{Cohen, L. \jmp{7}{1966}{781}.}
  \ref{AandW}{Argawal G.S. and Wolf, E. \prD{2}{1970}{2161,2187,2206}.}
  \ref{Jantzen}{Jantzen,R.T. \jmp{19}{1978}{1163}.}
  \ref{Moses2}{Moses,H.E. \aop{42}{1967}{343}.}
  \ref{Carmeli}{Carmeli,M. \jmp{9}{1968}{1987}.}
  \ref{SHS}{Siemans,M., Hancock,J. and Siminovitch,D. {\it Solid State
  Nuclear Magnetic Resonance} {\bf 31}(2007)35.}
 \ref{Dowk}{Dowker,J.S. \prD{28}{1983}{3013}.}
 \ref{Heine}{Heine, E. {\it Handbuch der Kugelfunctionen}
  (G.Reimer, Berlin. 1878, 1881).}
  \ref{Pockels}{Pockels, F. {\it \"Uber die Differentialgleichung $\De
  u+k^2u=0$} (Teubner, Leipzig. 1891).}
  \ref{Hamermesh}{Hamermesh, M., {\it Group Theory} (Addison--Wesley,
  Reading. 1962).}
  \ref{Racah}{Racah, G. {\it Group Theory and Spectroscopy}
  (Princeton Lecture Notes, 1951). }
  \ref{Gourdin}{Gourdin, M. {\it Basics of Lie Groups} (Editions
  Fronti\'eres, Gif sur Yvette. 1982.)}
  \ref{Clifford}{Clifford, W.K. \plms{2}{1866}{116}.}
  \ref{Story2}{Story, W.E. \plms{23}{1892}{265}.}
  \ref{Story}{Story, W.E. \ma{41}{1893}{469}.}
  \ref{Poole}{Poole, E.G.C. \plms{33}{1932}{435}.}
  \ref{Dickson}{Dickson, L.E. {\it Algebraic Invariants} (Wiley, N.Y.
  1915).}
  \ref{Dickson2}{Dickson, L.E. {\it Modern Algebraic Theories}
  (Sanborn and Co., Boston. 1926).}
  \ref{Hilbert2}{Hilbert, D. {\it Theory of algebraic invariants} (C.U.P.,
  Cambridge. 1993).}
  \ref{Olver}{Olver, P.J. {\it Classical Invariant Theory} (C.U.P., Cambridge.
  1999.)}
  \ref{AST}{A\v{s}erova, R.M., Smirnov, J.F. and Tolsto\v{i}, V.N. {\it
  Teoret. Mat. Fyz.} {\bf 8} (1971) 255.}
  \ref{AandS}{A\v{s}erova, R.M., Smirnov, J.F. \np{4}{1968}{399}.}
  \ref{Shapiro}{Shapiro, J. \jmp{6}{1965}{1680}.}
  \ref{Shapiro2}{Shapiro, J.Y. \jmp{14}{1973}{1262}.}
  \ref{NandS}{Noz, M.E. and Shapiro, J.Y. \np{51}{1973}{309}.}
  \ref{Cayley2}{Cayley, A. {\it Phil. Trans. Roy. Soc. Lond.}
  {\bf 144} (1854) 244.}
  \ref{Cayley3}{Cayley, A. {\it Phil. Trans. Roy. Soc. Lond.}
  {\bf 146} (1856) 101.}
  \ref{Wigner}{Wigner, E.P. {\it Gruppentheorie} (Vieweg, Braunschweig. 1931).}
  \ref{Sharp}{Sharp, R.T. \ajop{28}{1960}{116}.}
  \ref{Laporte}{Laporte, O. {\it Z. f. Naturf.} {\bf 3a} (1948) 447.}
  \ref{Lowdin}{L\"owdin, P-O. \rmp{36}{1964}{966}.}
  \ref{Ansari}{Ansari, S.M.R. {\it Fort. d. Phys.} {\bf 15} (1967) 707.}
  \ref{SSJR}{Samal, P.K., Saha, R., Jain, P. and Ralston, J.P. {\it
  Testing Isotropy of Cosmic Microwave Background Radiation},
  astro-ph/0708.2816.}
  \ref{Lachieze}{Lachi\'eze-Rey, M. {\it Harmonic projection and
  multipole Vectors}. astro- \break ph/0409081.}
  \ref{CHS}{Copi, C.J., Huterer, D. and Starkman, G.D.
  \prD{70}{2003}{043515}.}
  \ref{Jaric}{Jari\'c, J.P. {\it Int. J. Eng. Sci.} {\bf 41} (2003) 2123.}
  \ref{RandD}{Roche, J.A. and Dowker, J.S. \jpa{1}{1968}{527}.}
  \ref{KandW}{Katz, G. and Weeks, J.R. \prD{70}{2004}{063527}.}
  \ref{Waerden}{van der Waerden, B.L. {\it Die Gruppen-theoretische
  Methode in der Quantenmechanik} (Springer, Berlin. 1932).}
  \ref{EMOT}{Erdelyi, A., Magnus, W., Oberhettinger, F. and Tricomi, F.G. {
  \it Higher Transcendental Functions} Vol.1 (McGraw-Hill, N.Y. 1953).}
  \ref{Dowzilch}{Dowker, J.S. {\it Proc. Phys. Soc.} {\bf 91} (1967) 28.}
  \ref{DandD}{Dowker, J.S. and Dowker, Y.P. {\it Proc. Phys. Soc.}
  {\bf 87} (1966) 65.}
  \ref{DandD2}{Dowker, J.S. and Dowker, Y.P. \prs{}{}{}.}
  \ref{Dowk3}{Dowker,J.S. \cqg{7}{1990}{1241}.}
  \ref{Dowk5}{Dowker,J.S. \cqg{7}{1990}{2353}.}
  \ref{CoandH}{Courant, R. and Hilbert, D. {\it Methoden der
  Mathematischen Physik} vol.1 \break (Springer, Berlin. 1931).}
  \ref{Applequist}{Applequist, J. \jpa{22}{1989}{4303}.}
  \ref{Torruella}{Torruella, \jmp{16}{1975}{1637}.}
  \ref{Weinberg}{Weinberg, S.W. \pr{133}{1964}{B1318}.}
  \ref{Meyerw}{Meyer, W.F. {\it Apolarit\"at und rationale Curven}
  (Fues, T\"ubingen. 1883.) }
  \ref{Ostrowski}{Ostrowski, A. {\it Jahrsb. Deutsch. Math. Verein.} {\bf
  33} (1923) 245.}
  \ref{Kramers}{Kramers, H.A. {\it Grundlagen der Quantenmechanik}, (Akad.
  Verlag., Leipzig, 1938).}
  \ref{ZandZ}{Zou, W.-N. and Zheng, Q.-S. \prs{459}{2003}{527}.}
  \ref{Weeks1}{Weeks, J.R. {\it Maxwell's multipole vectors
  and the CMB}.  astro-ph/0412231.}
  \ref{Corson}{Corson, E.M. {\it Tensors, Spinors and Relativistic Wave
  Equations} (Blackie, London. 1950).}
  \ref{Rosanes}{Rosanes, J. \jram{76}{1873}{312}.}
  \ref{Salmon}{Salmon, G. {\it Lessons Introductory to the Modern Higher
  Algebra} 3rd. edn. \break (Hodges,  Dublin. 1876.)}
  \ref{Milnew}{Milne, W.P. {\it Homogeneous Coordinates} (Arnold. London. 1910).}
  \ref{Niven}{Niven, W.D. {\it Phil. Trans. Roy. Soc.} {\bf 170} (1879) 393.}
  \ref{Scott}{Scott, C.A. {\it An Introductory Account of
  Certain Modern Ideas and Methods in Plane Analytical Geometry,}
  (MacMillan, N.Y. 1896).}
  \ref{Bargmann}{Bargmann, V. \rmp{34}{1962}{300}.}
  \ref{Maxwell}{Maxwell, J.C. {\it A Treatise on Electricity and
  Magnetism} 2nd. edn. (Clarendon Press, Oxford. 1882).}
  \ref{BandL}{Biedenharn, L.C. and Louck, J.D.
  {\it Angular Momentum in Quantum Physics} (Addison-Wesley, Reading. 1981).}
  \ref{Weylqm}{Weyl, H. {\it The Theory of Groups and Quantum Mechanics}
  (Methuen, London. 1931).}
  \ref{Robson}{Robson, A. {\it An Introduction to Analytical Geometry} Vol I
  (C.U.P., Cambridge. 1940.)}
  \ref{Sommerville}{Sommerville, D.M.Y. {\it Analytical Conics} 3rd. edn.
   (Bell, London. 1933).}
  \ref{Coolidge}{Coolidge, J.L. {\it A Treatise on Algebraic Plane Curves}
  (Clarendon Press, Oxford. 1931).}
  \ref{SandK}{Semple, G. and Kneebone. G.T. {\it Algebraic Projective
  Geometry} (Clarendon Press, Oxford. 1952).}
  \ref{AandC}{Abdesselam A., and Chipalkatti, J. {\it The Higher
  Transvectants are redundant}, arXiv:0801.1533 [math.AG] 2008.}
  \ref{Elliott}{Elliott, E.B. {\it The Algebra of Quantics} 2nd edn.
  (Clarendon Press, Oxford. 1913).}
  \ref{Elliott2}{Elliott, E.B. \qjpam{48}{1917}{372}.}
  \ref{Howe}{Howe, R. \tams{313}{1989}{539}.}
  \ref{Clebsch}{Clebsch, A. \jram{60}{1862}{343}.}
  \ref{Prasad}{Prasad, G. \ma{72}{1912}{136}.}
  \ref{Dougall}{Dougall, J. \pems{32}{1913}{30}.}
  \ref{Penrose}{Penrose, R. \aop{10}{1960}{171}.}
  \ref{Penrose2}{Penrose, R. \prs{273}{1965}{171}.}
  \ref{Burnside}{Burnside, W.S. \qjm{10}{1870}{211}. }
  \ref{Lindemann}{Lindemann, F. \ma{23} {1884}{111}.}
  \ref{Backus}{Backus, G. {\it Rev. Geophys. Space Phys.} {\bf 8} (1970) 633.}
  \ref{Baerheim}{Baerheim, R. {\it Q.J. Mech. appl. Math.} {\bf 51} (1998) 73.}
  \ref{Lense}{Lense, J. {\it Kugelfunktionen} (Akad.Verlag, Leipzig. 1950).}
  \ref{Littlewood}{Littlewood, D.E. \plms{50}{1948}{349}.}
  \ref{Fierz}{Fierz, M. {\it Helv. Phys. Acta} {\bf 12} (1938) 3.}
  \ref{Williams}{Williams, D.N. {\it Lectures in Theoretical Physics} Vol. VII,
  (Univ.Colorado Press, Boulder. 1965).}
  \ref{Dennis}{Dennis, M. \jpa{37}{2004}{9487}.}
  \ref{Pirani}{Pirani, F. {\it Brandeis Lecture Notes on
  General Relativity,} edited by S. Deser and K. Ford. (Brandeis, Mass. 1964).}
  \ref{Sturm}{Sturm, R. \jram{86}{1878}{116}.}
  \ref{Schlesinger}{Schlesinger, O. \ma{22}{1883}{521}.}
  \ref{Askwith}{Askwith, E.H. {\it Analytical Geometry of the Conic
  Sections} (A.\&C. Black, London. 1908).}
  \ref{Todd}{Todd, J.A. {\it Projective and Analytical Geometry}.
  (Pitman, London. 1946).}
  \ref{Glenn}{Glenn. O.E. {\it Theory of Invariants} (Ginn \& Co, N.Y. 1915).}
  \ref{DowkandG}{Dowker, J.S. and Goldstone, M. \prs{303}{1968}{381}.}
  \ref{Turnbull}{Turnbull, H.A. {\it The Theory of Determinants,
  Matrices and Invariants} 3rd. edn. (Dover, N.Y. 1960).}
  \ref{MacMillan}{MacMillan, W.D. {\it The Theory of the Potential}
  (McGraw-Hill, N.Y. 1930).}
   \ref{Hobson}{Hobson, E.W. {\it The Theory of Spherical
   and Ellipsoidal Harmonics} (C.U.P., Cambridge. 1931).}
  \ref{Hobson1}{Hobson, E.W. \plms {24}{1892}{55}.}
  \ref{GandY}{Grace, J.H. and Young, A. {\it The Algebra of Invariants}
  (C.U.P., Cambridge, 1903).}
  \ref{FandR}{Fano, U. and Racah, G. {\it Irreducible Tensorial Sets}
  (Academic Press, N.Y. 1959).}
  \ref{TandT}{Thomson, W. and Tait, P.G. {\it Treatise on Natural Philosophy}
   (Clarendon Press, Oxford. 1867).}
  \ref{Brinkman}{Brinkman, H.C. {\it Applications of spinor invariants in
atomic physics}, North Holland, Amsterdam 1956.}
  \ref{Kramers1}{Kramers, H.A. {\it Proc. Roy. Soc. Amst.} {\bf 33} (1930) 953.}
  \ref{DandP2}{Dowker,J.S. and Pettengill,D.F. \jpa{7}{1974}{1527}}
  \ref{Dowk1}{Dowker,J.S. \jpa{}{}{45}.}
  \ref{Dowk2}{Dowker,J.S. \aop{71}{1972}{577}}
  \ref{DandA}{Dowker,J.S. and Apps, J.S. \cqg{15}{1998}{1121}.}
  \ref{Weil}{Weil,A., {\it Elliptic functions according to Eisenstein
  and Kronecker}, Springer, Berlin, 1976.}
  \ref{Ling}{Ling,C-H. {\it SIAM J.Math.Anal.} {\bf5} (1974) 551.}
  \ref{Ling2}{Ling,C-H. {\it J.Math.Anal.Appl.}(1988).}
 \ref{BMO}{Brevik,I., Milton,K.A. and Odintsov, S.D. \aop{302}{2002}{120}.}
 \ref{KandL}{Kutasov,D. and Larsen,F. {\it JHEP} 0101 (2001) 1.}
 \ref{KPS}{Klemm,D., Petkou,A.C. and Siopsis {\it Entropy
 bounds, monoticity properties and scaling in CFT's}. hep-th/0101076.}
 \ref{DandC}{Dowker,J.S. and Critchley,R. \prD{15}{1976}{1484}.}
 \ref{AandD}{Al'taie, M.B. and Dowker, J.S. \prD{18}{1978}{3557}.}
 \ref{Dow1}{Dowker,J.S. \prD{37}{1988}{558}.}
 \ref{Dow30}{Dowker,J.S. \prD{28}{1983}{3013}.}
 \ref{DandK}{Dowker,J.S. and Kennedy,G. \jpa{}{1978}{895}.}
 \ref{Dow2}{Dowker,J.S. \cqg{1}{1984}{359}.}
 \ref{DandKi}{Dowker,J.S. and Kirsten, K. {\it Comm. in Anal. and Geom.
 }{\bf7} (1999) 641.}
 \ref{DandKe}{Dowker,J.S. and Kennedy,G.\jpa{11}{1978}{895}.}
 \ref{Gibbons}{Gibbons,G.W. \pl{60A}{1977}{385}.}
 \ref{Cardy}{Cardy,J.L. \np{366}{1991}{403}.}
 \ref{ChandD}{Chang,P. and Dowker,J.S. \np{395}{1993}{407}.}
 \ref{DandC2}{Dowker,J.S. and Critchley,R. \prD{13}{1976}{224}.}
 \ref{Camporesi}{Camporesi,R. \prp{196}{1990}{1}.}
 \ref{BandM}{Brown,L.S. and Maclay,G.J. \pr{184}{1969}{1272}.}
 \ref{CandD}{Candelas,P. and Dowker,J.S. \prD{19}{1979}{2902}.}
 \ref{Unwin1}{Unwin,S.D. Thesis. University of Manchester. 1979.}
 \ref{Unwin2}{Unwin,S.D. \jpa{13}{1980}{313}.}
 \ref{DandB}{Dowker,J.S. and Banach,R. \jpa{11}{1978}{2255}.}
 \ref{Obhukov}{Obhukov,Yu.N. \pl{109B}{1982}{195}.}
 \ref{Kennedy}{Kennedy,G. \prD{23}{1981}{2884}.}
 \ref{CandT}{Copeland,E. and Toms,D.J. \np {255}{1985}{201}.}
  \ref{CandT2}{Copeland,E. and Toms,D.J. \cqg {3}{1986}{431}.}
 \ref{ELV}{Elizalde,E., Lygren, M. and Vassilevich,
 D.V. \jmp {37}{1996}{3105}.}
 \ref{Malurkar}{Malurkar,S.L. {\it J.Ind.Math.Soc} {\bf16} (1925/26) 130.}
 \ref{Glaisher}{Glaisher,J.W.L. {\it Messenger of Math.} {\bf18} (1889) 1.}
  \ref{Anderson}{Anderson,A. \prD{37}{1988}{536}.}
 \ref{CandA}{Cappelli,A. and D'Appollonio,\pl{487B}{2000}{87}.}
 \ref{Wot}{Wotzasek,C. \jpa{23}{1990}{1627}.}
 \ref{RandT}{Ravndal,F. and Tollesen,D. \prD{40}{1989}{4191}.}
 \ref{SandT}{Santos,F.C. and Tort,A.C. \pl{482B}{2000}{323}.}
 \ref{FandO}{Fukushima,K. and Ohta,K. {\it Physica} {\bf A299} (2001) 455.}
 \ref{GandP}{Gibbons,G.W. and Perry,M. \prs{358}{1978}{467}.}
 \ref{Dow4}{Dowker,J.S..}
  \ref{Rad}{Rademacher,H. {\it Topics in analytic number theory,}
Springer-Verlag,  Berlin,1973.}
  \ref{Halphen}{Halphen,G.-H. {\it Trait\'e des Fonctions Elliptiques},
  Vol 1, Gauthier-Villars, Paris, 1886.}
  \ref{CandW}{Cahn,R.S. and Wolf,J.A. {\it Comm.Mat.Helv.} {\bf 51}
  (1976) 1.}
  \ref{Berndt}{Berndt,B.C. \rmjm{7}{1977}{147}.}
  \ref{Hurwitz}{Hurwitz,A. \ma{18}{1881}{528}.}
  \ref{Hurwitz2}{Hurwitz,A. {\it Mathematische Werke} Vol.I. Basel,
  Birkhauser, 1932.}
  \ref{Berndt2}{Berndt,B.C. \jram{303/304}{1978}{332}.}
  \ref{RandA}{Rao,M.B. and Ayyar,M.V. \jims{15}{1923/24}{150}.}
  \ref{Hardy}{Hardy,G.H. \jlms{3}{1928}{238}.}
  \ref{TandM}{Tannery,J. and Molk,J. {\it Fonctions Elliptiques},
   Gauthier-Villars, Paris, 1893--1902.}
  \ref{schwarz}{Schwarz,H.-A. {\it Formeln und
  Lehrs\"atzen zum Gebrauche..},Springer 1893.(The first edition was 1885.)
  The French translation by Henri Pad\'e is {\it Formules et Propositions
  pour L'Emploi...},Gauthier-Villars, Paris, 1894}
  \ref{Hancock}{Hancock,H. {\it Theory of elliptic functions}, Vol I.
   Wiley, New York 1910.}
  \ref{watson}{Watson,G.N. \jlms{3}{1928}{216}.}
  \ref{MandO}{Magnus,W. and Oberhettinger,F. {\it Formeln und S\"atze},
  Springer-Verlag, Berlin 1948.}
  \ref{Klein}{Klein,F. {\it Lectures on the Icosohedron}
  (Methuen, London. 1913).}
  \ref{AandL}{Appell,P. and Lacour,E. {\it Fonctions Elliptiques},
  Gauthier-Villars,
  Paris. 1897.}
  \ref{HandC}{Hurwitz,A. and Courant,C. {\it Allgemeine Funktionentheorie},
  Springer,
  Berlin. 1922.}
  \ref{WandW}{Whittaker,E.T. and Watson,G.N. {\it Modern analysis},
  Cambridge. 1927.}
  \ref{SandC}{Selberg,A. and Chowla,S. \jram{227}{1967}{86}. }
  \ref{zucker}{Zucker,I.J. {\it Math.Proc.Camb.Phil.Soc} {\bf 82 }(1977)
  111.}
  \ref{glasser}{Glasser,M.L. {\it Maths.of Comp.} {\bf 25} (1971) 533.}
  \ref{GandW}{Glasser, M.L. and Wood,V.E. {\it Maths of Comp.} {\bf 25}
  (1971)
  535.}
  \ref{greenhill}{Greenhill,A,G. {\it The Applications of Elliptic
  Functions}, MacMillan. London, 1892.}
  \ref{Weierstrass}{Weierstrass,K. {\it J.f.Mathematik (Crelle)}
{\bf 52} (1856) 346.}
  \ref{Weierstrass2}{Weierstrass,K. {\it Mathematische Werke} Vol.I,p.1,
  Mayer u. M\"uller, Berlin, 1894.}
  \ref{Fricke}{Fricke,R. {\it Die Elliptische Funktionen und Ihre Anwendungen},
    Teubner, Leipzig. 1915, 1922.}
  \ref{Konig}{K\"onigsberger,L. {\it Vorlesungen \"uber die Theorie der
 Elliptischen Funktionen},  \break Teubner, Leipzig, 1874.}
  \ref{Milne}{Milne,S.C. {\it The Ramanujan Journal} {\bf 6} (2002) 7-149.}
  \ref{Schlomilch}{Schl\"omilch,O. {\it Ber. Verh. K. Sachs. Gesell. Wiss.
  Leipzig}  {\bf 29} (1877) 101-105; {\it Compendium der h\"oheren
  Analysis}, Bd.II, 3rd Edn, Vieweg, Brunswick, 1878.}
  \ref{BandB}{Briot,C. and Bouquet,C. {\it Th\`eorie des Fonctions
  Elliptiques}, Gauthier-Villars, Paris, 1875.}
  \ref{Dumont}{Dumont,D. \aim {41}{1981}{1}.}
  \ref{Andre}{Andr\'e,D. {\it Ann.\'Ecole Normale Superior} {\bf 6} (1877)
  265;
  {\it J.Math.Pures et Appl.} {\bf 5} (1878) 31.}
  \ref{Raman}{Ramanujan,S. {\it Trans.Camb.Phil.Soc.} {\bf 22} (1916) 159;
 {\it Collected Papers}, Cambridge, 1927}
  \ref{Weber}{Weber,H.M. {\it Lehrbuch der Algebra} Bd.III, Vieweg,
  Brunswick 190  3.}
  \ref{Weber2}{Weber,H.M. {\it Elliptische Funktionen und algebraische
  Zahlen},
  Vieweg, Brunswick 1891.}
  \ref{ZandR}{Zucker,I.J. and Robertson,M.M.
  {\it Math.Proc.Camb.Phil.Soc} {\bf 95 }(1984) 5.}
  \ref{JandZ1}{Joyce,G.S. and Zucker,I.J.
  {\it Math.Proc.Camb.Phil.Soc} {\bf 109 }(1991) 257.}
  \ref{JandZ2}{Zucker,I.J. and Joyce.G.S.
  {\it Math.Proc.Camb.Phil.Soc} {\bf 131 }(2001) 309.}
  \ref{zucker2}{Zucker,I.J. {\it SIAM J.Math.Anal.} {\bf 10} (1979) 192.}
  \ref{BandZ}{Borwein,J.M. and Zucker,I.J. {\it IMA J.Math.Anal.} {\bf 12}
  (1992) 519.}
  \ref{Cox}{Cox,D.A. {\it Primes of the form $x^2+n\,y^2$}, Wiley,
  New York, 1989.}
  \ref{BandCh}{Berndt,B.C. and Chan,H.H. {\it Mathematika} {\bf42} (1995)
  278.}
  \ref{EandT}{Elizalde,R. and Tort.hep-th/}
  \ref{KandS}{Kiyek,K. and Schmidt,H. {\it Arch.Math.} {\bf 18} (1967) 438.}
  \ref{Oshima}{Oshima,K. \prD{46}{1992}{4765}.}
  \ref{greenhill2}{Greenhill,A.G. \plms{19} {1888} {301}.}
  \ref{Russell}{Russell,R. \plms{19} {1888} {91}.}
  \ref{BandB}{Borwein,J.M. and Borwein,P.B. {\it Pi and the AGM}, Wiley,
  New York, 1998.}
  \ref{Resnikoff}{Resnikoff,H.L. \tams{124}{1966}{334}.}
  \ref{vandp}{Van der Pol, B. {\it Indag.Math.} {\bf18} (1951) 261,272.}
  \ref{Rankin}{Rankin,R.A. {\it Modular forms} C.U.P. Cambridge}
  \ref{Rankin2}{Rankin,R.A. {\it Proc. Roy.Soc. Edin.} {\bf76 A} (1976) 107.}
  \ref{Skoruppa}{Skoruppa,N-P. {\it J.of Number Th.} {\bf43} (1993) 68 .}
  \ref{Down}{Dowker.J.S. {\it Nucl.Phys.}B (Proc.Suppl) ({\bf 104})(2002)153;
  also Dowker,J.S. hep-th/ 0007129.}
  \ref{Eichler}{Eichler,M. \mz {67}{1957}{267}.}
  \ref{Zagier}{Zagier,D. \invm{104}{1991}{449}.}
  \ref{Lang}{Lang,S. {\it Modular Forms}, Springer, Berlin, 1976.}
  \ref{Kosh}{Koshliakov,N.S. {\it Mess.of Math.} {\bf 58} (1928) 1.}
  \ref{BandH}{Bodendiek, R. and Halbritter,U. \amsh{38}{1972}{147}.}
  \ref{Smart}{Smart,L.R., \pgma{14}{1973}{1}.}
  \ref{Grosswald}{Grosswald,E. {\it Acta. Arith.} {\bf 21} (1972) 25.}
  \ref{Kata}{Katayama,K. {\it Acta Arith.} {\bf 22} (1973) 149.}
  \ref{Ogg}{Ogg,A. {\it Modular forms and Dirichlet series} (Benjamin,
  New York,
   1969).}
  \ref{Bol}{Bol,G. \amsh{16}{1949}{1}.}
  \ref{Epstein}{Epstein,P. \ma{56}{1903}{615}.}
  \ref{Petersson}{Petersson.}
  \ref{Serre}{Serre,J-P. {\it A Course in Arithmetic}, Springer,
  New York, 1973.}
  \ref{Schoenberg}{Schoenberg,B., {\it Elliptic Modular Functions},
  Springer, Berlin, 1974.}
  \ref{Apostol}{Apostol,T.M. \dmj {17}{1950}{147}.}
  \ref{Ogg2}{Ogg,A. {\it Lecture Notes in Math.} {\bf 320} (1973) 1.}
  \ref{Knopp}{Knopp,M.I. \dmj {45}{1978}{47}.}
  \ref{Knopp2}{Knopp,M.I. \invm {}{1994}{361}.}
  \ref{LandZ}{Lewis,J. and Zagier,D. \aom{153}{2001}{191}.}
  \ref{DandK1}{Dowker,J.S. and Kirsten,K. {\it Elliptic functions and
  temperature inversion symmetry on spheres} hep-th/.}
  \ref{HandK}{Husseini and Knopp.}
  \ref{Kober}{Kober,H. \mz{39}{1934-5}{609}.}
  \ref{HandL}{Hardy,G.H. and Littlewood, \am{41}{1917}{119}.}
  \ref{Watson}{Watson,G.N. \qjm{2}{1931}{300}.}
  \ref{SandC2}{Chowla,S. and Selberg,A. {\it Proc.Nat.Acad.} {\bf 35}
  (1949) 371.}
  \ref{Landau}{Landau, E. {\it Lehre von der Verteilung der Primzahlen},
  (Teubner, Leipzig, 1909).}
  \ref{Berndt4}{Berndt,B.C. \tams {146}{1969}{323}.}
  \ref{Berndt3}{Berndt,B.C. \tams {}{}{}.}
  \ref{Bochner}{Bochner,S. \aom{53}{1951}{332}.}
  \ref{Weil2}{Weil,A.\ma{168}{1967}{}.}
  \ref{CandN}{Chandrasekharan,K. and Narasimhan,R. \aom{74}{1961}{1}.}
  \ref{Rankin3}{Rankin,R.A. {} {} ().}
  \ref{Berndt6}{Berndt,B.C. {\it Trans.Edin.Math.Soc}.}
  \ref{Elizalde}{Elizalde,E. {\it Ten Physical Applications of Spectral
  Zeta Function Theory}, \break (Springer, Berlin, 1995).}
  \ref{Allen}{Allen,B., Folacci,A. and Gibbons,G.W. \pl{189}{1987}{304}.}
  \ref{Krazer}{Krazer}
  \ref{Elizalde3}{Elizalde,E. {\it J.Comp.and Appl. Math.} {\bf 118}
  (2000) 125.}
  \ref{Elizalde2}{Elizalde,E., Odintsov.S.D, Romeo, A. and Bytsenko,
  A.A and
  Zerbini,S.
  {\it Zeta function regularisation}, (World Scientific, Singapore,
  1994).}
  \ref{Eisenstein}{Eisenstein}
  \ref{Hecke}{Hecke,E. \ma{112}{1936}{664}.}
  \ref{Hecke2}{Hecke,E. \ma{112}{1918}{398}.}
  \ref{Terras}{Terras,A. {\it Harmonic analysis on Symmetric Spaces} (Springer,
  New York, 1985).}
  \ref{BandG}{Bateman,P.T. and Grosswald,E. {\it Acta Arith.} {\bf 9}
  (1964) 365.}
  \ref{Deuring}{Deuring,M. \aom{38}{1937}{585}.}
  \ref{Mordell}{Mordell,J. \prs{}{}{}.}
  \ref{GandZ}{Glasser,M.L. and Zucker, {}.}
  \ref{Landau2}{Landau,E. \jram{}{1903}{64}.}
  \ref{Kirsten1}{Kirsten,K. \jmp{35}{1994}{459}.}
  \ref{Sommer}{Sommer,J. {\it Vorlesungen \"uber Zahlentheorie}
  (1907,Teubner,Leipzig).
  French edition 1913 .}
  \ref{Reid}{Reid,L.W. {\it Theory of Algebraic Numbers},
  (1910,MacMillan,New York).}
  \ref{Milnor}{Milnor, J. {\it Is the Universe simply--connected?},
  IAS, Princeton, 1978.}
  \ref{Milnor2}{Milnor, J. \ajm{79}{1957}{623}.}
  \ref{Opechowski}{Opechowski,W. {\it Physica} {\bf 7} (1940) 552.}
  \ref{Bethe}{Bethe, H.A. \zfp{3}{1929}{133}.}
  \ref{LandL}{Landau, L.D. and Lishitz, E.M. {\it Quantum
  Mechanics} (Pergamon Press, London, 1958).}
  \ref{GPR}{Gibbons, G.W., Pope, C. and R\"omer, H., \np{157}{1979}{377}.}
  \ref{Jadhav}{Jadhav,S.P. PhD Thesis, University of Manchester 1990.}
  \ref{DandJ}{Dowker,J.S. and Jadhav, S. \prD{39}{1989}{1196}.}
  \ref{CandM}{Coxeter, H.S.M. and Moser, W.O.J. {\it Generators and
  relations of finite groups} (Springer. Berlin. 1957).}
  \ref{Coxeter2}{Coxeter, H.S.M. {\it Regular Complex Polytopes},
   (Cambridge University Press, \break Cambridge, 1975).}
  \ref{Coxeter}{Coxeter, H.S.M. {\it Regular Polytopes}.}
  \ref{Stiefel}{Stiefel, E., J.Research NBS {\bf 48} (1952) 424.}
  \ref{BandS}{Brink, D.M. and Satchler, G.R. {\it Angular momentum theory}.
  (Clarendon Press, Oxford. 1962.).}
  \ref{Rose}{Rose}
  \ref{Schwinger}{Schwinger, J. {\it On Angular Momentum}
  in {\it Quantum Theory of Angular Momentum} edited by
  Biedenharn,L.C. and van Dam, H. (Academic Press, N.Y. 1965).}
  
  \ref{Ray}{Ray,D.B. \aim{4}{1970}{109}.}
  \ref{Ikeda}{Ikeda,A. {\it Kodai Math.J.} {\bf 18} (1995) 57.}
  \ref{Kennedy}{Kennedy,G. \prD{23}{1981}{2884}.}
  \ref{Ellis}{Ellis,G.F.R. {\it General Relativity} {\bf2} (1971) 7.}
  \ref{Dow8}{Dowker,J.S. \cqg{20}{2003}{L105}.}
  \ref{IandY}{Ikeda, A and Yamamoto, Y. \ojm {16}{1979}{447}.}
  \ref{BandI}{Bander,M. and Itzykson,C. \rmp{18}{1966}{2}.}
  \ref{Schulman}{Schulman, L.S. \pr{176}{1968}{1558}.}
  \ref{Bar1}{B\"ar,C. {\it Arch.d.Math.}{\bf 59} (1992) 65.}
  \ref{Bar2}{B\"ar,C. {\it Geom. and Func. Anal.} {\bf 6} (1996) 899.}
  \ref{Vilenkin}{Vilenkin, N.J. {\it Special functions},
  (Am.Math.Soc., Providence, 1968).}
  \ref{Talman}{Talman, J.D. {\it Special functions} (Benjamin,N.Y.,1968).}
  \ref{Miller}{Miller, W. {\it Symmetry groups and their applications}
  (Wiley, N.Y., 1972).}
  \ref{Dow3}{Dowker,J.S. \cmp{162}{1994}{633}.}
  \ref{Cheeger}{Cheeger, J. \jdg {18}{1983}{575}.}
  \ref{Cheeger2}{Cheeger, J. \aom {109}{1979}{259}.}
  \ref{Dow6}{Dowker,J.S. \jmp{30}{1989}{770}.}
  \ref{Dow20}{Dowker,J.S. \jmp{35}{1994}{6076}.}
  \ref{Dowjmp}{Dowker,J.S. \jmp{35}{1994}{4989}.}
  \ref{Dow21}{Dowker,J.S. {\it Heat kernels and polytopes} in {\it
   Heat Kernel Techniques and Quantum Gravity}, ed. by S.A.Fulling,
   Discourses in Mathematics and its Applications, No.4, Dept.
   Maths., Texas A\&M University, College Station, Texas, 1995.}
  \ref{Dow9}{Dowker,J.S. \jmp{42}{2001}{1501}.}
  \ref{Dow7}{Dowker,J.S. \jpa{25}{1992}{2641}.}
  \ref{Warner}{Warner.N.P. \prs{383}{1982}{379}.}
  \ref{Wolf}{Wolf, J.A. {\it Spaces of constant curvature},
  (McGraw--Hill,N.Y., 1967).}
  \ref{Meyer}{Meyer,B. \cjm{6}{1954}{135}.}
  \ref{BandB}{B\'erard,P. and Besson,G. {\it Ann. Inst. Four.} {\bf 30}
  (1980) 237.}
  \ref{PandM}{P\'{o}lya,G. and Meyer,B. \cras{228}{1948}{28}.}
  \ref{Springer}{Springer, T.A. Lecture Notes in Math. vol 585 (Springer,
  Berlin,1977).}
  \ref{SeandT}{Threlfall, H. and Seifert, W. \ma{104}{1930}{1}.}
  \ref{Hopf}{Hopf,H. \ma{95}{1925}{313}. }
  \ref{Dow}{Dowker,J.S. \jpa{5}{1972}{936}.}
  \ref{LLL}{Lehoucq,R., Lachi\'eze-Rey,M. and Luminet, J.--P. {\it
  Astron.Astrophys.} {\bf 313} (1996) 339.}
  \ref{LaandL}{Lachi\'eze-Rey,M. and Luminet, J.--P.
  \prp{254}{1995}{135}.}
  \ref{Schwarzschild}{Schwarzschild, K., {\it Vierteljahrschrift der
  Ast.Ges.} {\bf 35} (1900) 337.}
  \ref{Starkman}{Starkman,G.D. \cqg{15}{1998}{2529}.}
  \ref{LWUGL}{Lehoucq,R., Weeks,J.R., Uzan,J.P., Gausman, E. and
  Luminet, J.--P. \cqg{19}{2002}{4683}.}
  \ref{Dow10}{Dowker,J.S. \prD{28}{1983}{3013}.}
  \ref{BandD}{Banach, R. and Dowker, J.S. \jpa{12}{1979}{2527}.}
  \ref{Jadhav2}{Jadhav,S. \prD{43}{1991}{2656}.}
  \ref{Gilkey}{Gilkey,P.B. {\it Invariance theory,the heat equation and
  the Atiyah--Singer Index theorem} (CRC Press, Boca Raton, 1994).}
  \ref{BandY}{Berndt,B.C. and Yeap,B.P. {\it Adv. Appl. Math.}
  {\bf29} (2002) 358.}
  \ref{HandR}{Hanson,A.J. and R\"omer,H. \pl{80B}{1978}{58}.}
  \ref{Hill}{Hill,M.J.M. {\it Trans.Camb.Phil.Soc.} {\bf 13} (1883) 36.}
  \ref{Cayley}{Cayley,A. {\it Quart.Math.J.} {\bf 7} (1866) 304.}
  \ref{Seade}{Seade,J.A. {\it Anal.Inst.Mat.Univ.Nac.Aut\'on
  M\'exico} {\bf 21} (1981) 129.}
  \ref{CM}{Cisneros--Molina,J.L. {\it Geom.Dedicata} {\bf84} (2001)
  \ref{Goette1}{Goette,S. \jram {526} {2000} 181.}
  207.}
  \ref{NandO}{Nash,C. and O'Connor,D--J, \jmp {36}{1995}{1462}.}
  \ref{Dows}{Dowker,J.S. \aop{71}{1972}{577}; Dowker,J.S. and Pettengill,D.F.
  \jpa{7}{1974}{1527}; J.S.Dowker in {\it Quantum Gravity}, edited by
  S. C. Christensen (Hilger,Bristol,1984)}
  \ref{Jadhav2}{Jadhav,S.P. \prD{43}{1991}{2656}.}
  \ref{Dow11}{Dowker,J.S. \cqg{21}{2004}4247.}
  \ref{Dow12}{Dowker,J.S. \cqg{21}{2004}4977.}
  \ref{Dow13}{Dowker,J.S. \jpa{38}{2005}1049.}
  \ref{Zagier}{Zagier,D. \ma{202}{1973}{149}}
  \ref{RandG}{Rademacher, H. and Grosswald,E. {\it Dedekind Sums},
  (Carus, MAA, 1972).}
  \ref{Berndt7}{Berndt,B, \aim{23}{1977}{285}.}
  \ref{HKMM}{Harvey,J.A., Kutasov,D., Martinec,E.J. and Moore,G.
  {\it Localised Tachyons and RG Flows}, hep-th/0111154.}
  \ref{Beck}{Beck,M., {\it Dedekind Cotangent Sums}, {\it Acta Arithmetica}
  {\bf 109} (2003) 109-139 ; math.NT/0112077.}
  \ref{McInnes}{McInnes,B. {\it APS instability and the topology of the brane
  world}, hep-th/0401035.}
  \ref{BHS}{Brevik,I, Herikstad,R. and Skriudalen,S. {\it Entropy Bound for the
  TM Electromagnetic Field in the Half Einstein Universe}; hep-th/0508123.}
  \ref{BandO}{Brevik,I. and Owe,C.  \prD{55}{4689}{1997}.}
  \ref{Kenn}{Kennedy,G. Thesis. University of Manchester 1978.}
  \ref{KandU}{Kennedy,G. and Unwin S. \jpa{12}{L253}{1980}.}
  \ref{BandO1}{Bayin,S.S.and Ozcan,M.
  \prD{48}{2806}{1993}; \prD{49}{5313}{1994}.}
  \ref{Chang}{Chang, P., {\it Quantum Field Theory on Regular Polytopes}.
   Thesis. University of Manchester, 1993.}
  \ref{Barnesa}{Barnes,E.W. {\it Trans. Camb. Phil. Soc.} {\bf 19} (1903) 374.}
  \ref{Barnesb}{Barnes,E.W. {\it Trans. Camb. Phil. Soc.}
  {\bf 19} (1903) 426.}
  \ref{Stanley1}{Stanley,R.P. \joa {49Hilf}{1977}{134}.}
  \ref{Stanley2}{Stanley,R.P. {\it Enumerative Combinatorics} vols.1,2
  (C.U.P., Cambridge, 1997, 1999.}
  \ref{Stanley}{Stanley,R.P. \bams {1}{1979}{475}.}
  \ref{Hurley}{Hurley,A.C. \pcps {47}{1951}{51}.}
  \ref{IandK}{Iwasaki,I. and Katase,K. {\it Proc.Japan Acad. Ser} {\bf A55}
  (1979) 141.}
  \ref{IandT}{Ikeda,A. and Taniguchi,Y. {\it Osaka J. Math.} {\bf 15} (1978)
  515.}
  \ref{GandM}{Gallot,S. and Meyer,D. \jmpa{54}{1975}{259}.}
  \ref{Flatto}{Flatto,L. {\it Enseign. Math.} {\bf 24} (1978) 237.}
  \ref{OandT}{Orlik,P and Terao,H. {\it Arrangements of Hyperplanes},
  Grundlehren der Math. Wiss. {\bf 300}, (Springer--Verlag, 1992).}
  \ref{Shepler}{Shepler,A.V. \joa{220}{1999}{314}.}
  \ref{SandT}{Solomon,L. and Terao,H. \cmh {73}{1998}{237}.}
  \ref{Vass}{Vassilevich, D.V. \plb {348}{1995}39.}
  \ref{Vass2}{Vassilevich, D.V. \jmp {36}{1995}3174.}
  \ref{CandH}{Camporesi,R. and Higuchi,A. {\it J.Geom. and Physics}
  {\bf 15} (1994) 57.}
  \ref{Solomon2}{Solomon,L. \tams{113}{1964}{274}.}
  \ref{Solomon}{Solomon,L. {\it Nagoya Math. J.} {\bf 22} (1963) 57.}
  \ref{Obukhov}{Obukhov,Yu.N. \pl{109B}{1982}{195}.}
  \ref{BGH}{Bernasconi,F., Graf,G.M. and Hasler,D. {\it The heat kernel
  expansion for the electromagnetic field in a cavity}; math-ph/0302035.}
  \ref{Baltes}{Baltes,H.P. \prA {6}{1972}{2252}.}
  \ref{BaandH}{Baltes.H.P and Hilf,E.R. {\it Spectra of Finite Systems}
  (Bibliographisches Institut, Mannheim, 1976).}
  \ref{Ray}{Ray,D.B. \aim{4}{1970}{109}.}
  \ref{Hirzebruch}{Hirzebruch,F. {\it Topological methods in algebraic
  geometry} (Springer-- Verlag,\break  Berlin, 1978). }
  \ref{BBG}{Bla\v{z}i\'c,N., Bokan,N. and Gilkey, P.B. {\it Ind.J.Pure and
  Appl.Math.} {\bf 23} (1992) 103.}
  \ref{WandWi}{Weck,N. and Witsch,K.J. {\it Math.Meth.Appl.Sci.} {\bf 17}
  (1994) 1017.}
  \ref{Norlund}{N\"orlund,N.E. \am{43}{1922}{121}.}
   \ref{Norlund1}{N\"orlund,N.E. {\it Differenzenrechnung} (Springer--Verlag, 1924, Berlin.)}
  \ref{Duff}{Duff,G.F.D. \aom{56}{1952}{115}.}
  \ref{DandS}{Duff,G.F.D. and Spencer,D.C. \aom{45}{1951}{128}.}
  \ref{BGM}{Berger, M., Gauduchon, P. and Mazet, E. {\it Lect.Notes.Math.}
  {\bf 194} (1971) 1. }
  \ref{Patodi}{Patodi,V.K. \jdg{5}{1971}{233}.}
  \ref{GandS}{G\"unther,P. and Schimming,R. \jdg{12}{1977}{599}.}
  \ref{MandS}{McKean,H.P. and Singer,I.M. \jdg{1}{1967}{43}.}
  \ref{Conner}{Conner,P.E. {\it Mem.Am.Math.Soc.} {\bf 20} (1956).}
  \ref{Gilkey2}{Gilkey,P.B. \aim {15}{1975}{334}.}
  \ref{MandP}{Moss,I.G. and Poletti,S.J. \plb{333}{1994}{326}.}
  \ref{BKD}{Bordag,M., Kirsten,K. and Dowker,J.S. \cmp{182}{1996}{371}.}
  \ref{RandO}{Rubin,M.A. and Ordonez,C. \jmp{25}{1984}{2888}.}
  \ref{BaandD}{Balian,R. and Duplantier,B. \aop {112}{1978}{165}.}
  \ref{Kennedy2}{Kennedy,G. \aop{138}{1982}{353}.}
  \ref{DandKi2}{Dowker,J.S. and Kirsten, K. {\it Analysis and Appl.}
 {\bf 3} (2005) 45.}
  \ref{Dow40}{Dowker,J.S. \cqg{23}{2006}{1}.}
  \ref{BandHe}{Br\"uning,J. and Heintze,E. {\it Duke Math.J.} {\bf 51} (1984)
   959.}
  \ref{Dowl}{Dowker,J.S. {\it Functional determinants on M\"obius corners};
    Proceedings, `Quantum field theory under
    the influence of external conditions', 111-121,Leipzig 1995.}
  \ref{Dowqg}{Dowker,J.S. in {\it Quantum Gravity}, edited by
  S. C. Christensen (Hilger, Bristol, 1984).}
  \ref{Dowit}{Dowker,J.S. \jpa{11}{1978}{347}.}
  \ref{Kane}{Kane,R. {\it Reflection Groups and Invariant Theory} (Springer,
  New York, 2001).}
  \ref{Sturmfels}{Sturmfels,B. {\it Algorithms in Invariant Theory}
  (Springer, Vienna, 1993).}
  \ref{Bourbaki}{Bourbaki,N. {\it Groupes et Alg\`ebres de Lie}  Chap.III, IV
  (Hermann, Paris, 1968).}
  \ref{SandTy}{Schwarz,A.S. and Tyupkin, Yu.S. \np{242}{1984}{436}.}
  \ref{Reuter}{Reuter,M. \prD{37}{1988}{1456}.}
  \ref{EGH}{Eguchi,T. Gilkey,P.B. and Hanson,A.J. \prp{66}{1980}{213}.}
  \ref{DandCh}{Dowker,J.S. and Chang,Peter, \prD{46}{1992}{3458}.}
  \ref{APS}{Atiyah M., Patodi and Singer,I.\mpcps{77}{1975}{43}.}
  \ref{Donnelly}{Donnelly.H. {\it Indiana U. Math.J.} {\bf 27} (1978) 889.}
  \ref{Katase}{Katase,K. {\it Proc.Jap.Acad.} {\bf 57} (1981) 233.}
  \ref{Gilkey3}{Gilkey,P.B.\invm{76}{1984}{309}.}
  \ref{Degeratu}{Degeratu.A. {\it Eta--Invariants and Molien Series for
  Unimodular Groups}, Thesis MIT, 2001.}
  \ref{Seeley}{Seeley,R. \ijmp {A\bf18}{2003}{2197}.}
  \ref{Seeley2}{Seeley,R. .}
  \ref{melrose}{Melrose}
  \ref{DandW}{Douglas,R.G. and Wojciekowski,K.P. \cmp{142}{1991}{139}.}
  \ref{Dai}{Dai,X. \tams{354}{2001}{107}.}
\end{putreferences}

\bye